\documentclass[11pt]{article}

\usepackage{epsfig,epsfig}
\usepackage{amsmath}
\usepackage{mathrsfs}
\usepackage{bm}
\usepackage{amssymb}
\usepackage{booktabs}
\usepackage{color}
\usepackage{multirow}
\usepackage{graphicx}
\usepackage{subfigure}
\usepackage{float}
\usepackage{appendix}
\usepackage{tikz}
\usetikzlibrary{calc}
\usetikzlibrary{matrix}

\usepackage{lineno}
\usepackage{srcltx,graphicx}

\newtheorem{thm}{Theorem}[section]
\newtheorem{proposition}[thm]{Proposition}

\newtheorem{lemma}[thm]{Lemma}

\newtheorem{proof}[thm]{Proof}

\newtheorem{example}[thm]{Example}
\newtheorem{rem}[thm]{Remark}

\usepackage{hyperref}

\numberwithin{equation}{section} \topmargin=-2cm \oddsidemargin=0.5cm
\begin{document}
%\linenumbers
%\baselineskip=2pc
\title{\textbf{High-order conservative positivity-preserving DG-interpolation for deforming meshes and application to moving mesh DG simulation of radiative transfer}}

\author{Min Zhang\footnote{School of Mathematical Sciences, Xiamen University,
Xiamen, Fujian 361005, China.
E-mail: minzhang2015@stu.xmu.edu.cn.
},
~Weizhang Huang\footnote{Department of Mathematics, University of Kansas, Lawrence, Kansas 66045, USA. E-mail: whuang@ku.edu.
},
~and Jianxian Qiu\footnote{School of Mathematical Sciences and Fujian Provincial Key Laboratory of Mathematical Modeling and High-Performance Scientific Computing, Xiamen University, Xiamen, Fujian 361005, China.
E-mail: jxqiu@xmu.edu.cn.
}
}

\date{}
\maketitle
\begin{abstract}
Solution interpolation between deforming meshes is an important component for several applications
in scientific computing, including indirect arbitrary-Lagrangian-Eulerian
and rezoning moving mesh methods in numerical solution of partial differential equations.
In this paper, a high-order, conservative, and positivity-preserving interpolation scheme
is developed based on the discontinuous Galerkin solution of a linear time-dependent equation on deforming meshes. The scheme works for bounded but otherwise arbitrary mesh deformation from the old mesh to the new one. The cost and positivity preservation (with a linear scaling limiter) of the DG-interpolation are investigated.
Numerical examples are presented to demonstrate the properties of the interpolation scheme.
The DG-interpolation is applied to the rezoning moving mesh DG solution of
the radiative transfer equation, an integro-differential equation modeling the conservation of photons and involving time, space, and angular variables.
Numerical results obtained for examples in one and two spatial dimensions with various settings show that the resulting rezoning moving mesh DG method maintains the same convergence order as the standard DG method, is more efficient than the method with a fixed uniform mesh, and is able to preserve the positivity of the radiative intensity.
\end{abstract}

\noindent\textbf{The 2010 Mathematics Subject Classification:} 65M50, 65M60, 65M70, 65R05, 65.75
\vspace{5pt}


\newcommand{\h}{\hspace{1.cm}}
\newcommand{\hh}{\hspace{2.cm}}
\newtheorem{yl}{\hspace{1.cm}Lemma}
\newtheorem{dl}{\hspace{1.cm}Theorem}
\renewcommand{\sec}{\section*}
\renewcommand{\l}{\langle}
\renewcommand{\r}{\rangle}
\newcommand{\be}{\begin{eqnarray}}
\newcommand{\ee}{\end{eqnarray}}

\normalsize \vskip 0.2in
\newpage
%\baselineskip=2pc
% section 1
%\tableofcontents

\section{Introduction}
Solution interpolation or remapping between two deforming meshes is an
important component for several applications in scientific computing,
including arbitrary-Lagrangian-Eulerian (ALE) methods in computational fluid dynamics
\cite{Adam-Pavlidis-Percival-et2016,Anderson-Dobrev-et2015,Anderson-Dobrev-et2018, Barlow-Maire-et2016,
Bo-Shashkov-2015, ChengShu2008, Dukowicz1984, Dukowicz-Baumgardener2000,Dukowicz-Kodis1987,
Hirt-Amsden-Cook1974,  Kucharik-Shashkov-2014, Kucharik-Shashkov-Wendroff2003,Margolin-Shashkov2003}
and rezoning moving mesh (MM) methods in general numerical solution of partial differential equations (PDEs)
\cite{DiLiTangZhang2005, LiTang2006, LiTangZhang2001,  TangTang2003, Zhangzr2006}.
If not designed properly, a scheme for the interpolation can lead to violation of conservation
of some important physical quantities, deterioration of accuracy, and/or
introduction of spurious negative values in variables supposed to be nonnegative.

Some of the earliest work on conservative interpolation between two deformating meshes
grew out of the development of ALE methods \cite{Hirt-Amsden-Cook1974}.
Depending on the relation between the old (Lagrangian) and new (rezoned) meshes,
we can classify mesh-to-mesh interpolation algorithms as integral-remapping or advection-remapping ones.
If the two meshes are completely independent of one another or have the same connectivity but are
arbitrarily displaced with respect to each other, one needs to use integral-remapping interpolation
which involves finding the intersections of the cells of two meshes.
In \cite{Dukowicz1984}, a conservative interpolation method is proposed that assumes
piecewise constant fields and simplifies the problem of computing the volume of intersection
of old and new cells into a surface integral by invoking the divergence theorem.
However, the first-order nature of the method leads to excessive diffusion.
The approach is extended to second order to improve its diffusive characteristics in \cite{Dukowicz-Kodis1987}.
In \cite{Margolin-Shashkov2003}, a second-order accurate, conservative, and sign-preserving local remapping
algorithm for a positive, scalar, cell-centered function is developed based on the intersection,
which can be written in flux form if two meshes have the same connectivity.
Then the authors simplify it as a face-based donor-cell method,
which avoids finding the cell intersections but requires the displacements to be small
to maintain the positivity of the remapping variables.
The main drawback of this type of method is the difficulty of evaluating the integrals for arbitrary meshes,
especially in higher dimensions.

When the old and new meshes have the same connectivity, they can be viewed as a deformation
of each other and advection-remapping can be used. It is shown in \cite{Dukowicz-Baumgardener2000}
that if the physical time step is made sufficiently short such that node trajectories are confined to the nearest neighbor
cells (and thus the magnitude of the deformation is small), then the remapping can be written as
a flux-form convection algorithm. An incremental remapping method based on the solution of convection equations
is developed in \cite{Dukowicz-Baumgardener2000}.
A linearity-and-bound preserving conservative interpolation scheme is introduced
in \cite{Kucharik-Shashkov-Wendroff2003}.
A main advantage of an advection-based scheme is that it does not require finding
the intersections of old and new mesh cells.
However, the connection between advection equations and conservative interpolation/remapping does not seem
to be well understood as assumptions and discretization errors of using advection methods
for interpolation/remapping are not easily identified.

There is a different approach of advection-remapping where the interpolation
is viewed as solving a linear convection PDE over a pseudo-time interval.
For example, Li et al. \cite{LiTangZhang2001} use such an interpolation scheme in an MM finite element method.
A conservative interpolation scheme is proposed and used by Tang and Tang \cite{TangTang2003} for finite volume computation of hyperbolic equations. This scheme seems to work only for small mesh deformation.
A divergence-free-preserving interpolation algorithm is developed in \cite{DiLiTangZhang2005}
for the MM finite element computation of the incompressible Navier-Stokes equations.
It is worth pointing out that only one pseudo-time step is used in their computation
since the mesh deformation is very small.
The idea of \cite{DiLiTangZhang2005} is extended in \cite{LiTang2006}
to develop a second-order conservative interpolation scheme for use with an MM-DG method.
Anderson et al. \cite{Anderson-Dobrev-et2015} propose a method for remapping the state variables
of single material ALE based on solving convection equations using semi-discrete DG methods
and three nonlinear approaches to enforce monotonicity of the remapping variables.
Its multi-material extension and combination with the Lagrangian phase can be found
in \cite{Anderson-Dobrev-et2018}.
For more remapping/interpolation methods, the interested reader is referred to
\cite{Barlow-Maire-et2016, ChengShu2008, Kucharik-Shashkov-2014, Zhangzr2006} and the references therein.

The objective of this paper is to develop an arbitrary high-order
conservative interpolation scheme and present an analysis for its cost
and positivity preservation, two issues that have hardly been
studied for interpolation/remapping for deforming meshes.
The scheme is based on solving a linear convection equation with an MM-DG method
for spatial discretization and an explicit third-order Runge-Kutta scheme for time discretization.
This DG-interpolation scheme is shown to be mass-conservative and applicable for bounded but otherwise arbitrary
mesh deformation. Moreover, it is shown that the cost of the DG-interpolation is in the order of the number
of mesh vertices multiplied by the number of pseudo-time steps needed to integrate
the convection equation from the pseudo-time zero corresponding to the old mesh and the pseudo-time one
corresponding to the new mesh. This number of pseudo-time steps depends on the magnitude
of mesh deformation relative to the size of mesh elements in general.
It stays constant as the mesh is being refined when the mesh deformation is in the order of
the minimum element height, a typical situation in the MM solution of conservation laws
with an explicit scheme. On the other hand, the number of pseudo-time steps increases as the mesh is being
refined if the mesh deformation only stays bounded. A typical scenario of this is in the MM solution
of PDEs with a fixed physical time step size or with an implicit scheme.
Another issue is positivity preservation. Generally speaking,
the DG-interpolation alone may not preserve the positivity/nonnegativity
of the function to be interpolated. We consider a limiter
\cite{Liu-Osher1996,ZhangShu2010,ZhangXiaShu2012} that
uses a linear scaling around the positive cell average while conserving the cell average
and maintaining the convergence order of the DG discretization. We
show analytically and verify numerically that the DG-interpolation with the limiter can preserve the positivity
of the function to be interpolated.

As an application example, we study the use of the DG-interpolation scheme
in the rezoning MM-DG solution of the radiative transfer equation (RTE).
The RTE is an integro-differential equation modeling the interaction of radiation
with scattering and absorbing media and having important applications in various fields
in science and engineering. It involves time, space and angular variables
and contains an integral term in angular directions while being hyperbolic in space.
The challenges for its numerical solution include the needs to handle with its
high dimensionality, the presence of the integral term, the development of
discontinuities and sharp layers in its solution along spatial directions, and
the appearance of spurious negative values in the nonnegative radiative intensity.
These challenges make adaptive high-order DG methods amenable to the numerical
solution of RTE. Indeed, DG methods have been considered for RTE.
For example, a quasi-Lagrangian MM-DG method is proposed in \cite{ZhangChengHuangQiu} for RTE
and the preservation of nonnegativity of the radiative
intensity is investigated in \cite{LingChengShu2018,YuanChengShu2016,ZhangChengQiu2019}
for the DG solution of RTE on a fixed mesh.

We consider a rezoning MM-DG method (instead of a quasi-Lagrangian one)
for the numerical solution of RTE. It typically includes three steps,
mesh redistribution/adaptation, solution interpolation from the old mesh to the new one,
and solution of the physical equation on the new mesh.
The method has the advantages that these steps are independent of each other
and existing schemes can be used for each step.
Moreover, the task seems to be simpler here than that with a quasi-Lagrangian MM method
that strongly couples the effects of mesh movement with the discretization of RTE.
For the current situation, we deal separately with a scalar function/equation on
a moving mesh for the second step (interpolation) and the discretization of RTE on a fixed mesh for the third step.
In our computation, we use the positivity-preserving DG method of \cite{LingChengShu2018,YuanChengShu2016,ZhangChengQiu2019}
for spatial variables and the discrete-ordinate method (DOM) \cite{DOM} for angular variables.
For adaptive mesh generation (the first step), we use  an MM method
\cite{Huang2015,HKR2015,Huang2010}
which is known to produces a nonsingular moving mesh \cite{Huang2018}.
We use the DG-interpolation for the second step.
The whole computation can be made positivity preserving
when the computation at the second and third steps can be made positivity preserving.

To conclude the introduction, we would like to emphasize that the current work
contains a few new contributions.
As mentioned earlier, a number of remapping or interpolation schemes for deforming meshes
have been developed (e.g., see
\cite{DiLiTangZhang2005,Dukowicz-Baumgardener2000,Kucharik-Shashkov-Wendroff2003,
Margolin-Shashkov2003,TangTang2003}) however none of the existing schemes seems to work well
with large mesh deformation. Here we propose to use multiple pseudo-time steps and demonstrate
that the resulting DG-interpolation scheme works for meshes of large or small bounded deformation.
Moreover, we give a cost analysis for the scheme and particularly obtain an estimate
of the number of pseudo-time steps needed for each interpolation in terms of mesh deformation.
Furthermore, we consider high order accuracy, conservation, geometric conservation law,
and positivity preservation in the construction of the scheme. The proposed scheme
appears to be the first interpolation/remapping scheme taking all of those properties into consideration.
Finally, RTE proves to be a right application for the DG interpolation scheme.
Its implicit time integration means large time stepsize which in turn leads to large mesh
deformation between time steps. In addition, the positivity of the radiative intensity needs to
be preserved in the computation. Positivity preservation in the MM solution of RTE is new too.

The outline of the paper is as follows.
The high-order DG-interpolation is developed and
and its cost, mass conservation, and positivity preservation are analyzed in \S\ref{sec:DG-interpolation}.
The moving mesh PDE (MMPDE) method is described in \S\ref{sec:mmpde}.
Numerical results obtained for one- and two-dimensional examples are presented
in \S\ref{sec:interp-tests} to demonstrate the high-order accuracy and positivity-preserving features
and the cost of the DG-interpolation.
A rezoning MM-DG method for RTE is described in \S\ref{sec:DOM-MM-DG-RTE} and numerical examples
with various settings in one and two spatial dimensions are given in \S\ref{sec:RTE-tests}.
Finally, \S\ref{sec:conclusion} contains the conclusions.

% section 2
%\section{Arbitrary high-order conservative positivity-preserving DG-interpolation}
\section{High-order conservative positivity-preserving DG-interpol-\\ation}
\label{sec:DG-interpolation}

In this section we present an interpolation scheme from an old simplicial mesh to
a new one with high-order accuracy, mass conservation, and positivity preservation.
The scheme works in any dimension although we restrict our discussion in one and two dimensions
for notational simplicity.

Let $\mathcal{D} \subset \mathbb{R}^d$ ($d = 1$ and 2) be a polygonal bounded domain.
Assume that we are given nonsingular simplicial meshes $\mathcal{T}_h^{old}$
and $\mathcal{T}_h^{new}$ on $\mathcal{D}$
that have the same number of elements and vertices and the same connectivity.
They differ only in the location of vertices and can be considered a deformation of each other.
They can also be regarded as a moving mesh at different time instants.
In this work, we use the MMPDE method (see \S\ref{sec:mmpde}) to generate such a mesh.

The interpolation problem between $\mathcal{T}_h^{old}$ and $\mathcal{T}_h^{new}$
is equivalent to the numerical solution of the differential equation
\cite{Anderson-Dobrev-et2015, Huang2010, LiTang2006, LiTangZhang2001}
\begin{equation}
\label{pde}
\frac{\partial u}{\partial \varsigma} (\bm{x},\varsigma) =0,\quad (\bm{x},\varsigma)\in
\mathcal{D} \times(0,1]
\end{equation}
on the moving mesh $\mathcal{T}_h(\varsigma)$
obtained as a linear interpolant of $\mathcal{T}_h^{old}$
and $\mathcal{T}_h^{new}$ in the pseudo-time $\varsigma \in [0,1]$.
In particular, $\mathcal{T}_h(\varsigma)$ has the same number of elements and vertices
and the same connectivity as $\mathcal{T}_h^{old}$ and $\mathcal{T}_h^{new}$ and its nodal positions and
velocities (which can also be interpreted as deformation) are given by
\begin{align}
\label{location+speed}
&\bm{x}_i(\varsigma) = (1-\varsigma)\bm{x}_i^{old}+\varsigma \bm{x}_i^{new}, \quad i =
1,...,N_v
\\
&\dot{\bm{x}}_i = \bm{x}_i^{new}-\bm{x}_i^{old}, \quad i = 1,...,N_v .
\label{location+speed+1}
\end{align}
The initial condition is
\begin{equation}\label{u0}
u(\bm{x},0) = u_{0}(\bm{x}),\quad \bm{x}\in\mathcal{D}
\end{equation}
where $u_{0}(\bm{x})$ is the original function defined on $\mathcal{T}_h^{old}$.
We define the piecewise linear mesh velocity function as
\begin{equation}\label{Xdot-1}
\dot{\bm{X}}(\bm{x},\varsigma) = \sum_{i=1}^{N_v} \dot{\bm{x}}_i \phi_i(\bm{x},\varsigma),
\end{equation}
where $ \phi_i$ is the linear basis function associated with the vertex $\bm{x}_i$.

We consider the numerical solution of \eqref{pde} using a quasi-Lagrangian MM-DG method
\cite{LHQ2019,ZhangChengHuangQiu}.
Let $K$ be an arbitrary element of $\mathcal{T}_h(\varsigma)$.
Denote the basis functions of degree up to $r \geq 1$ on $K$ by $\phi_{K}^{[j]},~
j=1,..., n_b$, where $n_b\equiv  (r+d)!/(d!r!)$ is the number of basis functions.
Notice that $n_b =r+1$ for $d=1$ and $n_b =(r+1)(r+2)/2$ for $d=2$.
The DG finite element space is defined as
\begin{equation}\label{Vh}
V^{r}_h(\varsigma) = \{v\in L^2(\mathcal{D}):\; v |_{K}\in P^{r}(K),
\; \forall K\in \mathcal{T}_h(\varsigma) \},
\end{equation}
where $P^{r}(K)$ stands for the space of polynomials of degree at most $r$ on $K$.
Then any DG approximation polynomial $u_h\in V^{r}_h(\varsigma)$ can be expressed as
\begin{equation}
\label{uh}
u_h(\bm{x},\varsigma)= \sum_{j=1}^{n_b} u^{[j]}_{K}(\varsigma) \phi^{[j]}_{K}(\bm{x},\varsigma),
\quad \bm{x} \in K, \quad K \in \mathcal{T}_h(\varsigma)
\end{equation}
where $u^{[j]}_{K},~j=1,...,n_b$, are the degrees of freedom.
Without causing confusion, hereafter we will suppress the subscript ``$h$'' in $u_h$, i.e.,
we will write $u_h$ as $u$.
We note that the basis functions depend on $\varsigma$ due to the movement of the vertices.
From the fact that $K$ is a simplex,
it is not difficult to show  that
\begin{equation}
\label{d-phi}
\frac{\partial \phi_{K}^{[j]}}{\partial \varsigma} (\bm{x},\varsigma)=
-\nabla \phi_{K}^{[j]}(\bm{x},\varsigma)\cdot \dot{\bm{X}}(\bm{x},\varsigma),\quad a.e. ~\text{in}~ \mathcal{D}.
\end{equation}

For the weak formulation of \eqref{pde}, multiplying it by a test function $v \in V^{r}_h(\varsigma)$ and
integrating the resulting equation over $K$, we obtain
\begin{equation}
\label{test}
\int_{K}\frac{\partial u}{\partial \varsigma}v d\bm{x}=0.
\end{equation}
On the other hand, from the Reynolds transport theorem we have
\begin{align*}
\frac{d}{d \varsigma}\int_{K} u v d\bm{x} & =
\int_{K} \left ( v \frac{\partial u}{\partial \varsigma} + u \frac{\partial v}{\partial \varsigma} \right )d\bm{x}
+\int_{\partial K} uv \dot{\bm{X}} \cdot\bm{n}_K ds,
\end{align*}
where $\bm{n}_K$ is the outward unit normal to the boundary $\partial K$.
Using \eqref{d-phi} (with $\phi_{K}^{[j]}$ being replaced by $v$) and \eqref{test} in the above equation, we get
\begin{equation}
\label{div}
\frac{d}{d \varsigma}\int_{K} u v d\bm{x}
+\int_{\partial K} v\left (-u \dot{\bm{X}}\cdot \bm{n}_K\right )   ds
+\int_{K} (u \dot{\bm{X}}) \cdot \nabla  v d\bm{x}=0.
\end{equation}
The boundary integral term is replaced by a numerical flux in the DG approximation.
Thus, the semi-discrete MM-DG solution for \eqref{pde} is
to seek $u \in V^{r}_h(\varsigma)$, $0 < \varsigma \le 1$ such that
\begin{equation}
\label{semi-DG}
\frac{d}{d \varsigma}\int_{K} u v d\bm{x}
+\sum_{e \in \partial K} \int_{e} v F_{e}(u_{K}^{in}, u_{K}^{out}) ds
+\int_{K} (u \dot{\bm{X}}) \cdot \nabla  v d\bm{x}=0, \quad
\forall  v\in V^{r}_h(\varsigma)
\end{equation}
where $F_{e}(u_{K}^{in}, u_{K}^{out}) \approx -u \dot{\bm{X}} \cdot \bm{n}_K$ is
a numerical flux defined on $e \in \partial K$, $u_{K}^{in}$ denotes the value of $u$ on $K$,
and $u_{K}^{out}$ is the value of $u$ on the element (denoted by $K'$) sharing the common edge $e$
with $K$.
We use the local Lax-Friedrichs numerical flux, viz.,
\begin{align}
F_{e}(u_{K}^{in}, u_{K}^{out}) =
\frac{1}{2}\Big{(}
\big{(}-u_K^{in}\dot{\bm{X}}^{e}-u_K^{out} \dot{\bm{X}}^{e} \big{)} \cdot \bm{n}^e_K
-\alpha_{e} (u_{K}^{out}-u_{K}^{in})\Big{)},\quad \forall e \in \partial K
\label{llf-flux}
\end{align}
where $\dot{\bm{X}}^{e}$ denotes the restriction of $\dot{\bm{X}}$ on $e$ and
\begin{equation}\label{alpha-local}
\alpha_{e} = \max\big{(}
|\dot{\bm{X}}^{e}\cdot \bm{n}^e_K |,
|\dot{\bm{X}}^{e}\cdot \bm{n}^e_{K'} |\big{)} .
\end{equation}
Note that this numerical flux is actually an upwind flux and vanishes on the boundary
of the domain due to the fact that the boundary does not move.
It satisfies several properties including consistency, monotonicity, Lipschitz continuity, and conservativeness, with the last property being expressed as
\begin{equation}\label{flux-c4}
  F_{e}(u_{K}^{in}, u_{K}^{out}) + F_{e}(u_{K'}^{in}, u_{K'}^{out})=0.
\end{equation}
In our computation, the second and third terms in the left of \eqref{semi-DG} are computed using Gaussian quadrature rules.

The third-order explicit total variation diminishing (TVD) Runge-Kutta scheme is used
to discretize \eqref{semi-DG} in time. To describe the scheme, we rewrite \eqref{semi-DG} into
\begin{align}
\label{ode}
& \frac{d}{d \varsigma}\int_{K}u vd\bm{x} = - \mathcal{A}(u ,v)|_{K}
\equiv - \sum_{e\in \partial K} \int_{e} vF_{e}(u_{K}^{in}, u_{K}^{out})  ds
-\int_{K} (u \dot{\bm{X}}) \cdot \nabla  v d\bm{x} .
\end{align}
Let the time instants be
\[
0=\varsigma^0<\varsigma^1<\cdots<\varsigma^{\nu}<\varsigma^{\nu+1}<\cdots <\varsigma^{N_\varsigma}=1,
\quad \text{and}\quad \Delta \varsigma^{\nu} = \varsigma^{\nu+1}- \varsigma^{\nu}.\]
The third-order explicit TVD Runge-Kutta scheme for \eqref{ode} reads as
\begin{equation}
\label{third}
\begin{cases}
& \int_{K^{\nu,(1)}}u ^{(1)}v^{\nu,(1)} d\bm{x}
= \int_{K^{\nu}}u ^{\nu}v^{\nu} d\bm{x}-\Delta \varsigma^\nu \mathcal{A}(u ^{\nu},v^{\nu})|_{K^\nu},
\\
& \int_{K^{\nu,(2)}}u ^{(2)}v^{\nu,(2)} d\bm{x}
= \frac{3}{4}\int_{K^{\nu}}u ^{\nu}v^{\nu}d\bm{x} \\
& \qquad \qquad \qquad
+\frac{1}{4}\Big (\int_{K^{\nu,(1)}}u^{(1)}v^{\nu,(1)} d\bm{x}
-\Delta \varsigma^\nu \mathcal{A}(u^{(1)},v^{\nu,(1)})|_{K^{\nu,(1)}}\Big ),
\\
& \int_{K^{\nu+1}}u^{\nu+1}v^{\nu+1} d\bm{x}
= \frac{1}{3}\int_{K^{\nu}}u^{\nu}v^{\nu} d\bm{x} \\
& \qquad \qquad \qquad
+\frac{2}{3}\Big (\int_{K^{\nu,(2)}}u^{(2)} v^{\nu,(2)} d\bm{x}
-\Delta \varsigma^\nu \mathcal{A}(u^{(2)},v^{\nu,(2)})|_{K^{\nu,(2)}}\Big ) ,
\end{cases}
\end{equation}
where $u^{(1)}$, $v^{\nu,(1)}$, $K^{\nu,(1)}$ are stage values
at $\varsigma = \varsigma^{\nu + 1}$,
$u^{(2)}$, $v^{\nu,(2)}$, $K^{\nu,(2)}$ are the values
at $\varsigma = \varsigma^{\nu+\frac{1}{2}}$, and
$u^{\nu+1}$, $v^{\nu+1}$, $K^{\nu+1}$ are at $\varsigma = \varsigma^{\nu + 1}$.
It is emphasized that {\em the coordinates of the vertices and the volume of $K$ need to
be updated at these stages}.
Especially, as will be seen in \S\ref{sec:GCL}, a special update scheme for the
element volume may be needed
for the scheme to satisfy the so-called geometric conservation law \cite{TL79,TT61}.
It is also worth pointing out that the test functions at $K^{\nu}$, $K^{\nu+1}$, $K^{\nu,(1)}$,
and $K^{\nu,(2)}$ are connected through their counterparts on the reference element $\widehat{K}$.
Indeed, for any $\hat{v} \in P^r(\widehat{K})$, we have
\begin{equation}
v^{\nu} = \hat{v} \circ {F}^{-1}_{K^{\nu}},\quad
v^{\nu+1} = \hat{v} \circ {F}^{-1}_{K^{\nu+1}},\quad
v^{\nu,(1)} = \hat{v} \circ {F}^{-1}_{K^{\nu,(1)}},\quad
v^{\nu,(2)} = \hat{v} \circ {F}^{-1}_{K^{\nu, (2)}} ,
\label{test-function-1}
\end{equation}
where $F_K$ is the affine mapping from $\widehat{K}$ to $K$ for $K = K^{\nu}$, $K^{\nu+1}$, $K^{\nu,(1)}$,
or $K^{\nu,(2)}$.

The time step size $\Delta \varsigma$ is chosen to ensure the stability of the scheme \cite{DG-review}, i.e.,
\begin{equation}
\label{DG-Interp-cfl-0}
\Delta \varsigma = \frac{C_{cfl}}{\max\limits_{e,K} |\dot{\bm{X}}^e \cdot \bm{n}^e_K | }\cdot
\min (h_{min}^{old}, h_{min}^{new}),
\end{equation}
where $C_{cfl}$ is a constant typically chosen to be less than $1/(2r+1)$ and
$h_{min}^{old}$ and $h_{min}^{new}$ are the minimum element height for the old and new meshes, respectively.

From the theory of DG and TVD Runge-Kutta scheme (e.g., see \cite{ZhangQiang-Shu2010}),
we can expect that the above described DG-interpolation scheme is
$(r+1)$th order in space and third order in time for problems with smooth solutions, viz.,
$
\mathcal{O} (\Delta \varsigma^3) + \mathcal{O} (h^{r+1}),
$
where $h$ denotes the maximum element diameter. Particularly, the scheme is
second order for $r=1$ and third order for $r=2$.
For $r>2$, we can choose a smaller $\Delta \varsigma$ or a higher-order time scheme
such that the temporal error is negligible.

It is emphasized that the above described scheme does not require any prior conditions
on the meshes $\mathcal{T}_h^{old}$ and $\mathcal{T}_h^{new}$. Particularly, it works
when the mesh has large deformation although more time steps may be needed.
The cost of the scheme is discussed in \S\ref{sec:Nsteps}.

\subsection{The geometric conservation law (GCL)}
\label{sec:GCL}

GCL stands for geometric identities that hold in continuous form.
They may no longer hold in a discrete setting especially in the computation with moving meshes \cite{TL79,TT61}.
A simple verification for satisfying GCL is to use uniform flow reproduction,
i.e., to check if the underlying scheme
produces a uniform flow if the initial flow is uniform. Theoretical and numerical analysis
(e.g., see \cite{BG04,FN04}) shows that satisfying GCL is neither a necessary nor
a sufficient condition for the stability of a scheme but often helps improve the accuracy
and stability of the computation.
We study \eqref{third} here for the satisfaction of GCL.

Taking $u=1$ in \eqref{ode} and using $F_{e}(1,1)= -\dot{\bm{X}}^{e}\cdot \bm{n}^e_K$ and
the divergence theorem, we have
\begin{align}
\mathcal{A}(1,v)|_{K}
&=\sum_{e\in \partial K} \int_{e} v(-\dot{\bm{X}}^{e}\cdot \bm{n}^e_K )ds
+\int_{K} \dot{\bm{X}} \cdot \nabla  v d\bm{x}
\notag \\
%&=\sum_{e\in \partial K} \int_{e} v(-\dot{\bm{X}}^{e}\cdot \bm{n}^e_K )ds
%+\sum_{e\in \partial K} \int_{e} v(\dot{\bm{X}}^{e}\cdot \bm{n}^e_K )ds
%-\int_{K} v\nabla\cdot\dot{\bm{X}} d\bm{x}
%\notag \\
&=-\int_{K} v\nabla \cdot\dot{\bm{X}}  d\bm{x} = - \nabla \cdot\dot{\bm{X}}|_{K} \int_{K} v  d\bm{x}
=-\frac{|K|}{|\hat{K}|} \nabla\cdot\dot{\bm{X}}|_{K} \int_{\widehat{K}}\hat{v}d\bm{\hat{x}} .
\label{A-1v}
\end{align}
Combining this equation with \eqref{ode} and taking $u= 1$ and $v=1$, we get
\begin{equation}
\label{area-K}
\frac{d}{d \varsigma}|K| = |K|\nabla\cdot \dot{\bm{X}}|_{K},
\end{equation}
which is the GCL governing the evolution of the volume of element $K$.
On the other hand, taking $u^{\nu} = 1$, $u^{(1)} = 1$, $u^{(2)} = 1$, $u^{\nu+1} = 1$, and
$\hat{v} = 1$ in \eqref{third}, we obtain
\begin{equation}
\begin{cases}
|K^{\nu,(1)}|
= |K^{\nu}|+\Delta \varsigma^\nu |K^{\nu}|\nabla\cdot\dot{\bm{X}}|_{K^{\nu}},
\\
|K^{\nu,(2)}|
= \frac{3}{4}|K^{\nu}|
+\frac{1}{4}\big (|K^{\nu,(1)}|
+\Delta \varsigma^\nu |K^{\nu,(1)}|\nabla\cdot\dot{\bm{X}}|_{K^{\nu,(1)}}\big ),
\\
|K^{\nu+1}|
= \frac{1}{3}|K^{\nu}|
+\frac{2}{3}\big ( |K^{\nu,(2)}|
+\Delta \varsigma^\nu |K^{\nu,(2)}|\nabla\cdot\dot{\bm{X}}|_{K^{\nu,(2)}}\big ),
\end{cases}
\label{Karea-1}
\end{equation}
which can be used to update the volume of $K$ at the three Runge-Kutta stages.
The above equation can also be obtained by applying the third-order Runge-Kutta scheme
directly to \eqref{area-K}.
The time-stepping (\ref{Karea-1}) has been derived by Cheng and Shu \cite{ChengShu2007JCP}
for a GCL-preserving ALE formulation and extended to ALE-DG and Lagrangian-DG
more recently by Pandare et al. \cite{Pandare-etal-2016,Pandare-etal-2018}.

\begin{lemma}
\label{GCL}
The fully-discrete MM-DG scheme \eqref{third} reproduces the uniform flow,
i.e., $u^{\nu} \equiv 1$ implies $u^{\nu+1} \equiv 1$, if the element volume is
updated according to \eqref{Karea-1}.
\end{lemma}

This lemma can be proved by taking $u^{\nu}\equiv 1$ in \eqref{third}
and using \eqref{A-1v} and \eqref{Karea-1}.

As mentioned above, the volume of $K$ at different Runge-Kutta stages can be obtained using \eqref{Karea-1}.
It can also be calculated directly using the coordinates of the vertices. Interestingly, it can be verified
that these two approaches are the same in one dimension but different in two
and higher dimensions. In the latter case, \eqref{Karea-1} needs to be used for uniform flow reproduction
and thus GCL satisfaction.

% section 2.2
\subsection{Mass conservation}
\label{mass-cons}

In this subsection we show that the DG-interpolation scheme \eqref{semi-DG} and
\eqref{third} conserves the mass.

\begin{lemma}
\label{semi-mass}
The semi-discrete MM-DG scheme \eqref{semi-DG} conserves the mass.
\end{lemma}

This lemma can be proved by
taking $v=1$ in \eqref{semi-DG}, summing the resulting equation over all elements,
re-arranging the terms according to interior and boundary edges, and using \eqref{flux-c4}
and the fact that the numerical flux vanishes on the boundary.

\begin{lemma}
\label{fully-mass}
The fully discrete MM-DG scheme \eqref{third} conserves the mass.
\end{lemma}

This lemma can be proved similarly as for Lemma~\ref{semi-mass}.

%\begin{proof}
%Taking the test function $v=1$ in \eqref{third}, summing the resulting equations over all of the elements,
%and using a similar argument as in the proof of Lemma~\ref{semi-mass}, we can get
%\begin{align}
%&\sum_{K^{\nu,(1)}}\int_{K^{\nu,(1)}} u^{(1)}  d\bm{x}
%=\sum_{K^{\nu}}\int_{K^{\nu}} u^{\nu}  d\bm{x},
%\label{f3-3-1}
%\\
%&\sum_{K^{\nu,(2)}}\int_{K^{\nu,(2)}} u^{(2)}  d\bm{x}
%=\frac{3}{4}\sum_{K^{\nu}}\int_{K^{\nu}} u^{\nu}  d\bm{x}
%+\frac{1}{4}\sum_{K^{\nu,(1)}}\int_{K^{\nu,(1)}} u^{(1)}  d\bm{x},
%\label{f3-3-2}
%\\
%&\sum_{K^{\nu+1}}\int_{K^{\nu+1}} u^{\nu+1}  d\bm{x}
%=\frac{1}{3}\sum_{K^{\nu}}\int_{K^{\nu}} u^{\nu}  d\bm{x}
%+\frac{2}{3}\sum_{K^{\nu,(2)}}\int_{K^{\nu,(2)}} u^{(2)} d\bm{x}.
%\label{f3-3-3}
%\end{align}
%Combining these we obtain \eqref{f3-4}, which implies that the mass is conserved
%by the scheme \eqref{third}.
%\end{proof}

\begin{rem}
\label{other-time}
{\em
Similarly, we can prove that the
first-order forward Euler scheme and the second-order explicit TVD Runge-Kutta
scheme also conserve the mass when applied to \eqref{semi-DG}.
}
\end{rem}

% section 2.3
\subsection{Cost of the DG-interpolation}
\label{sec:Nsteps}

We now investigate the cost of the DG-interpolation scheme \eqref{third}.
We start with noticing that the cost of each time step of the scheme
is $\mathcal{O}(N_v)$ and the total cost is
$\mathcal{O}(N_v N_\varsigma)$, where $N_v$ is the number of the mesh vertices and $N_\varsigma$
is the number of  time steps to reach $\varsigma = 1$.
Note that this total cost is the cost for each interpolation of the function from the old mesh to the new one.
The key to the estimation of this cost is to estimate $N_\varsigma$.

To this end, we recall that the CFL stability condition (\ref{DG-Interp-cfl-0}).
Since $\dot{\bm{X}}$ is piecewise linear, from \eqref{location+speed+1} we have
\[
\max_{e,K}|\dot{\bm{X}}\cdot \bm{n}^e_K| \sim \max_{i}|\bm{x}^{old}_i - \bm{x}^{new}_i | .
\]
Then, \eqref{DG-Interp-cfl-0} becomes
\begin{equation}
\label{DG-Interp-cfl}
\Delta \varsigma = \frac{C_{cfl}}{\max_{i}|\bm{x}^{old}_i - \bm{x}^{new}_i | }
\cdot \min (h_{min}^{old}, h_{min}^{new}) .
\end{equation}
This indicates that {\em $\Delta \varsigma$ and thus $N_\varsigma$ depend on the magnitude
of mesh deformation relative to the size of mesh elements}.
In the following we consider two special cases.

{\bf Case 1.} In the first case we consider the situation where
\begin{equation}
\label{DG-Interp-cfl+1}
\max_{i}|\bm{x}^{old}_i - \bm{x}^{new}_i |=\mathcal{O}(\min (h_{min}^{old}, h_{min}^{new})).
\end{equation}
Then, \eqref{DG-Interp-cfl} implies that the DG-interpolation only takes
a constant number of time steps to reach $\varsigma = 1$ and
its total cost is $\mathcal{O}(N_v)$.

An extreme situation for \eqref{DG-Interp-cfl+1} is that the mesh is fixed.
Then we have $\max_{i}|\bm{x}^{old}_i$ $-\bm{x}^{new}_i | = 0$ and the upper bound
of \eqref{DG-Interp-cfl} becomes infinity, which means just one step is needed
for the DG-interpolation.

In the context of the MM solution of first-order hyperbolic equations,
$\mathcal{T}_h^{old}$ and $\mathcal{T}_h^{new}$ correspond to meshes at consecutive time steps,
i.e., $\mathcal{T}_h^{old} = \mathcal{T}_h^{n}$ and $\mathcal{T}_h^{new} = \mathcal{T}_h^{n+1}$,
where $n$ stands for the index for the physical time step, and
the time step size used for integrating the physical equations is typically chosen as
\begin{equation}
\Delta t = \mathcal{O}\big(\min (h_{min}^{old}, h_{min}^{new})\big)
\label{dt-explicit}
\end{equation}
to ensure stability. If the mesh velocities are bounded, i.e.,
\begin{equation}
\max_{i}\left | \frac{\bm{x}^{n}_i - \bm{x}^{n+1}_i }{\Delta t}\right | = \mathcal{O}(1)
\quad \text{or} \quad \max_{i}|\bm{x}^{n}_i - \bm{x}^{n+1}_i |=\mathcal{O}(\Delta t),
\label{DG-Interp-cfl+2}
\end{equation}
then (\ref{dt-explicit}) implies \eqref{DG-Interp-cfl+1}. As a consequence, we can expect
that the cost for each DG-interpolation in the MM solution of hyperbolic equations
is $\mathcal{O}(N_v)$.

{\bf Case 2.} In this case we consider the situation with
\begin{equation}
\label{DG-Interp-cfl+1+2}
\max_{i}|\bm{x}^{old}_i - \bm{x}^{new}_i |=\mathcal{O}(1) .
\end{equation}
Then (\ref{DG-Interp-cfl}) means that the number of the time steps needed is
\begin{equation}
\label{DG-Interp-cfl+1+3}
N_\varsigma = \mathcal{O}\left ( \frac{1}{\min (h_{min}^{old}, h_{min}^{new})} \right ),
\end{equation}
which is $\mathcal{O}(N^{\frac{1}{d}})$ at the minimum (where $N$ is the number of elements).
It clearly indicates that $N_\varsigma$ increases
as the mesh is being refined.

A typical scenario for this case is when the physical PDE is integrated with an implicit scheme
and the physical time step size $\Delta t$ is taken independent of the mesh size
(in contrast to \eqref{dt-explicit}).
Then we have $\max_{i}|\bm{x}^{old}_i - \bm{x}^{new}_i |=\mathcal{O}(\Delta t)$
and
\begin{equation}
\label{DG-Interp-cfl+1+1}
N_\varsigma = \mathcal{O}\left ( \frac{\Delta t}{\min (h_{min}^{old}, h_{min}^{new})} \right ),
\end{equation}
which increases as the mesh is being refined.

%It can be shown that
%\begin{equation}
%\label{DG-Interp-cfl+1}
%\max_{i}|\bm{x}^{old}_i - \bm{x}^{new}_i |=\mathcal{O}(h_{min}^{old})
%\end{equation}
%implies (\ref{DG-Interp-cfl+1+1}).

\begin{rem}
{\em
The condition \eqref{DG-Interp-cfl+1} has been used in \cite{YangHuangQiu2012} to restrict the mesh movement in the MM WENO solution of conservation laws.
}
\end{rem}

\begin{rem}
{\em
It is interesting to mention that the interpolation schemes in \cite{DiLiTangZhang2005,TangTang2003}
for the rezoning MM methods and
in \cite{Dukowicz-Baumgardener2000,Kucharik-Shashkov-Wendroff2003,Margolin-Shashkov2003}
for ALE methods can be viewed as the one-step implementation of
some explicit schemes for integrating \eqref{pde} on a moving mesh.
These schemes have been observed \cite{DiLiTangZhang2005,Dukowicz-Baumgardener2000,Kucharik-Shashkov-Wendroff2003,
Margolin-Shashkov2003,TangTang2003} to work only for small mesh deformation.
This may be explained using \eqref{DG-Interp-cfl} and \eqref{DG-Interp-cfl+1},
i.e., \eqref{DG-Interp-cfl+1} (which implies small mesh deformation)
needs to be held if we want the right-hand side of \eqref{DG-Interp-cfl} to be constant.
The analysis in this subsection also shows that multiple steps are needed if large
mesh deformation is allowed.
}
\end{rem}

\subsection{Preservation of positivity}
\label{positivity}

It should be pointed out that the above described DG-interpolation scheme \eqref{third} cannot preserve
the positivity of the solution in general.
In this subsection, we consider a positivity-preserving (PP) limiter that
uses a linear scaling around nonnegative cell averages, conserves the cell averages,
and maintains the accuracy order of the original DG-interpolation.
The approach we use here is similar to the general techniques developed in \cite{Liu-Osher1996,ZhangShu2010,ZhangXiaShu2012}
for constructing high-order PP DG schemes on fixed meshes for scalar conservation laws.
To save space, we only discuss the forward Euler time discretization here.
The conclusion will hold over for the third-order explicit TVD Runge-Kutta method
since it is a convex combination of the forward Euler scheme.

The Euler scheme for the semi-discrete MM-DG scheme \eqref{semi-DG} is given by
\begin{align}
& \int_{K^{\nu+1}} (uv)^{\nu+1}d\bm{x} = \int_{K^{\nu}} (uv)^{\nu}d\bm{x}
\notag \\
& \qquad -\Delta \varsigma^\nu \Big{(}\sum_{e \in \partial K^{\nu}} \int_{e} v^{\nu} F_{e}((u^{\nu}_{K})^{int},
(u^{\nu}_{K})^{ext}) ds
+\int_{K^{\nu}} (u^\nu \dot{\bm{X}}^\nu) \cdot \nabla  v^{\nu} d\bm{x}\Big{)}.
\label{first}
\end{align}
%Denote the cell average of $u$ on $K$ by $\bar{u}_{K}$, i.e.,
%\begin{equation}\label{u-ave}
%\bar{u}_{K}= \frac{1}{|K|}\int_{K} u  d \bm{x}.
%\end{equation}
Taking $v=1$ in \eqref{first}, we obtain the evolution equation of the cell average $\bar{u}$ as
% (denoted $\bar{u}_{K}= \frac{1}{|K|}\int_{K} u  d \bm{x}$)
\begin{equation}\label{pp-cell}
|K^{\nu+1}| \bar{u}_{K^{\nu+1}}
= |K^{\nu}| \bar{u}_{K^{\nu}}
-\Delta \varsigma^\nu \sum_{e \in \partial K^{\nu}} \int_{e} F_{e}((u^{\nu}_{K})^{int}, (u^{\nu}_{K})^{ext}) ds .
\end{equation}

\begin{proposition}
\label{DG-interp-pp-2d}
For \eqref{first} and \eqref{pp-cell}, if $\bar{u}_{K^\nu} \ge 0$ and $u_{K^\nu}(\bm{\hat{x}})\ge 0$
for all $\bm{\hat{x}}\in G_{K^\nu}$ and all $K^\nu \in \mathcal{T}_h^{\nu}$, where $G_{K^\nu}$ is
a set of special quadrature points (e.g., \cite{ZhangXiaShu2012}) on $K^\nu$,
then $ \bar{u}_{K^{\nu+1}} \geq 0$ hold for all $K^{\nu +1} \in \mathcal{T}_h^{\nu + 1}$
under the CFL condition
\begin{equation}\label{pp-2d-cfl}
\Delta \varsigma^\nu \leq \frac{\frac{2}{3}\hat{w}_{1}}{
\max\limits_{K^\nu} \max\limits_{e\in \partial K^\nu}|\dot{\bm{X}}^{e}\cdot \bm{n}^e|}\cdot \min_{K^\nu}
\frac{| K^{\nu}|}{|\partial K^\nu|} ,
\end{equation}
where $\hat{w}_{1}$ is the first point weight of the $n_g$-point Gauss-Lobatto quadrature ($2n_g-3\geq r$)
and has the value of $1/2$ and $1/6$ for $r=1$ and $r=2$, respectively.
\end{proposition}

This proposition can be proved by following \cite{ZhangXiaShu2012} and using the mesh nonsingularity
assumption $|K^{\nu}|>0$ which is warranted by the MMPDE method; see \S\ref{sec:mmpde}.

Once we have $ \bar{u}_{K^{\nu+1}} \geq 0$, we can apply the linear scaling PP limiter as
\begin{equation}\label{scaling-2d}
 \begin{split}
\hat{u}_{K^{\nu + 1 }}
= \lambda_{K^{\nu + 1 }} (u_{K^{\nu + 1 }} -\bar{u}_{K^{\nu + 1 }})+ \bar{u}_{K^{\nu + 1 }},
\quad \forall \bm{x}\in K^{\nu + 1}, \quad K^{\nu + 1} \in \mathcal{T}_h^{\nu + 1}
  \end{split}
 \end{equation}
 where
\begin{align}\label{lambda-pp-2d}
\lambda_{K^{\nu + 1 }} =
\frac{\bar{u}_{K^{\nu + 1 }}}{\bar{u}_{K^{\nu + 1 }}-z_{K^{\nu + 1 }}},
\quad
z_{K^{\nu + 1 }} = \min_{\bm{\hat{x}}\in G_{K^{\nu + 1}}}
\big{\{}u_{K^{\nu + 1 }}(\bm{\hat{x}}),0\big{\}}.
 \end{align}
It can be verified that $\bar{\hat{u}}_{K^{\nu + 1 }} = \bar{u}_{K^{\nu + 1 }}$,
$\hat{u}_{K^{\nu + 1 }} (\bm{\hat{x}}) \ge 0$ for all
$\bm{\hat{x}}\in G_{K^{\nu + 1}}$, and $\hat{u}_{K^{\nu + 1 }}$
maintains the DG convergence order \cite{Liu-Osher1996,ZhangXiaShu2012}.

Finally, we note that if the initial solution is nonnegative, we have $\bar{u}_{K^0} \ge 0$.
By applying the linear scaling PP limiter, we can obtain an initial approximation that meets the assumption
of Proposition~\ref{DG-interp-pp-2d}. Hence, we conclude that the DG-interpolation
preserves the nonnegativity of the solution when the PP limiter is applied.

% section 4
\section{The MMPDE moving mesh method}
\label{sec:mmpde}

In this section we describe the generation of the new mesh $\mathcal{T}^{new}_{h}$
from the old one $\mathcal{T}^{old}_{h}$ using the MMPDE method \cite{Huang2010}.
We use here a new implementation of the method proposed in \cite{Huang2015}.
Adaptive meshes generated using this method are used in \S\ref{sec:interp-tests}
for the numerical examination of the DG-interpolation scheme and in \S\ref{sec:RTE-tests}
for the numerical solution of the radiative transfer equation.

To describe the MMPDE method, we introduce a computational mesh
$\mathcal{T}_c= \{\bm{\xi}_1, \;...,\;\bm{\xi}_{N_v} \}$ which serves as an intermediate variable,
and an almost uniform reference computational mesh
$\hat{\mathcal{T}}_c = \{\hat{ \bm{\xi}}_1,\; ...,\; \hat{ \bm{\xi}}_{N_v} \}$
which keeps fixed in the computation.
A key idea of the MMPDE method is to view any nonuniform mesh as a uniform one
in some metric $\mathbb{M}$ \cite{Huang2010}.
$\mathbb{M} = \mathbb{M}(\bm{x})$ is a symmetric and uniformly positive definite
matrix-valued function defined on $\mathcal{D}$. It provides the information needed for determining the size, shape, and orientation of the mesh elements throughout the domain.
Various metric tensors have been proposed; e.g., see \cite{Huang2010,Huang2003}.
We use here a metric tensor based on the Hessian of the computed solution.
To be specific, we consider a physical variable $u$ and denote its finite element approximation
by $u_h$. Let $H_K$ be a recovered Hessian of $u_h$ on $K\in \mathcal{T}_{h}$
for a mesh $\mathcal{T}_{h}$.
A number of strategies can be used for Hessian recovery for finite element approximations; e.g., see \cite{Buscaglia1998, Kamenski-PhD2009,Zhang-Naga2005,Zienkiewicz-Zhu-1992-1}.
Least square fitting \cite{Zhang-Naga2005} is used in our computation.
Denoting
\[
|H_K| = Q\hbox{diag}(|\lambda_1|,...,|\lambda_d|)Q^T,
\]
where $Q\hbox{diag}(\lambda_1,...,\lambda_d)Q^T$ is the eigen-decomposition of $H_K$,
the metric tensor is then defined as
\begin{equation}\label{mer}
\mathbb{M}_{K} =\det \big{(}\mathbb{I}+|H_K|\big{)}^{-\frac{1}{d+4}}
\big{(}\mathbb{I}+|H_K|\big{)},
\quad \forall K \in \mathcal{T}_h
\end{equation}
where $\mathbb{I}$ is the identity matrix and $\det(\cdot)$ is the determinant of a matrix.
The metric tensor \eqref{mer} is known \cite{Huang2003} to be optimal for the $L^2$-norm of
linear interpolation error.
%We use a least squares fitting strategy for Hessian recovery
%and take $\alpha_{h}=1$ in our computation.
For situations with several physical variables, we first compute the metric tensor for each of the variables
and then obtain the final metric tensor by matrix intersection.
%\cite[\S5.3.4]{Huang2010}.
%Let $A$ and $B$ be two symmetric and positive definite matrices.
%There exists a nonsingular matrix
%$P$ such that $P A P^T = \mathbb{I}$ and $P B P^T = \text{diag}(b_1,..., b_d)$.
%The intersection of matrices $A$ and $B$ is then defined as
%\begin{equation}\label{intersection}
%A\cap B = P^{-1} \text{diag}(\max(1, b_1), ...,\max(1, b_d)) P^{-T}.
%\end{equation}

When $\mathcal{T}_h$ is uniform
in the metric $\mathbb{M}$ in reference to the computational mesh $\mathcal{T}_c$,
it is known \cite{Huang2010}
%\cite{Huang2006, Huang2010}
that it satisfies
the equidistribution and alignment conditions,
\begin{align}
\label{ec}
& |K|\sqrt{\det(\mathbb{M}_K)}=\frac{\sigma_h|K_c|}{|\mathcal{D}_c|},&\quad \forall K\in \mathcal{T}_h
\\
\label{al}
& \frac{1}{d}\hbox{tr}\big{(} (F'_K)^{-1}\mathbb{M}_K^{-1}(F'_K)^{-T}\big{)}
= \hbox{det}\big{(}
(F'_K)^{-1}\mathbb{M}_K^{-1}(F'_K)^{-T}\big{)}^{\frac{1}{d}},&\quad \forall K \in \mathcal{T}_h
\end{align}
where $F'_K$ is the Jacobian matrix of the affine mapping: $F_K: K_c\in \mathcal{T}_c \rightarrow K \in \mathcal{T}_h$,
$\mathbb{M}_K$ is the average of $\mathbb{M}$ over $K$,
$\hbox{tr}(\cdot)$ denotes the trace of a matrix, and
\[
|\mathcal{D}_c|=\sum\limits_{K_c\in\mathcal{T}_c}|K_c|,\quad \sigma_h
=\sum\limits_{K\in\mathcal{T}_h}|K|\hbox{det}(\mathbb{M}_K)^{\frac{1}{2}} .
\]
The condition \eqref{ec} determines the size of elements through
the metric tensor $\mathbb{M}$. On the other hand, \eqref{al}, derived from requiring
$K$ (measured in the metric $\mathbb{M}_K$) to be similar to $K_c$ (measured in the Euclidean metric),
determines the shape and orientation of $K$ through $\mathbb{M}$ and shape of $K_c$.
An energy function associated with these conditions is given by
\begin{equation}\label{energy}
\begin{split}
\mathcal{I}_h(\mathcal{T}_h,\mathcal{T}_c)
=&\frac{1}{3}\sum_{K\in\mathcal{T}_h}|K|\hbox{det}(\mathbb{M}_K)^{\frac{1}{2}}\big{(}
\hbox{tr}((F'_K)^{-1}\mathbb{M}^{-1}_K(F'_K)^{-T})\big{)}^{\frac{3 d}{4}}
\\&+\frac{1}{3} d^{\frac{3 d}{4}}\sum_{K\in\mathcal{T}_h}|K|\hbox{det}(\mathbb{M}_K)^{\frac{1}{2}}
\left (\hbox{det}(F'_K) \hbox{det}(\mathbb{M}_K)^{\frac{1}{2}} \right )^{-\frac{3}{2}},
\end{split}
\end{equation}
which is a Riemann sum of a continuous functional developed in \cite{Huang2010}
based on mesh equidistribution and alignment.

Note that $\mathcal{I}_h(\mathcal{T}_h,\mathcal{T}_c)$ is a function of
the vertices $\bm{\xi}_i,\; i=1,...,N_v$, of $\mathcal{T}_c$
and the vertices $\bm{x}_i,\;i=1,...,N_v$, of $\mathcal{T}_h$.
Here, we adopt an indirect approach with which we take $\mathcal{T}_h$ as $\mathcal{T}_h^{old}$,
minimize $\mathcal{I}_h(\mathcal{T}_h^{old},\mathcal{T}_c)$ with respect to $\mathcal{T}_c$,
and obtain the new physical mesh through the relation between
$\mathcal{T}_h^{old}$ and newly obtained $\mathcal{T}_c$.
The mesh equation is defined as the gradient system of the energy function (the MMPDE approach), i.e.,
\begin{equation}\label{MM}
\begin{split}
\frac{d \bm{\xi}_i }{d\varsigma}
=-\frac{\hbox{det}(\mathbb{M}(\bm{x_i}))^{\frac{1}{2}} }{\tau}
\Big{(}\frac{\partial \mathcal{I}_h(\mathcal{T}^{old}_h,\mathcal{T}_c)}{\partial \pmb{\xi}_i}\Big{)}^T,
\quad i=1,...,N_v
\end{split}
\end{equation}
where ${\partial \mathcal{I}_h }/{\partial \bm{\xi}_i}$ is considered as a row vector and
$\tau>0$ is a parameter used to adjust the response time of mesh movement to
the changes in $\mathbb{M}$.

Let $\mathbb{J}=(F'_K)^{-1} = E_{K_c}E_K^{-1}$ with $
E_K=[\bm{x}_1^K-\bm{x}_0^K, ..., \bm{x}_d^K-\bm{x}_0^K]$ and
$E_{K_c}=[\bm{\xi}_1^K - \bm{\xi}_0^K, ..., \bm{\xi}_d^K -\bm{\xi}_0^K ]$, and define the function $G$ associated with the energy \eqref{energy} as
\begin{equation}\label{G}
G(\mathbb{J},\hbox{det}(\mathbb{J}))=
\frac{1}{3}\hbox{det}(\mathbb{M}_{K})^{\frac{1}{2}}
(\hbox{tr}(\mathbb{J}\mathbb{M}_{K}^{-1}\mathbb{J}^T))^{\frac{3 d}{4}}
+\frac{d^{\frac{3 d}{4}}}{3} \hbox{det}(\mathbb{M}_{K})^{\frac{1}{2}}
\left (\frac{\hbox{det}(\mathbb{J})}{\hbox{det}(\mathbb{M}_{K})^{\frac{1}{2}}} \right )^{\frac{3}{2}}.
\end{equation}

Using the notion of scalar-by-matrix differentiation, the derivatives of $G$ with respect
to $\mathbb{J}$ and $\hbox{det}(\mathbb{J})$ can be found \cite{Huang2015} as
\begin{align}
\frac{\partial G}{\partial\mathbb{J}} &=
\frac{d}{2}\hbox{det}(\mathbb{M}_{K})^{\frac{1}{2}}
(\hbox{tr}(\mathbb{J}\mathbb{M}_K^{-1}\mathbb{J}^T))^{\frac{3 d}{4}-1}\mathbb{M}_K^{-1}\mathbb{J}^T,
\label{partial-J}
\\
\frac{\partial G}{\partial \hbox{det}(\mathbb{J})} &=
\frac{1}{2} d^\frac{3 d}{4}\hbox{det}(\mathbb{M}_{K})^{-\frac{1}{4}}\hbox{det}(\mathbb{J})^{\frac{1}{2}}
\label{partial-detJ}.
\end{align}
With these formulas, we can rewrite \eqref{MM} as (cf. \cite{Huang2015})
\begin{equation}\label{xim}
\begin{split}
\frac{d\bm{\xi}_i}{d\varsigma}=
\frac{\hbox{det}(\mathbb{M}(\bm{x_i}))^{\frac{1}{2}} }{\tau}
\sum_{K\in\omega_i}|K|\bm{v}^K_{i_K},
\quad i=1,...,N_v
\end{split}
\end{equation}
where $\omega_i$ is the element patch associated with the vertex $\bm{x}_i$,
$i_K$ is the local index of $\bm{x}_i$ on $K$,
and $\bm{v}^K_{i_K}$ is the local velocity contributed by the element $K$ to the vertex $i_K$.
The local velocities $\bm{v}^K_{i_K},\; i_K=1,...,d$, are given by
\begin{equation}\label{vjk}
\begin{split}
\left[
  \begin{array}{c}
    ( \bm{v}_1^K  )^T  \\
    ( \bm{v}_2^K  )^T   \\
    \vdots\\
    ( \bm{v}_d^K  )^T   \\
   \end{array}
 \right]
=
- E_K^{-1}\frac{\partial G }{\partial \mathbb{J} }
- \frac{\partial G}{\partial \hbox{det}(\mathbb{J})}\frac{ \hbox{det}( E_{K_c} )}
{\hbox{det}(E_K)} E_{K_c}^{-1},
\quad
\bm{v}^K_0=-\sum_{i_K=1}^d\bm{v}^K_{i_K} .
\end{split}
\end{equation}
Note that the velocities for the boundary nodes need to be modified properly.
For example, the velocities for the corner vertices should be set to be zero.
For other boundary vertices,
the velocities should be modified such that they only slide along the boundary and
do not move out of the domain.

Starting with the reference computational mesh $\hat{\mathcal{T}}_c$ as the initial mesh,
the mesh equation \eqref{xim} is integrated over a physical time step for the case
with numerical solution of RTE (cf. \S\ref{sec:RTE-tests}) or from $\varsigma = 0$ to $\varsigma=1$
for DG-interpolation testing (cf. \S\ref{sec:interp-tests}).
The obtained new mesh is denoted by $\mathcal{T}_c^{new}$.
Note that $\mathcal{T}_h^{old}$ is kept fixed during the integration
and forms a correspondence with $\mathcal{T}_c^{new}$,
i.e., $\mathcal{T}_h^{old}=\Phi_h(\mathcal{T}_c^{new})$.
Then the new physical mesh $\mathcal{T}_h^{new}$ is defined
 as $\mathcal{T}_h^{new}=\Phi_h(\hat{\mathcal{T}}_c)$,
which can be computed using linear interpolation.

% section 4
\section{Numerical results for DG-interpolation}
\label{sec:interp-tests}

\begin{figure}[h]
\centering
\subfigure[Without PP limiter]{
\includegraphics[width=0.30\textwidth]{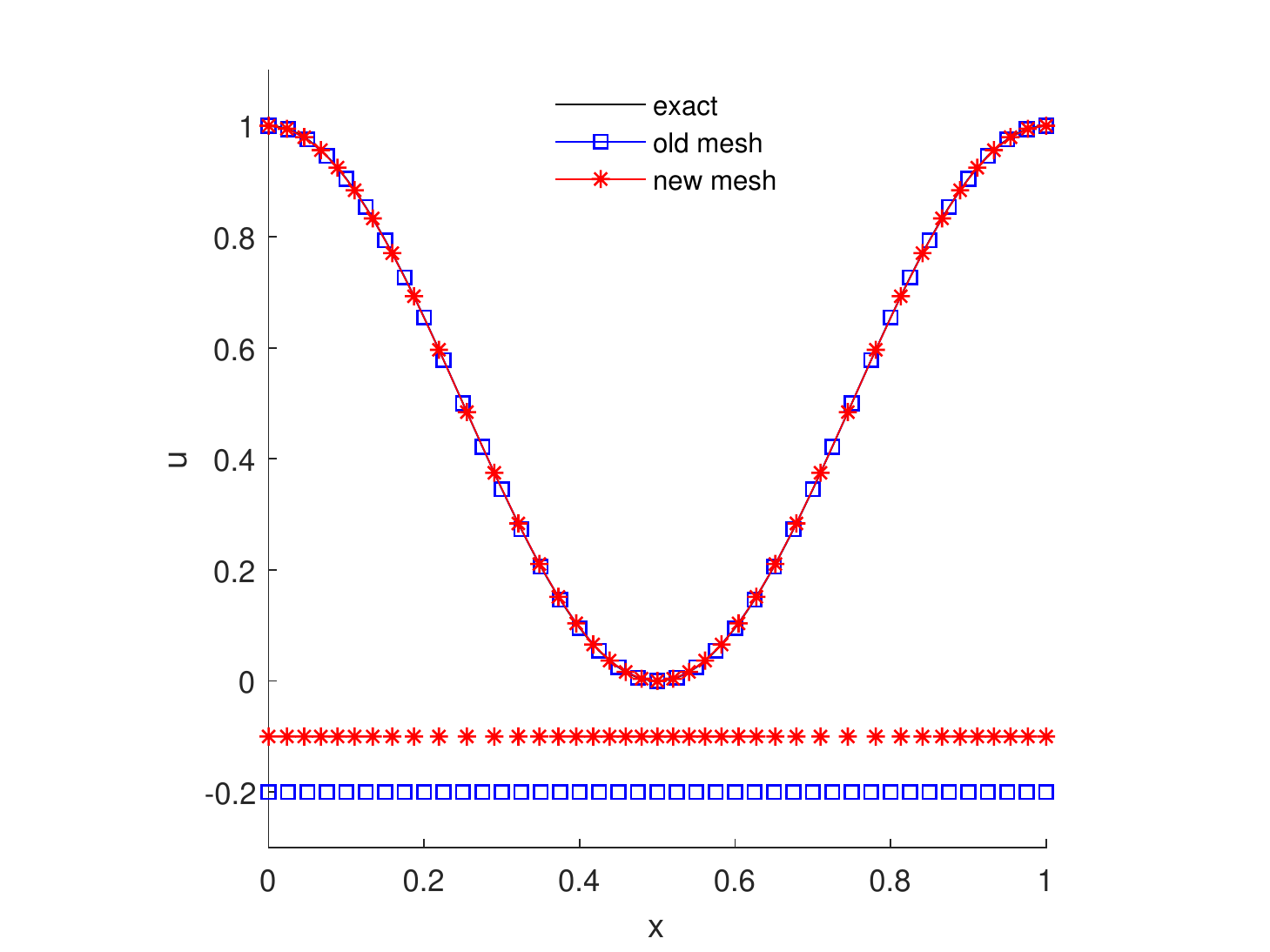}}
\subfigure[Close view of (a) near $u$=0]{
\includegraphics[width=0.30\textwidth]{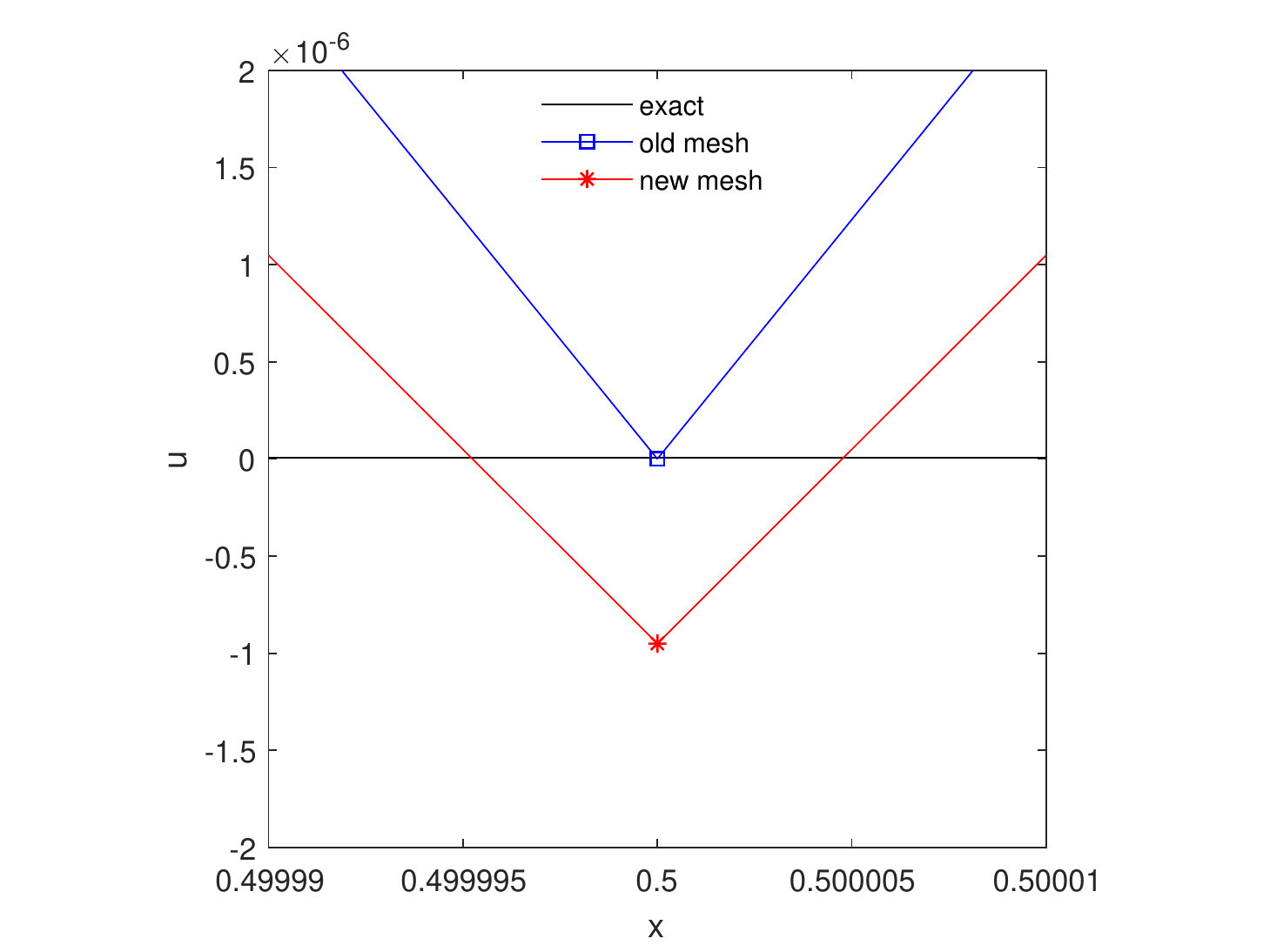}}
\\
\subfigure[With PP limiter]{
\includegraphics[width=0.30\textwidth]{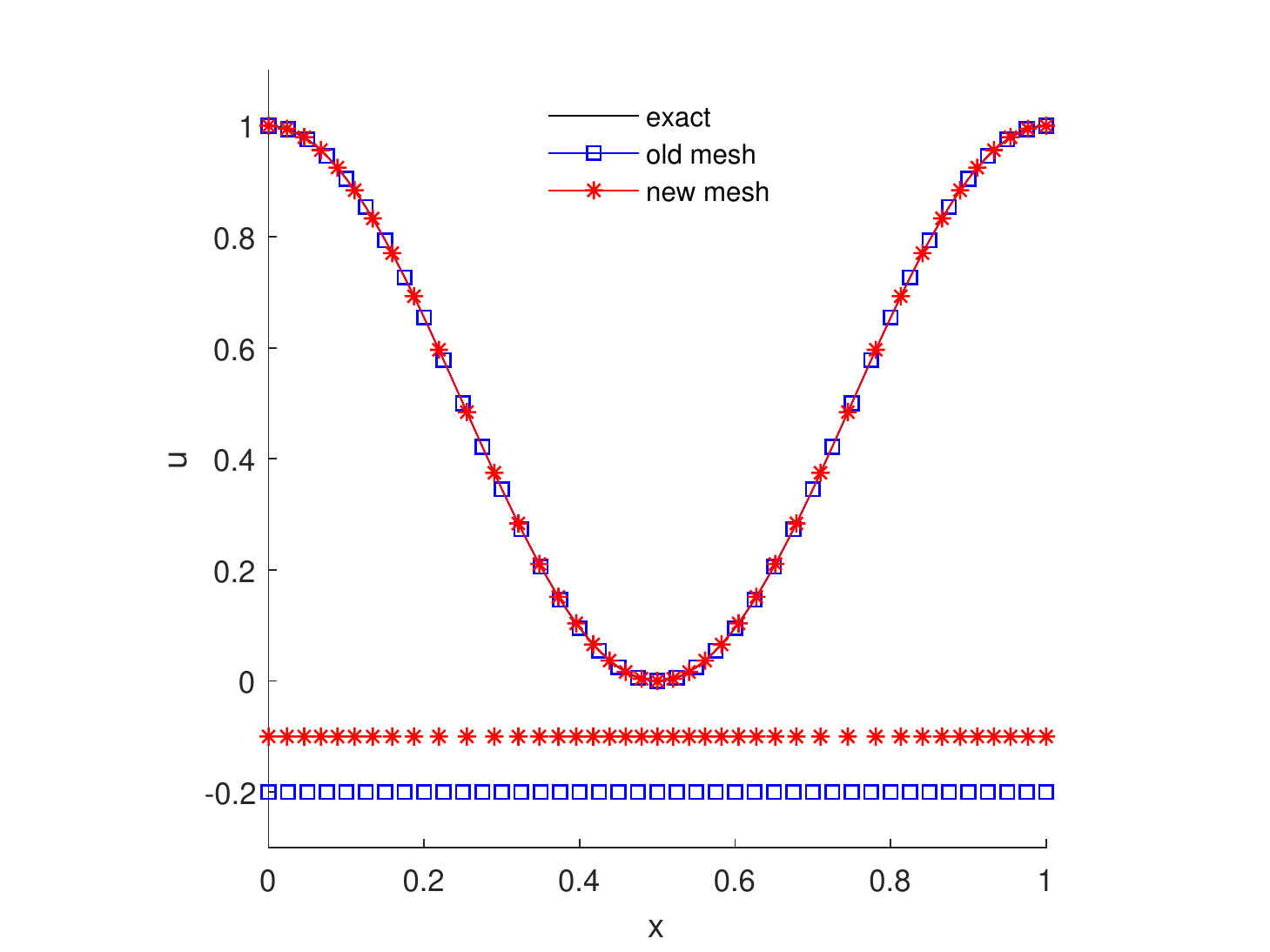}}
\subfigure[Close view of (c) near $u$=0]{
\includegraphics[width=0.30\textwidth]{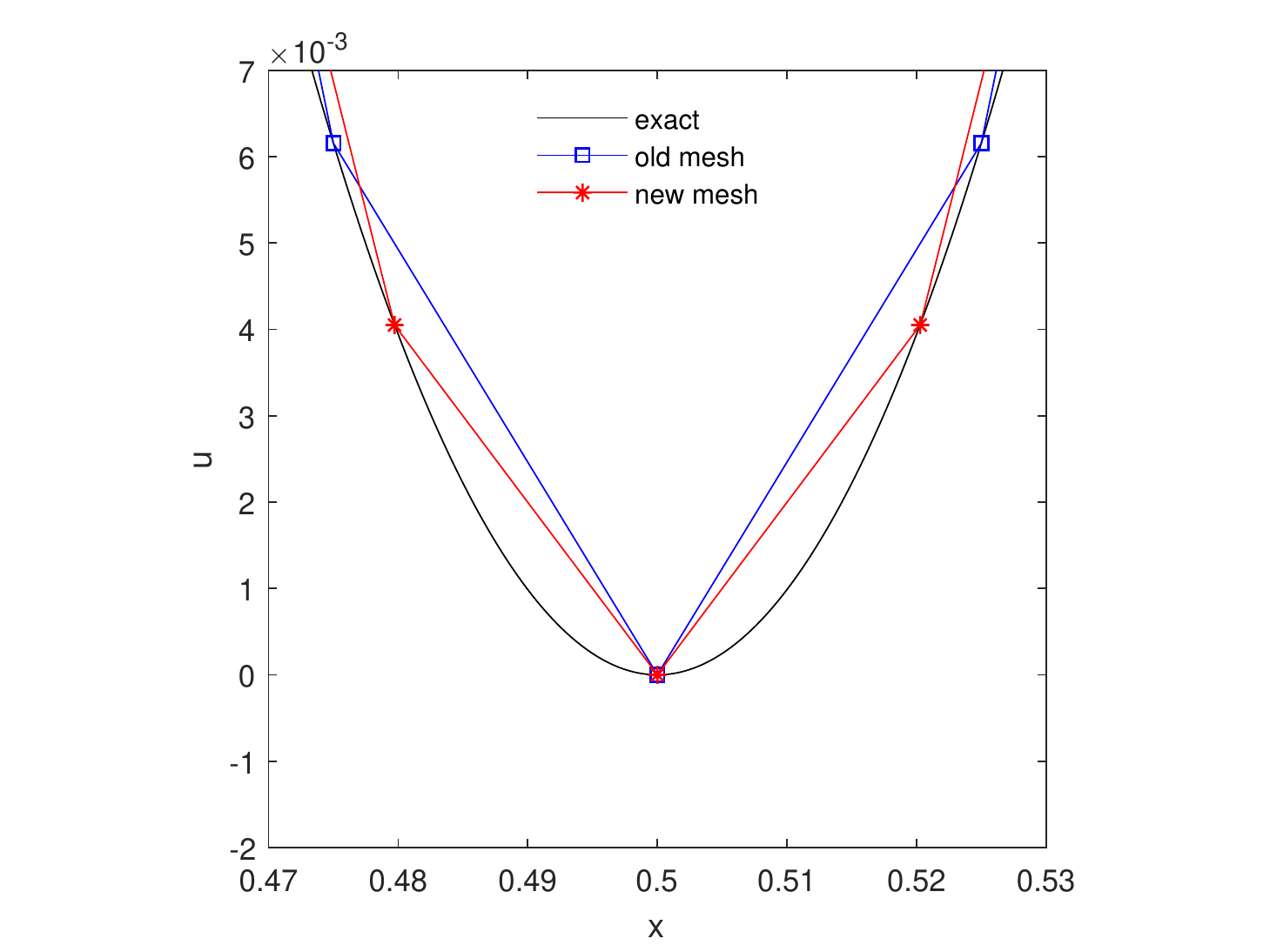}}
\caption{Example \ref{Ex1-interp-1d}. The meshes ($N=40$) and numerical solutions obtained by
$P^2$-DG interpolation with or without PP limiter.
}\label{Fig:ex1-ld-p2}
\end{figure}

\begin{figure}[h]
\centering
\subfigure[$L^1$ norm of error]{
\includegraphics[width=0.30\textwidth]{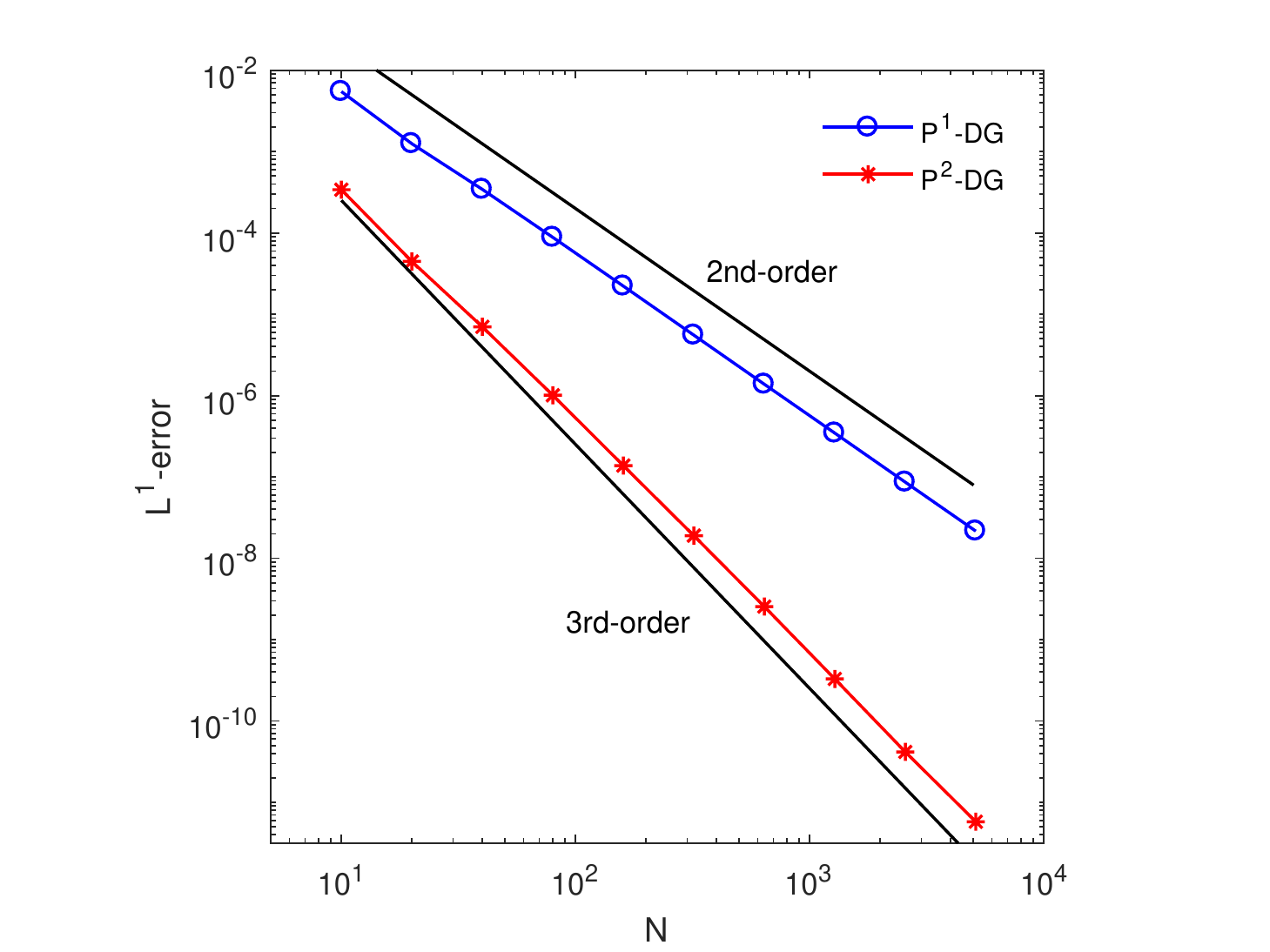}}
\subfigure[$L^{\infty}$ norm of error]{
\includegraphics[width=0.30\textwidth]{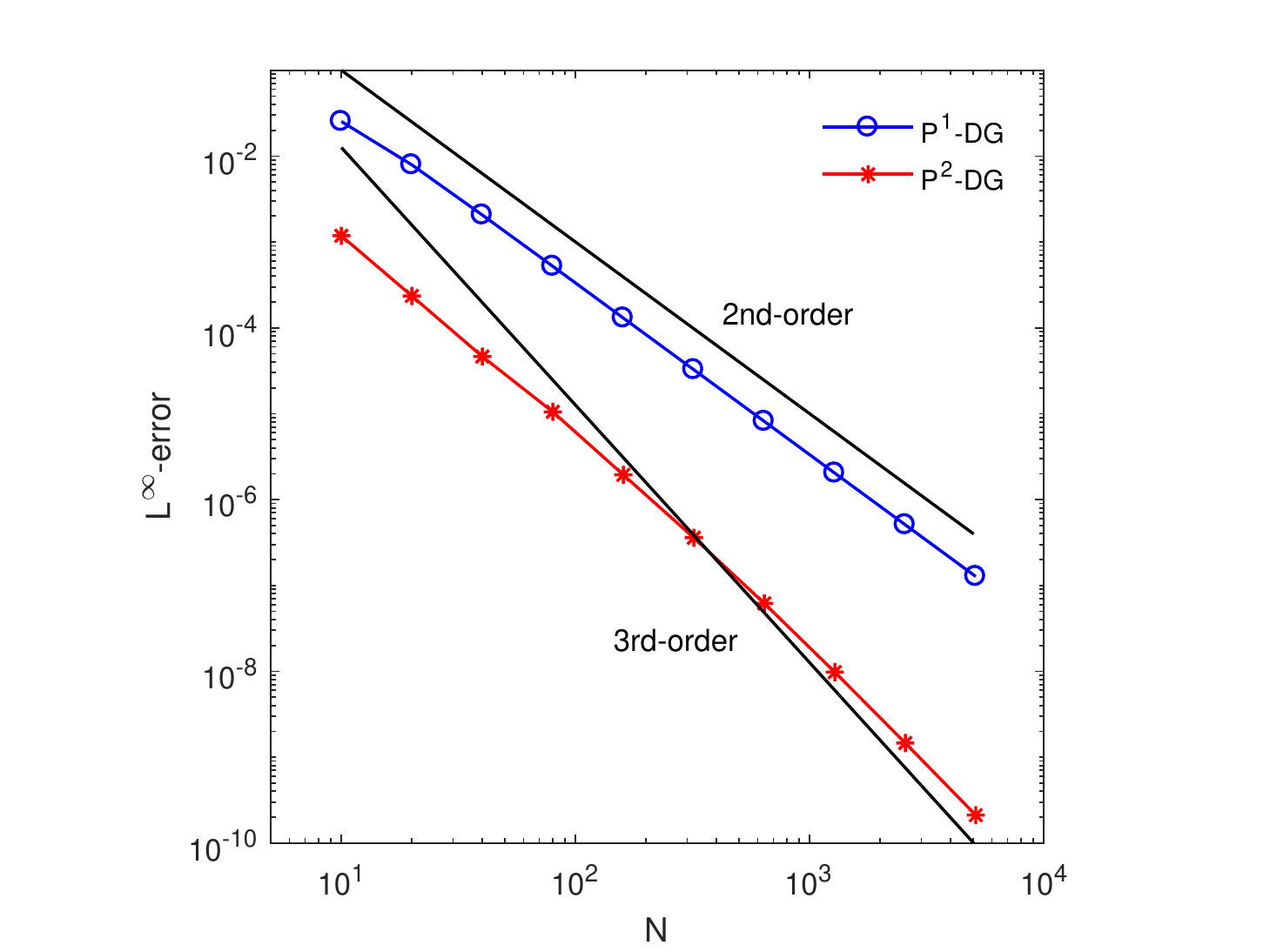}}
\subfigure[Number of time steps]{
\includegraphics[width=0.30\textwidth]{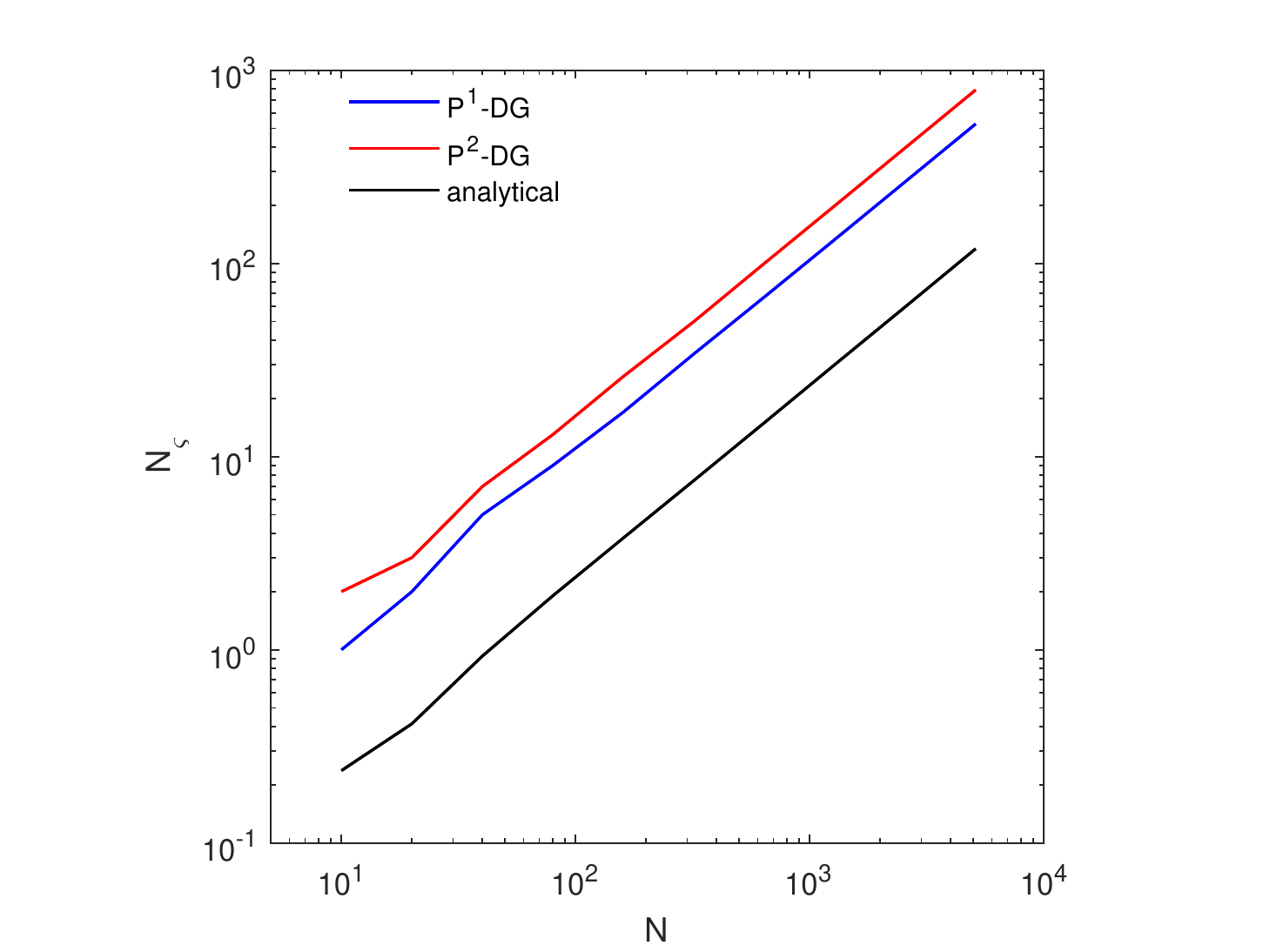}}
\caption{Example \ref{Ex1-interp-1d}. (a) and (b): The convergence history.
(c): The number of time steps used to reach $\varsigma=1$
is plotted against $N$ for the PP DG-interpolation.
The ``analytical" stands for the curve $N_\varsigma = N\max_{i}|\bm{x}^{old}_i - \bm{x}^{new}_i |$.
}\label{Fig:ex1-ld-order}
\end{figure}

In this section we present numerical results in one and two dimensions to demonstrate
the accuracy and positivity preservation property of the DG-interpolation scheme with PP limiter.
The CFL number for pseudo-time stepping taken
to be $1/4$ for $P^1$-DG and $1/6$ for $P^2$-DG in one dimension, and $1/4$ for $P^1$-DG and $1/9$ for $P^2$-DG in two dimensions.
In the computation, $\mathcal{T}^{old}_{h}$ is a given initial mesh and $\mathcal{T}^{new}_{h}$
is obtained using the MMPDE method described in the previous section.
More specifically, the metric tensor is first computed on the current mesh for a given function
(with its Hessian recovered based on nodal values via quadratic least squares fitting) and
then a new physical mesh is obtained through integrating the mesh equation \eqref{xim}
from $\varsigma = 0$ to $\varsigma=1$ with $\tau = 0.01$ and linear interpolation.
This procedure is iterated five times. No restriction is imposed on the mesh deformation.
% example
\begin{example}\label{Ex1-interp-1d}
In this test, we choose the function as
\begin{equation*}
u(x)=\cos^2(\pi x)+10^{-14},\qquad x\in(0,1) .
\end{equation*}
\end{example}
%Fig.~\ref{Fig:ex1-ld-p1} and
Fig.~\ref{Fig:ex1-ld-p2} shows the meshes and numerical solutions
obtained by $P^2$-DG interpolation with or without PP limiter
from the old mesh to the new one. It demonstrates that the PP limiter
is able to maintain the positivity of the solution.
The convergence history is plotted in Fig.~\ref{Fig:ex1-ld-order} (a,b),
which show that the PP DG-interpolation has the expected convergence order in both $L^1$ and $L^{\infty}$ norm.
Fig.~\ref{Fig:ex1-ld-order} (c) shows the number ($N_\varsigma$) of time steps used to
reach $\varsigma=1$ as $N$ increases for the PP DG-interpolation.
One can see that the curves for $P^1$-DG and $P^2$-DG are almost parallel to
$N_\varsigma = N\max_{i}|\bm{x}^{old}_i - \bm{x}^{new}_i |$, which is consistent with
the analysis for Case 2 in \S\ref{sec:Nsteps} (cf. \eqref{DG-Interp-cfl+1+3}).
For this example, $\max_{i}|\bm{x}^{old}_i - \bm{x}^{new}_i |$ stays almost constant
(about 0.023 for large $N$) as the mesh is being refined.

% example
\begin{example}\label{Ex2-interp-2d}
%(An accuracy and positivity-preserving test in two dimensions)
We consider
\begin{equation*}
u(x,y)=1- \tanh\Big{(}50\big{(}(x-0.5)^2+(y-0.5)^2-\frac{1}{16}\big{)}\Big{)}+10^{-14},
~ (x,y)\in(0,1)\times(0,1)
\end{equation*}
which has a sharp jump around the circle $(x-0.5)^2+(y-0.5)^2 = 1/16$.
\end{example}

\begin{figure}[h]
\centering
\subfigure[The $\mathcal{T}_h^{old}$ with $N$=1600]{
\includegraphics[width=0.30\textwidth]{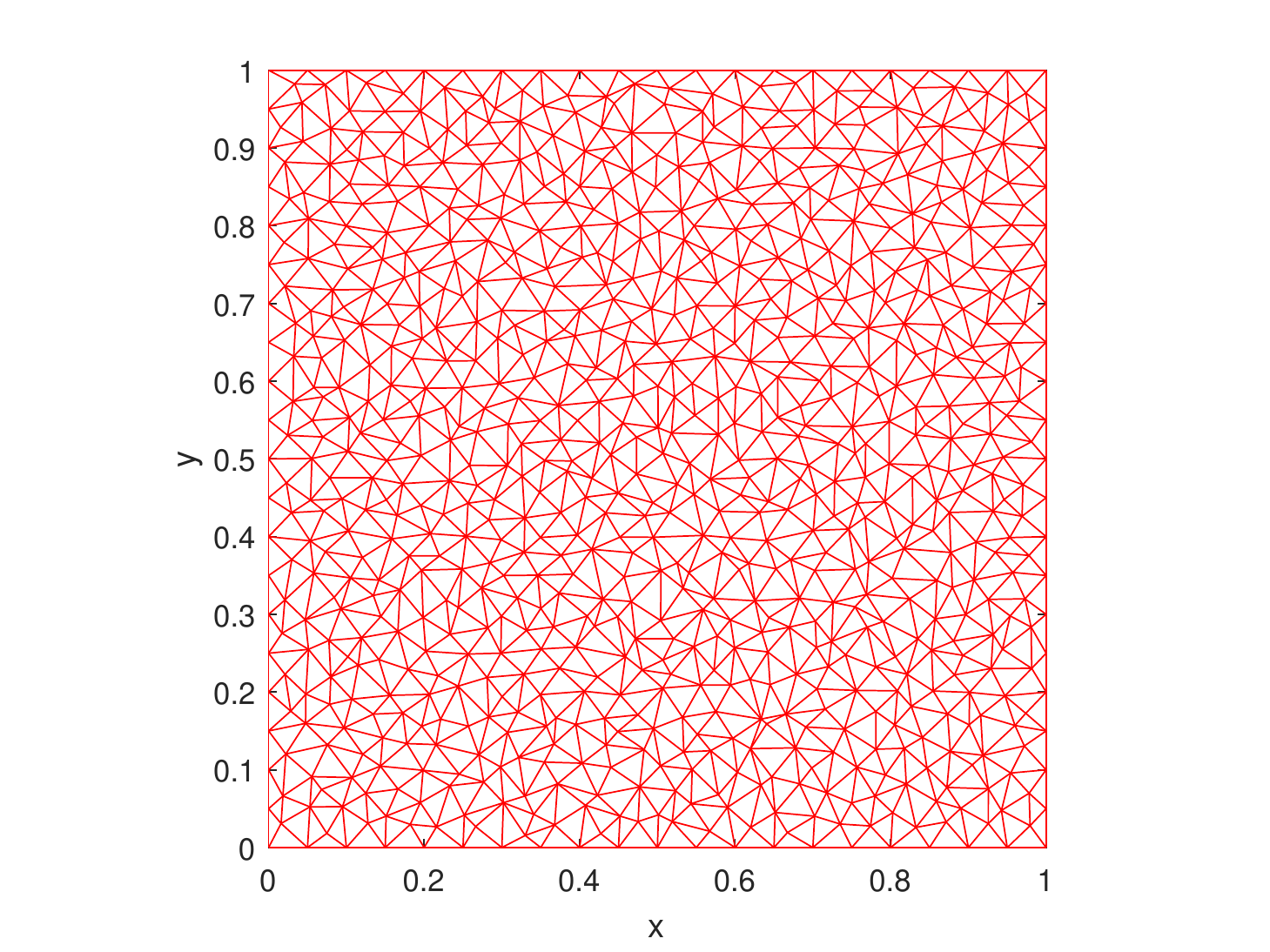}}
\subfigure[Solution contours on $\mathcal{T}_h^{old}$]{
\includegraphics[width=0.30\textwidth]{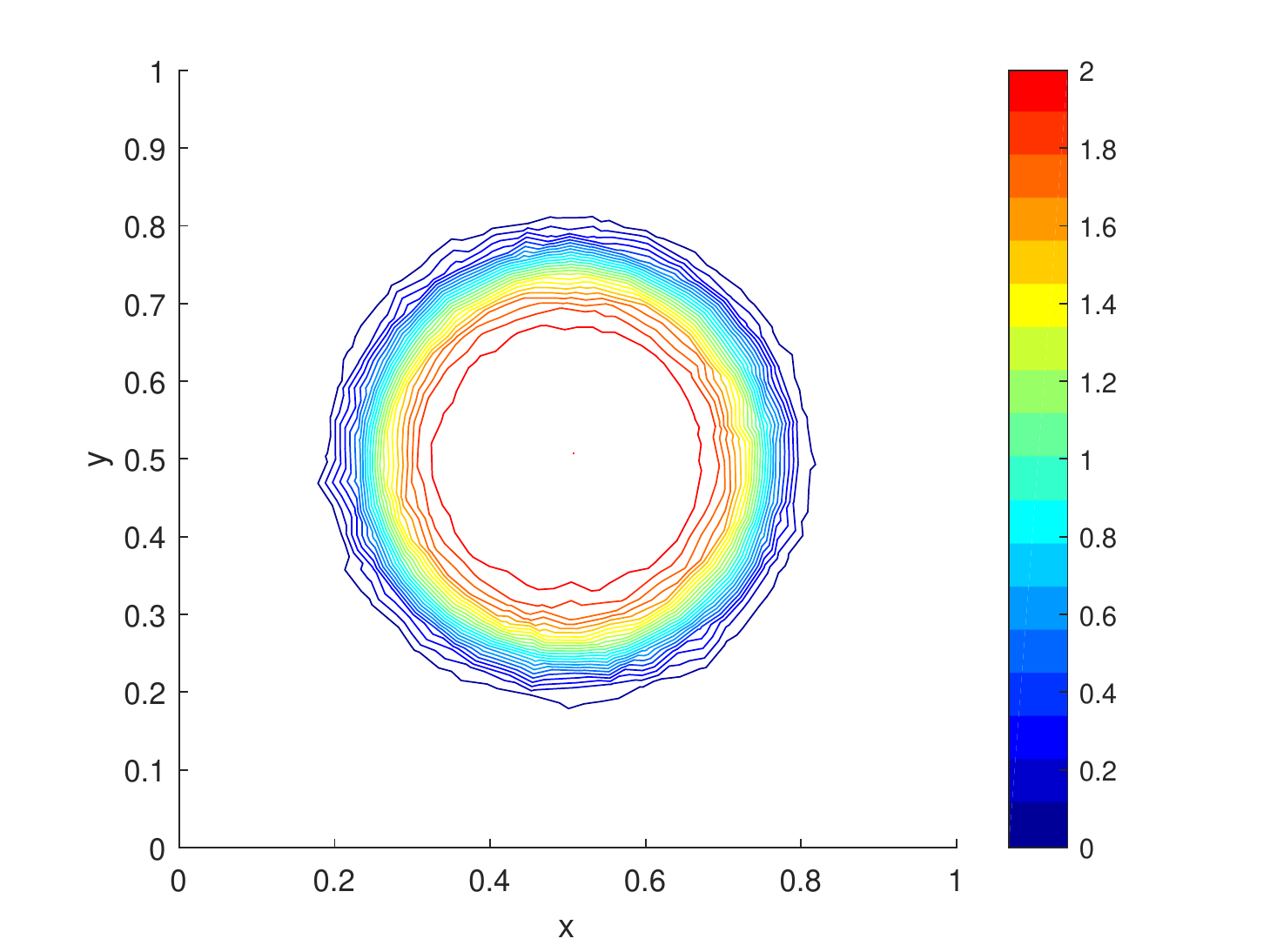}}
\\
\subfigure[The $\mathcal{T}_h^{new}$ with $N$=1600]{
\includegraphics[width=0.30\textwidth]{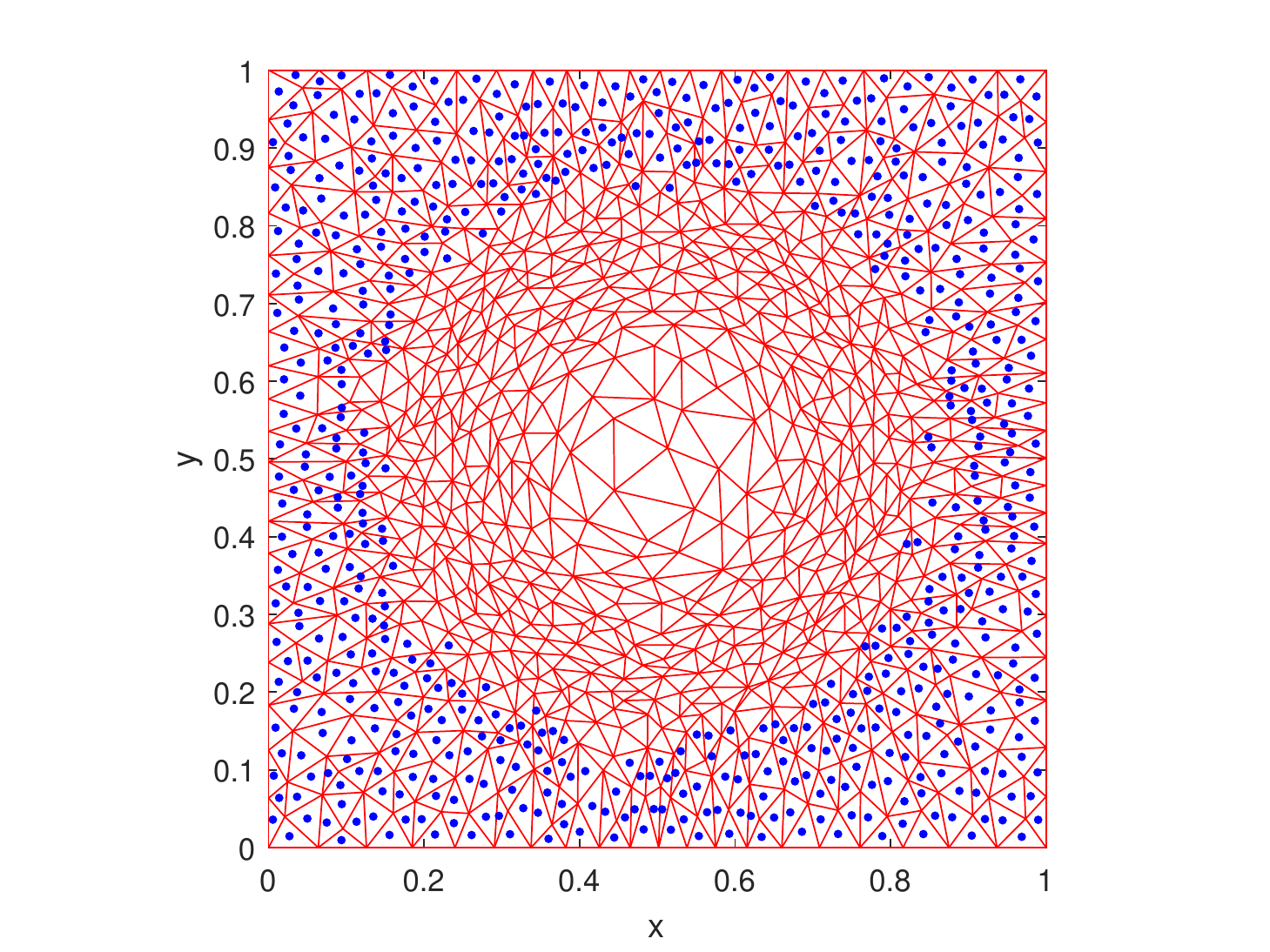}}
\subfigure[Solution contours on $\mathcal{T}_h^{new}$]{
\includegraphics[width=0.30\textwidth]{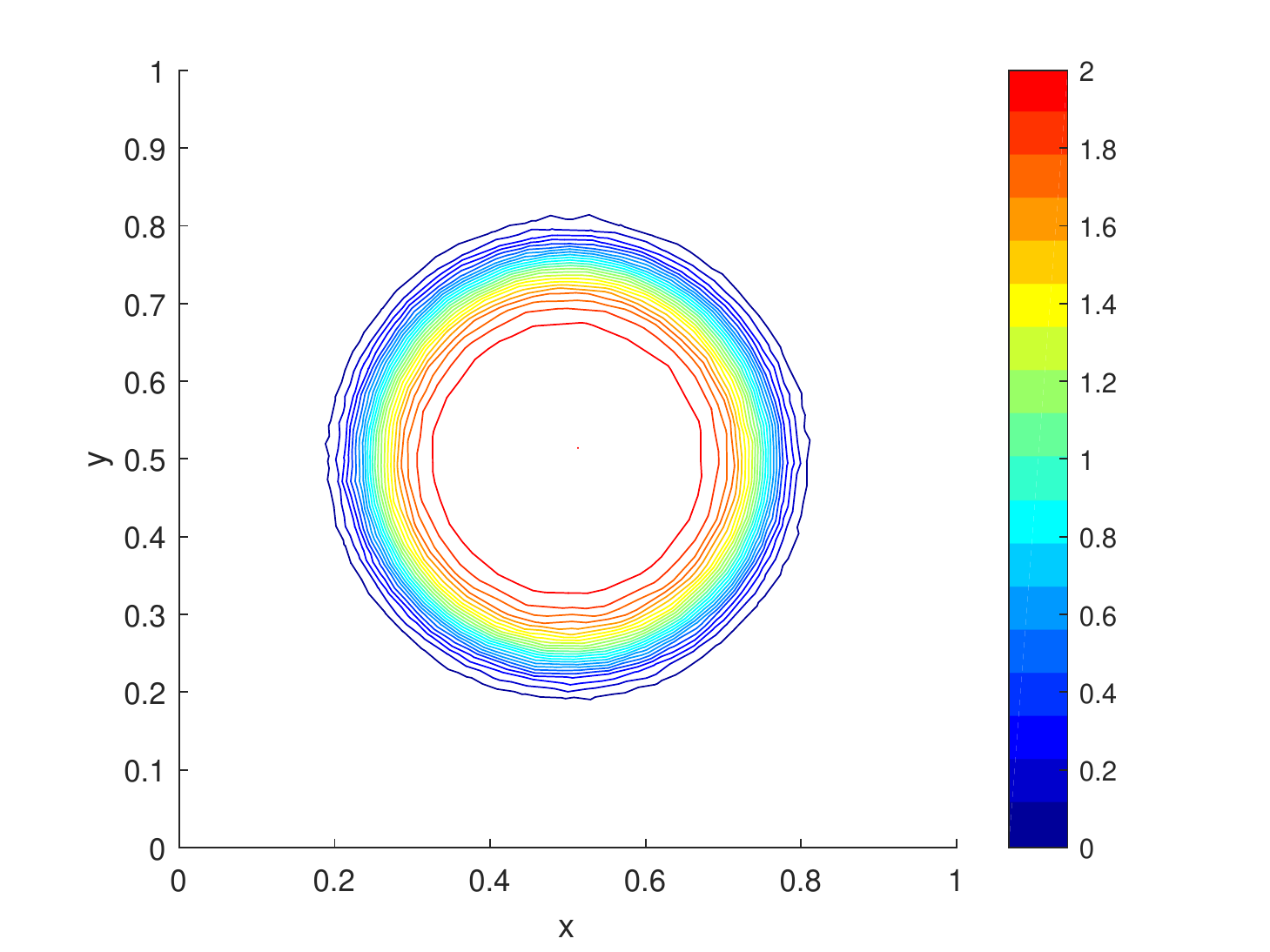}}
\caption{Example \ref{Ex2-interp-2d}.
The meshes of $N=1600$ and solution contours for $P^2$-DG interpolation with PP limiter.
The blue dots on the new mesh represent the cells where the PP limiter has been applied.
}\label{Fig:ex2-2d-p2}
\end{figure}

\begin{figure}[h]
\centering
\subfigure[$L^1$ norm of error]{
\includegraphics[width=0.30\textwidth]{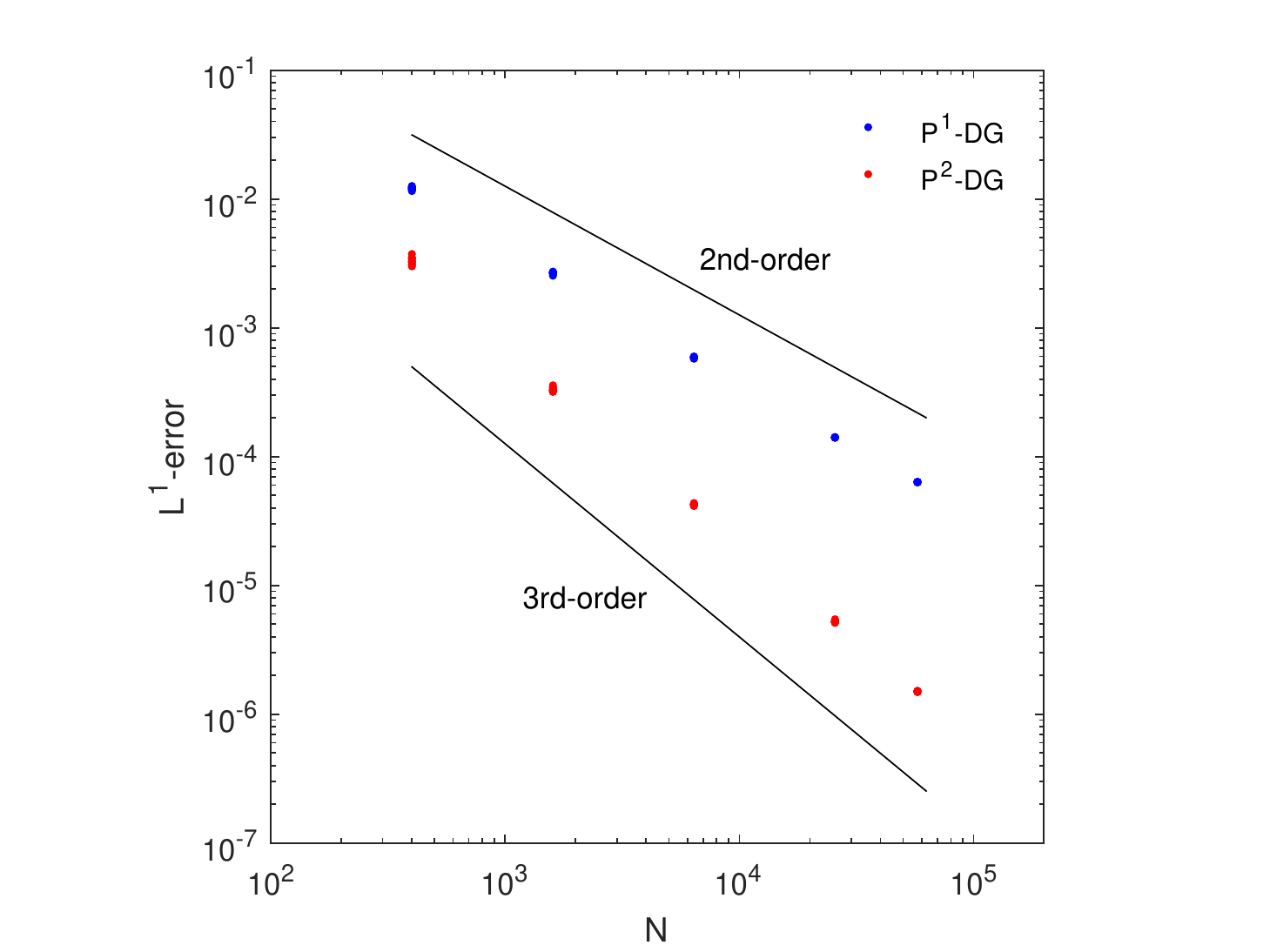}}
\subfigure[$L^{\infty}$ norm of error]{
\includegraphics[width=0.30\textwidth]{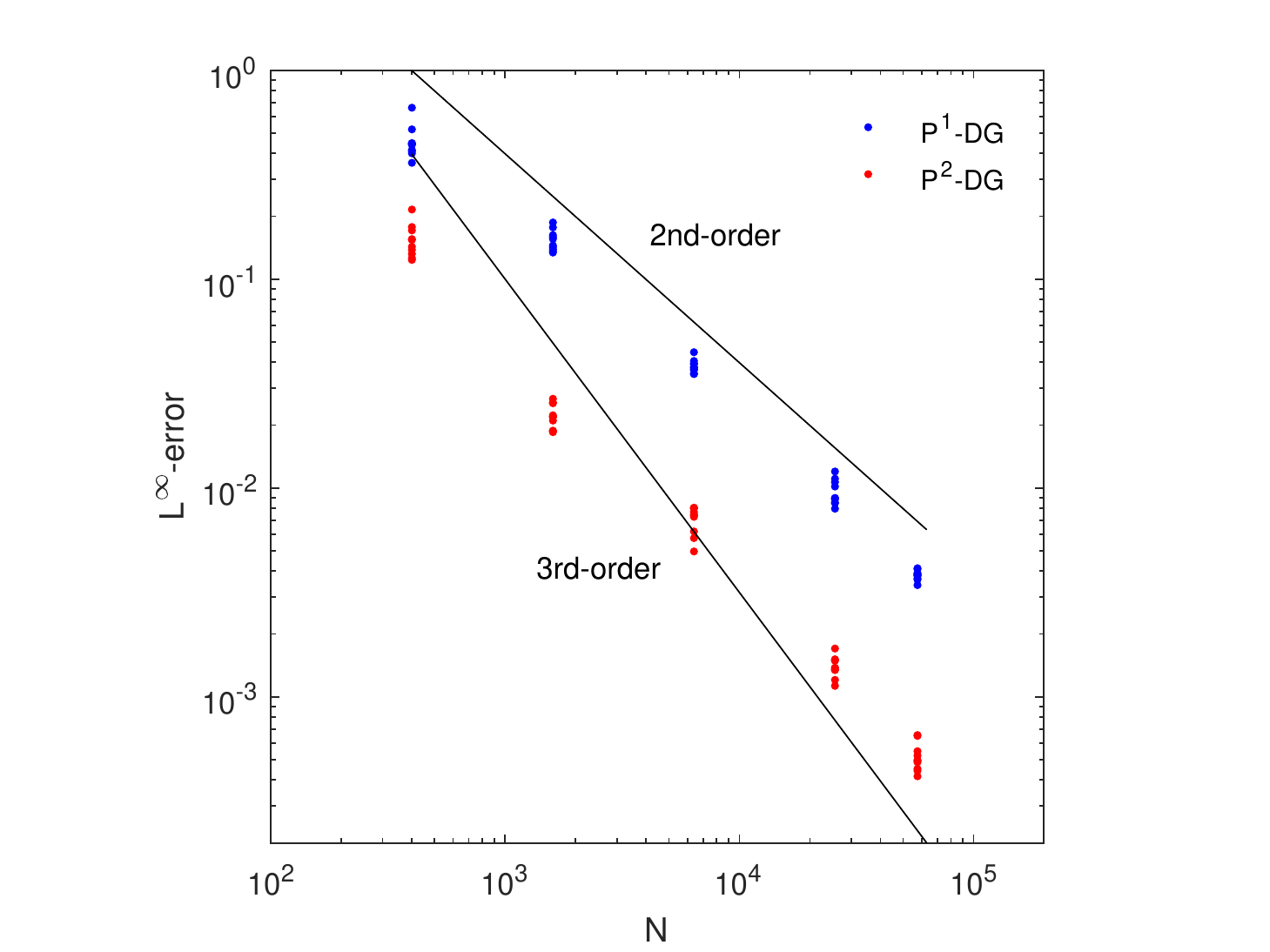}}
\subfigure[Number of time steps]{
\includegraphics[width=0.30\textwidth]{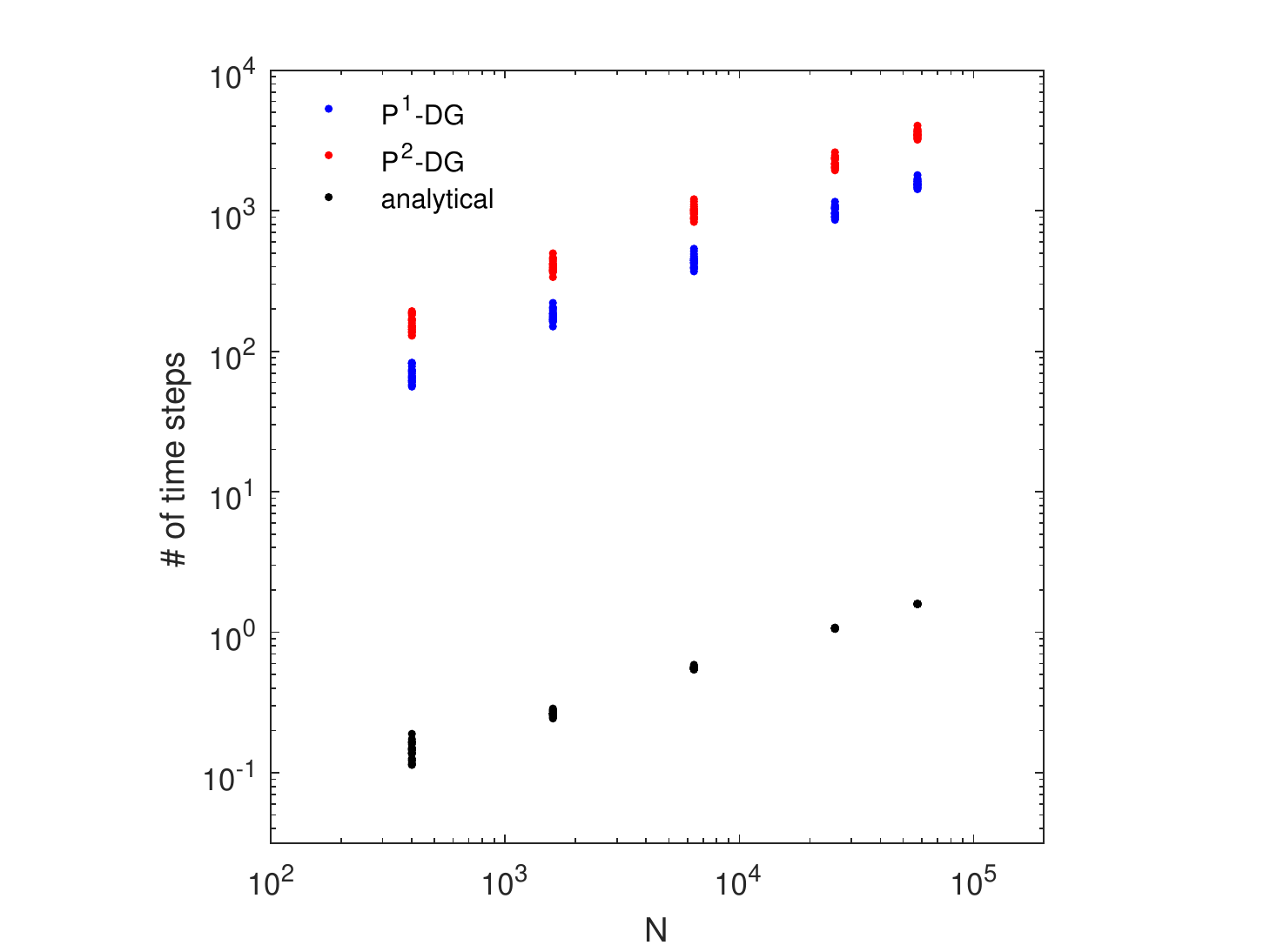}}
\caption{Example \ref{Ex2-interp-2d}. (a) and (b): The convergence history.
(c): The number of time steps used to reach $\varsigma=1$
is plotted against $N$ for the PP DG-interpolation.
The ``analytical" stands for the curve $N_\varsigma = \sqrt{N}\max_{i}|\bm{x}^{old}_i - \bm{x}^{new}_i |$.
}\label{Fig:ex2-2d-order}
\end{figure}

For this example, we start with a rectangular mesh, randomly perturb the interior vertices
by 40\% of the average element diameter, and obtain the final initial mesh
as the Delaunay mesh associated with the perturbed vertices. For each mesh resolution,
we carry out 20 runs. Fig.~\ref{Fig:ex2-2d-p2} shows a typical mesh of $N=1600$ (starting from
a $20\times 20$ rectangular mesh) and
corresponding solution contours obtained by PP $P^2$-DG interpolation.
The mesh elements where the solution becomes negative and the PP limiter
has been applied are indicated by blue dots.
The convergence history in $L^1$ and $L^\infty$ norm  is plotted in Fig.~\ref{Fig:ex2-2d-order} (a,b).
We can see that the error and $N_\varsigma$ have different values for
different runs for the same $N$ due to the randomness of initial meshes.
Moreover, initial meshes are nonsmooth, which leads to a less-than-third-order convergence
rate for the $L^\infty$ norm of the $P^2$-DG error. Nevertheless, the results show that
$P^1$-DG is second-order in both $L^1$ and $L^\infty$ norm and $P^2$-DG is third-order in $L^1$ norm.
Fig.~\ref{Fig:ex2-2d-order} (c) shows that the number of time steps to reach $\varsigma=1$
for $P^1$-DG and $P^2$-DG have a similar increase rate as
$N_\varsigma = \sqrt{N} \max_{i}|\bm{x}^{old}_i - \bm{x}^{new}_i |$,
verifying the estimate (\ref{DG-Interp-cfl+1+3}).
We note that $\max_{i}|\bm{x}^{old}_i - \bm{x}^{new}_i |$ stays almost constant (about 0.007)
for large $N$. This is large compared to the average element diameter $N^{-\frac{1}{2}}$
which decreases as $N$ increases, indicating a large deformation between the old and new meshes.
A larger number of pseudo-time steps is required for larger mesh deformation.

% section 5
\section{Application of DG-interpolation to MM-DG simulation of RTE}
\label{sec:DOM-MM-DG-RTE}

In this section, as an application, we consider the use of the DG-interpolation
in a rezoning MM-DG method for the numerical solution of RTE in one and two spatial dimensions.
Our goal is to show that the method maintains high-order accuracy of DG schemes
while preserving the positivity of the radiative intensity.

\begin{figure}[h]
\centering
\tikzset{my node/.code=\ifpgfmatrix\else\tikzset{matrix of nodes}\fi}
{\footnotesize
\begin{tikzpicture}[every node/.style={my node},scale=0.45]
\draw[thick] (0,0) rectangle (4,5.5);
\draw[thick] (6,0) rectangle (12,5.5);
\draw[thick] (14,0) rectangle (20,5.5);
\draw[thick] (22,0) rectangle (28,5.5);
\draw[->,thick] (-2 ,2.75)--(0 ,2.75);
\draw[->,thick] (4 ,2.75)--(6 ,2.75);
\draw[->,thick] (12,2.75)--(14,2.75);
\draw[->,thick] (20,2.75)--(22,2.75);
\draw[->,thick] (28,2.75)--(30,2.75);
%\node[above] at (15,-2) {next time step\\};
\node (node1) at (2,2.75) {Given \\ $\mathcal{T}_h^n$, $u^n$ \\};
\node (node2) at (9,2.75) {Generate \\$\mathcal{T}_h^{n+1}$\\ using MMPDE \\ moving mesh\\ method\\};
\node (node3) at (17,2.75) {DG-interpolate\\ $u^n$ (on $\mathcal{T}_h^n$)\\ $\rightarrow$ $\tilde{u}^n$ (on $\mathcal{T}_h^{n+1}$)\\};
\node (node4) at (25,2.75) {Solve RTE on \\ fixed mesh $\mathcal{T}_h^{n+1}$: \\
                            $\tilde{u}^n\rightarrow u^{n+1}$\\};
\end{tikzpicture}
}
\caption{Illustration of the rezoning moving mesh method.}
\label{fig:Ibvp-solver}
\end{figure}
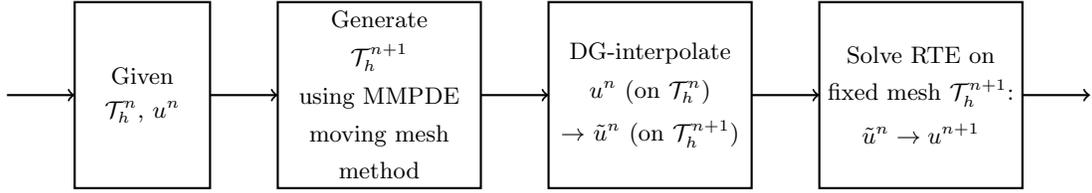

The rezoning MM method is illustrated in Fig.~\ref{fig:Ibvp-solver}.
As one can see, it involves three independent steps, generating the new mesh, interpolating
the solution from the old mesh to the new one, and solving the RTE on the new mesh.
In this work, we use the MMPDE method described in \S\ref{sec:mmpde} to
generate the new mesh, the DG-interpolation scheme of \S\ref{sec:DG-interpolation}
to interpolate the physical variables between the old and new meshes,
and a high-order PP DG scheme of \cite{LingChengShu2018,YuanChengShu2016,ZhangChengQiu2019}
to solve the RTE on the new mesh $\mathcal{T}_h^{n+1}$. Since the first two steps have been
discussed in previous sections, we focus on the last step in this section.
The RTE is an integro-differential equation modeling the conservation of photons \cite{gc}.
We consider a case with isotropically scattering radiative transfer. The governing equation
for this case reads as
\begin{equation}
\label{pde-RTE}
\frac{1}{c}\frac{\partial I}{\partial t}
+\Omega \cdot \nabla I+\sigma_tI
=\frac{\sigma_s}{4\pi}\int_S I(\bm{x},\tilde{\Omega},t)d\tilde{\Omega}+q,
\end{equation}
where $c$ is the speed of photons, $\bm{x}$ is the spatial variable,
$\nabla$ is the gradient operator with respect to $\bm{x}$,
$\Omega$ is the unit angular variable,
$S$ is the unit sphere, $t$ is time,
$I(\bm{x},\Omega,t)$ is the radiative intensity in the direction $\Omega$,
$\sigma_s \geq 0 $ is the scattering coefficient of the medium,
$ \sigma_t$ is the extinction coefficient of the medium which includes absorption and scattering,
and $ q(\bm{x}, \Omega,t)$ is a given source term.
The vector $\bm{x}$ is described by the Cartesian coordinates $x, y, z$
while $\Omega=(\zeta,\eta,\mu)$ is usually described by a polar angle $\beta$ measured with respect to the $z$ axis
and a corresponding azimuthal angle $\varphi$.
Denoting $\mu = \cos\beta$, $\zeta = \sin\beta \cos\varphi$, $\eta =\sin\beta \sin\varphi$
then
\begin{equation*}
 d\bm{x} = dxdydz, \quad d\Omega = \sin\beta d\beta d\varphi = -d\mu d\varphi.
\end{equation*}
In this work we consider the numerical solution of \eqref{pde-RTE} in one and two spatial dimensions.

\subsection{Positivity-preserving DG scheme for RTE in one dimension}
\label{PPDG-RTE-1d}
The one-dimensional form of \eqref{pde-RTE} reads as
\begin{equation}
\label{rte-1d-pde}
\frac{1}{c}\frac{\partial I}{\partial t}
+\mu\frac{\partial I}{\partial x}+\sigma_tI
=\frac{\sigma_s}{2}\int_{-1}^{1}I(x,\tilde{\mu},t)d\tilde{\mu}+q,
\end{equation}
where $x \in (a, b)$, $\mu \in (-1, 1)$, and $t \in (0, T]$.
The initial and boundary conditions are
\begin{align*}
%\label{rte-1d-initial}
& I(x,\mu,0)=I_0(x,\mu),\quad  x\in (a, b)
\\
%\label{rte-1d-boundy}
& \begin{cases}
I(a,\mu,t)=I_l(\mu,t),~~~ \quad 0 <\mu \leq 1 ,~0< t\leq T
\\
I(b,\mu,t)=I_r(\mu,t), \quad -1 \leq \mu < 0 ,~0< t\leq T.
\end{cases}
\end{align*}

We first use the discrete ordinate method (DOM) \cite{DOM} to discretize
\eqref{rte-1d-pde} in the angular variable. Consider a Gauss-Legendre quadrature rule
with weights $w_m$ and nodes $\mu_m$, $m = 1, ..., N_a$.
We define the discrete-ordinate approximation for RTE as
\begin{equation}
\label{rte-1d-dom}
\frac{1}{c}\frac{\partial I_m}{\partial t}
+\mu_m\frac{\partial I_m}{\partial x}
+\sigma_tI_m
=\sigma_s\sum_{m'=1}^{N_a}w_{m'}I_{m'}+q_m,\quad
m=1,..., N_a
\end{equation}
where $I_m = I_m(x,t) \approx I(x,\mu_m, t)$.

For temporal discretization, if we use an explicit scheme, we will have to take a small time step as
$\mathcal{O}(1/c)$ to ensure stability. To avoid this, we use the backward Euler scheme and have
\begin{equation}
\label{rte-1d-dom-dt}
\tilde{\sigma}_t I_m^{n+1}
+\mu_m\frac{\partial I_m^{n+1}}{\partial x}
=\sigma_s\sum_{m'=1}^{N_a}w_{m'}I_{m'}^{n+1}
+\tilde{q}_m^{n+1},\quad m=1,..., N_a
\end{equation}
where $I_m^{n+1} \approx I_m(x,t_{n+1})$ and
\begin{equation}
\tilde{\sigma}_t = \sigma_t + \frac{1}{c\Delta t},
\quad
\tilde{q}_m^{n+1} = q_m^{n+1}+ \frac{1}{c\Delta t}I_m^{n},
\quad
\Delta t=t^{n+1}-t^{n}.
\label{sigma-1}
\end{equation}

We now consider the DG spatial discretization for \eqref{rte-1d-dom-dt} on $\mathcal{T}^{n+1}_h$.
We only consider here for the case with $\mu_m > 0$, as a similar procedure can be used for $\mu_m < 0$.
Assume that the cells of $\mathcal{T}^{n+1}_h$ can be written as
\[
K^{n+1}_i=(x^{n+1}_{i-1/2}, x^{n+1}_{i+1/2}),\quad i=1,...,N.
\]
Multiplying (\ref{rte-1d-dom-dt}) with a test function, integrating the resulting
equation over $K^{n+1}_i$, taking integration by part for the second term,
and applying the upwind numerical flux at the cell boundaries, we obtain
the DG formulation as: find $I^{n+1}_{m}\in V^{r}_h(t_{n+1})$ such that,
for $\forall \phi\in P^{r}(K^{n+1}_i)$, $i = 1, ..., N$,
\begin{align}
\label{rte-1d-DG}
&\int_{K^{n+1}_i}\tilde{\sigma}_tI^{n+1}_{m,i}\phi dx
-\mu_m\int_{K^{n+1}_i}I^{n+1}_{m,i}\phi'dx
+\mu_m \big(I^{n+1}_{m,i}\phi\big) |_{x^{n+1}_{i+1/2}}
\\
&=\int_{K^{n+1}_i}\sigma_s\Psi^{n+1}_{i}\phi dx
+\int_{K^{n+1}_i}\tilde{q}^{n+1}_{m,i}\phi dx
+\mu_m \big(I^{n+1}_{m,i-1}\phi\big) |_{x^{n+1}_{i-1/2}},
\notag
\end{align}
where $I^{n+1}_{m,i} = I_{m}^{n+1}(x)|_{K^{n+1}_i}$,
$\tilde{I}^{n}_{m}$ is the DG-interpolant of $I_m^{n}$
from $\mathcal{T}^{n}_h$ to $\mathcal{T}^{n+1}_h$, and
\[
\Psi^{n+1}_{i} = \sum_{m=1}^{N_a}w_{m}I^{n+1}_{m,i},
\quad
\tilde{q}^{n+1}_{m,i} = q^{n+1}_{m,i}+ \frac{1}{c\Delta t}\tilde{I}_{m,i}^{n}.
\]

\begin{figure}[h]
\centering
\begin{tikzpicture}[scale = 0.8]
\draw [-,thick] (-4,1)--(-3.3,1)-- (-2.2,1)-- (-1,1) -- (0.5, 1)-- (1.7, 1)-- (3, 1)--(4,1);
\draw [left] (-4.3,1.3) node {$\mu_m>0$:};
\draw [fill] (-4,1) circle (0.05);
\draw [fill] (-3.3,1) circle (0.05);
\draw [fill] (-2.2,1) circle (0.05);
\draw [fill] (-1,1) circle (0.05);
\draw [fill] (0.5,1) circle (0.05);
\draw [fill] (1.7,1) circle (0.05);
\draw [fill] (3,1) circle (0.05);
\draw [fill] (4,1) circle (0.05);
\draw [->,thick,red ] (-3.9,1.3)--(-3.4,1.3);
\draw [->,thick,red ] (-3.1,1.3)--(-2.4,1.3);
\draw [->,thick,red ] (-2.0,1.3)--(-1.2,1.3);
\draw [->,thick,red ] (-0.8,1.3)--( 0.3,1.3);
\draw [->,thick,red ] ( 0.7,1.3)--( 1.5,1.3);
\draw [->,thick,red ] ( 1.9,1.3)--( 2.8,1.3);
\draw [->,thick,red ] ( 3.2,1.3)--( 3.8,1.3);
%\draw [below] (-0.25,1) node {$K^{n+1}_{i}$};
%
\draw [-,thick] (-4,0)--(-3.3,0)-- (-2.2,0)-- (-1,0) -- (0.5, 0)-- (1.7, 0)-- (3, 0)--(4,0);
\draw [left] (-4.3,0.3) node {$\mu_m<0$:};
\draw [fill] (-4,0) circle (0.05);
\draw [fill] (-3.3,0) circle (0.05);
\draw [fill] (-2.2,0) circle (0.05);
\draw [fill] (-1,0) circle (0.05);
\draw [fill] (0.5,0) circle (0.05);
\draw [fill] (1.7,0) circle (0.05);
\draw [fill] (3,0) circle (0.05);
\draw [fill] (4,0) circle (0.05);
\draw [<-,thick,red ] (-3.9,0.3)--(-3.4,0.3);
\draw [<-,thick,red ] (-3.1,0.3)--(-2.4,0.3);
\draw [<-,thick,red ] (-2.0,0.3)--(-1.2,0.3);
\draw [<-,thick,red ] (-0.8,0.3)--( 0.3,0.3);
\draw [<-,thick,red ] ( 0.7,0.3)--( 1.5,0.3);
\draw [<-,thick,red ] ( 1.9,0.3)--( 2.8,0.3);
\draw [<-,thick,red ] ( 3.2,0.3)--( 3.8,0.3);
\draw [below] (-0.25,0) node {$K^{n+1}_{i}$};
\end{tikzpicture}
\caption{Mesh sweeping directions for $\mu_m>0$ (top) and $\mu_m<0$ (bottom).}
\label{Fig:sweep-direction}
\end{figure}
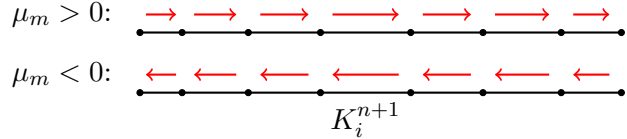

Notice that the unknown variables in different angular directions are coupled in \eqref{rte-1d-DG}
through $\Psi^{n+1}_{i}$. The so-called source iteration (SI) \cite{SI}
is commonly employed to solve the equations separately.
Denote the $\ell$-th iterate of the solution by $I_{m,i}^{n+1,(\ell)}$. Then the scheme reads as
\begin{align}
\label{rte-1d-DG-SI}
&\int_{K^{n+1}_i}I^{n+1,(\ell+1)}_{m,i}
\big{(}\tilde{\sigma}_t\phi -\mu_m\phi'\big{)} dx
+\mu_m \big(I^{n+1,(\ell+1)}_{m,i}\phi \big)|_{x^{n+1}_{i+1/2}}
 \\
&=\int_{K^{n+1}_i}
\big{(}\sigma_s\Psi^{n+1,(*)}_{i}+\tilde{q}^{n+1}_{m,i}\big{)}\phi dx
+\mu_m \big(I^{n+1,(\ell+1)}_{m,i-1}\phi\big) |_{x^{n+1}_{i-1/2}},
~ \forall \phi\in P^{r}(K^{n+1}_i)
\notag
\end{align}
where $\Psi^{n+1,(*)}_{i}= \sum_{m=1}^{N_a}w_{m}I^{n+1,(*)}_{m,i}$
and $I^{n+1,(*)}_{m,i}$ is taken as $I^{n+1,(\ell+1)}_{m,i}$
when it is available and otherwise as $I^{n+1,(\ell)}_{m,i}$.
The sweeping direction in space is indicated in Fig.~\ref{Fig:sweep-direction}.
The iteration is stopped when the maximum norm of the difference between
any two consecutive iterates is smaller than $10^{-12}$.
The radiative intensity is positive in physics. However, a numerical approximation
may contain negative values especially for high-order methods.
The appearance of spurious negative values could lead
to instability in the computation and slow iterative convergence.
Thus, it is important to develop schemes that preserve the positivity of the radiative intensity.
To this end, we mention that  it has been proved in \cite{LingChengShu2018} any $P^r$-DG scheme
(including the one described above) produces the positive cell averages for the one-dimensional RTE
on fixed meshes provided that both the inflow boundary condition from the upstream cell (including
the physical boundary condition for the first cell) and the source term are positive
and the initial condition is nonnegative.
As a consequence, the linear scaling PP limiter \cite{Liu-Osher1996,ZhangShu2010}
can be used to preserve the positivity of the radiative intensity.
The reader is referred to \S\ref{positivity} and \cite{LingChengShu2018,YuanChengShu2016} for detail.

With the positivity preserving property of the DG scheme (\ref{rte-1d-DG-SI}) and that of the PP DG-interpolation,
we can claim (cf. Fig.~\ref{fig:Ibvp-solver}) that the rezoning MM-DG method preserves the nonnegativity
of the radiative intensity.

% section 5.2
\subsection{Positivity-preserving DG scheme for RTE on triangular meshes}
\label{PPDG-RTE-2d}
The two-dimensional form of \eqref{pde-RTE} reads as
\begin{equation}
\label{rte-2d-pde}
\frac{1}{c}\frac{\partial I}{\partial t}
+\Omega \cdot \nabla I+\sigma_tI
= \frac{\sigma_s}{4\pi}\int_S I(x,y,\tilde{\Omega},t)d\tilde{\Omega} + q,
\end{equation}
where $(x,y)\in \mathcal{D}$, $t\in(0,T]$, $\Omega = (\zeta,\eta)$, and
\[
\zeta =\sqrt{1-\mu^2} \cos\varphi \in (-1, 1),
~ \eta =\sqrt{1-\mu^2} \sin\varphi \in (-1, 1), \quad \mu \in (-1, 1),
~ \varphi \in (0, 2\pi).
\]
The initial and inflow boundary conditions are
\begin{align*}
%\label{s3.2i}
& I(x,y,\Omega,0) = I_0(x,y,\Omega), \quad (x,y)\in \mathcal{D}, \quad \Omega \in S%[-1, 1]\times [-1, 1]
\\
%\label{s3.2b}
& I(x,y,\Omega,t) = I_b(x,y,\Omega,t) ,\quad (x,y)\in \partial \mathcal{D}_{in},\quad \Omega \in S, \quad t\in(0,T] .
\end{align*}
Here, $I_0(x,y,\Omega)$ and $I_b(x,y,\Omega,t)$ are given functions, $\partial \mathcal{D}_{in} = \{(x,y)\in \partial\mathcal{D} \; |\; \bm{n}(x,y) \cdot\Omega<0\}$, and
$\bm{n}(x,y)$ is the unit outward normal vector of the boundary.
It is worth pointing out that no boundary condition is needed in $\Omega$-direction.

Once again, we use the DOM for the discretization in $\Omega$.
Specifically, a Legendre-Chebyshev quadrature rule with weights $w_m$'s and nodes
$\Omega_m = (\zeta_m,\eta_m)$'s, $m=1,...,N_a \equiv N_{l} N_{c}$ is used to approximate the integral
in (\ref{rte-2d-pde}). (The meanings of $N_l$ and $N_c$ are given below.)
The nodes $\Omega_m = (\zeta_m,\eta_m)$'s are given by
\[
\zeta_m = \sqrt{1-\mu_i^2} \cos\varphi_j,\quad \eta_m =\sqrt{1-\mu_i^2} \sin\varphi_j,
\quad m = (i-1)N_{c}+j,
\]
where $\mu_i,~i=1,...,N_{l}$ denote the roots of the Legendre polynomial of degree $N_l$ and $\varphi_j=(2j-1)\pi /N_c,~j=1,...,N_{c}$ are the nodes based on a Chebyshev polynomial.
Once the discrete angles are defined, the DOM approximation in $(\zeta, \eta)$, the DG discretization
in $(x,y)$, and the PP limiter for \eqref{rte-2d-pde} are similar to those in one dimension. To save space,
we omit the detail here. The interested reader is referred to \cite{ZhangChengHuangQiu,ZhangChengQiu2019}.
We remark that the PP limiter uses a set of special quadrature points on triangle $K$ \cite{ZhangXiaShu2012}.
The limiter guarantees the nonnegativity of the approximate radiative intensity $\hat{I}^{n+1,(\ell+1)}_{m,K}$
at the quadrature points while maintaining the mass conservation and high-order accuracy
if the cell averages are nonnegative.
Ling et al. \cite{LingChengShu2018} give a counterexample showing that $P^r$- or $Q^r$-DG
schemes on rectangular meshes can result in negative cell averages
for the two-dimensional RTE even if both the inflow boundary value and
the source term are positive and the initial condition is nonnegative.
On the other hand, we have not observed in our limited numerical experience
that $P^r$-DG schemes lead to negative cell averages on triangular meshes
(cf.~\S\ref{sec:RTE-tests}) and thus we use the linear scaling PP limiter (cf. \S\ref{positivity})
in our computation.
It is interesting to point out that the rotational PP limiter
on triangular meshes \cite{ZhangChengQiu2019} can be used for situations with negative
cell averages. Since this limiter is non-conservative, we will not discuss it further in this work.

To conclude this section, we emphasize that, since the DG-interpolation with PP limiter is positivity-preserving,
our rezoning MM-DG method with DG-interpolation is positivity-preserving as long as
the physical PDE solver on a fixed mesh is positivity-preserving.

% section 6
\section{Numerical results for RTE}
\label{sec:RTE-tests}
In this section we present numerical results obtained for the one- and two-dimensional versions of RTE
using the rezoning MM-DG method with and without the positivity-preserving (PP) limiter
as described in the previous section. For comparison purpose, we consider three variants of the DG method.
\begin{itemize}
\item The fixed mesh (FM) DG method with PP limiter: The PP limiter
is applied to the DG solution of RTE;
\item The MM-DG method with PP limiter: The PP limiter
is applied to both the DG solution of RTE and the DG-interpolation;
\item The MM-DG method without PP limiter: The PP limiter
is applied to neither the DG solution of RTE nor the DG-interpolation.
\end{itemize}
The numerical results are presented to demonstrate the performance of the DG-interpolation scheme
in the adaptive MM solution of RTE. They also show that the proposed MM-DG method with PP limiter can maintain high-order accuracy of the DG method,
preserve the positivity of radiative intensity, and be able to adapt the mesh to the dynamic
structures in the solution.

\begin{figure}[h]
\centering
\subfigure[FM $N$=1600]{
\includegraphics[width=0.30\textwidth]{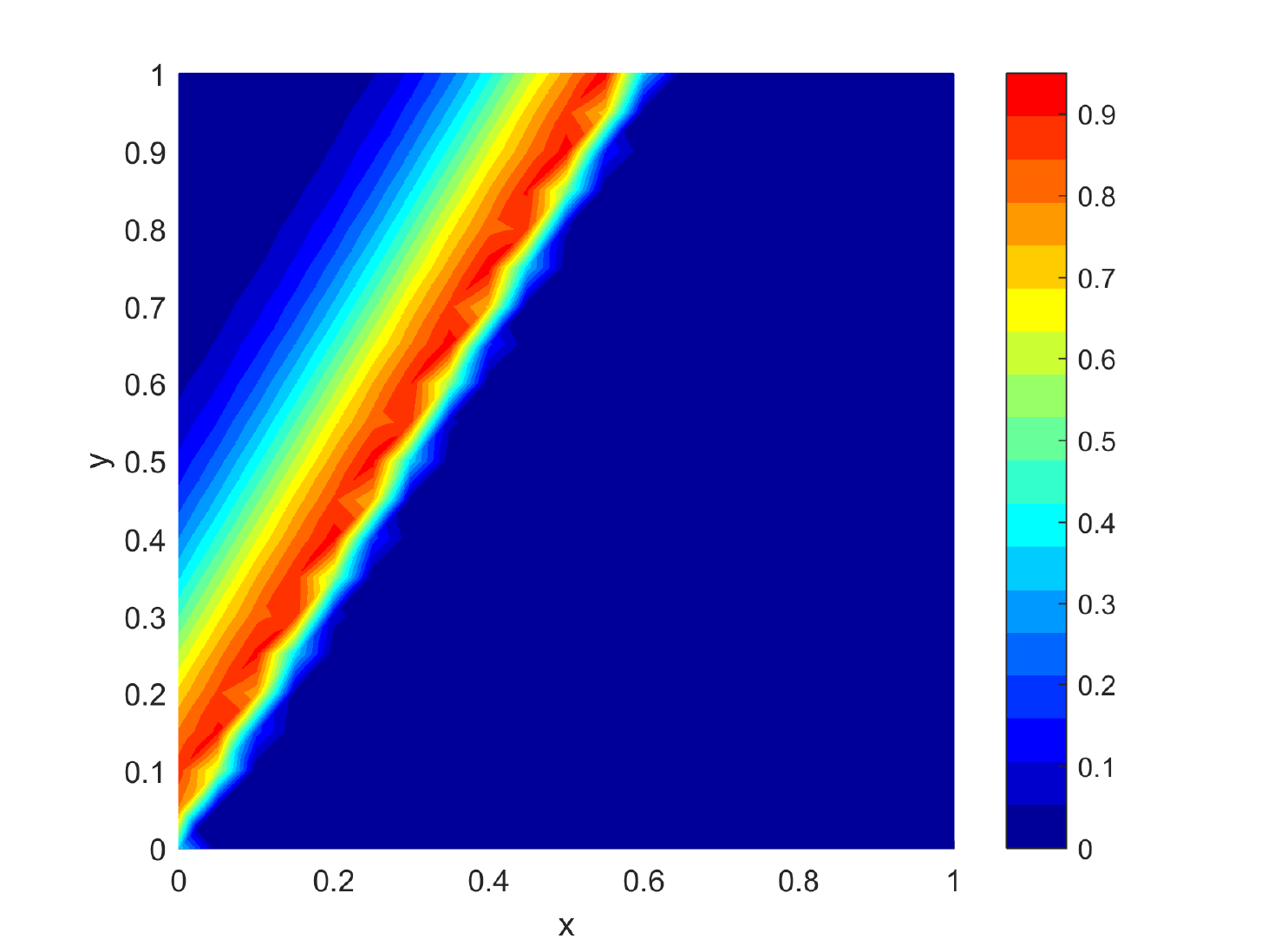}}
\subfigure[FM $N$=57600]{
\includegraphics[width=0.30\textwidth]{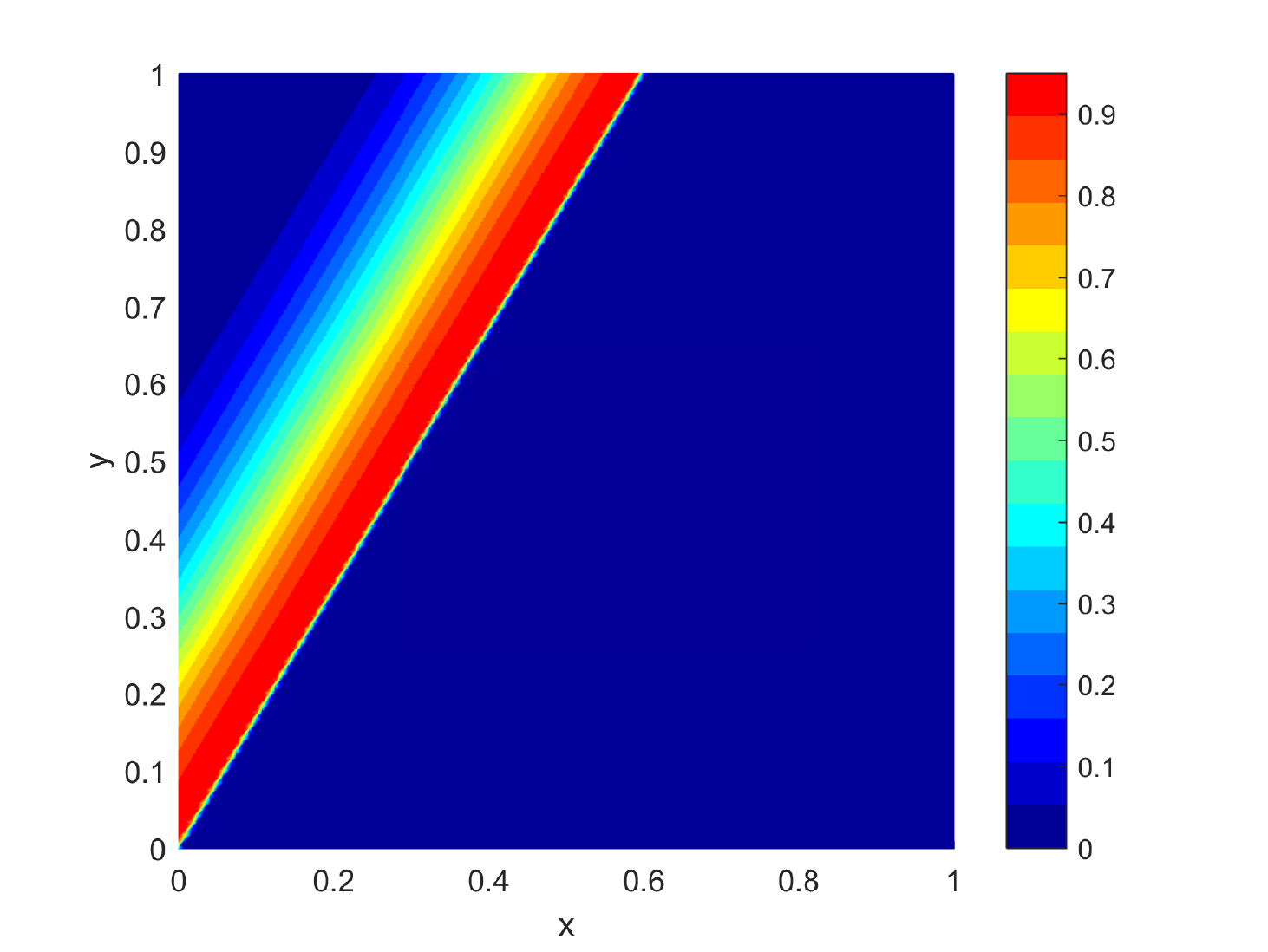}}
\subfigure[MM $N$=1600]{
\includegraphics[width=0.30\textwidth]{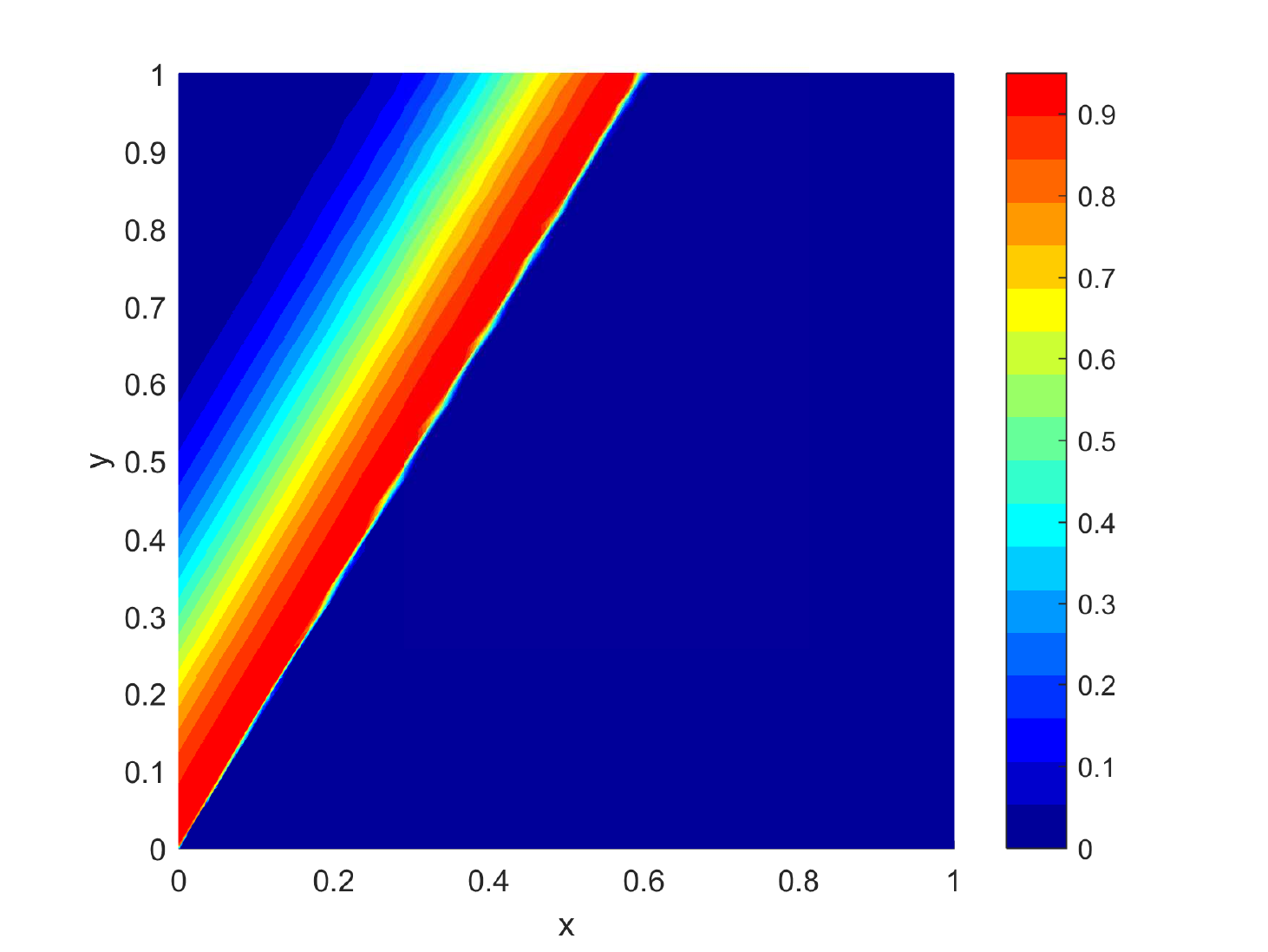}}
\\
\subfigure[MM 1600 vs FM 1600]{
\includegraphics[width=0.30\textwidth]{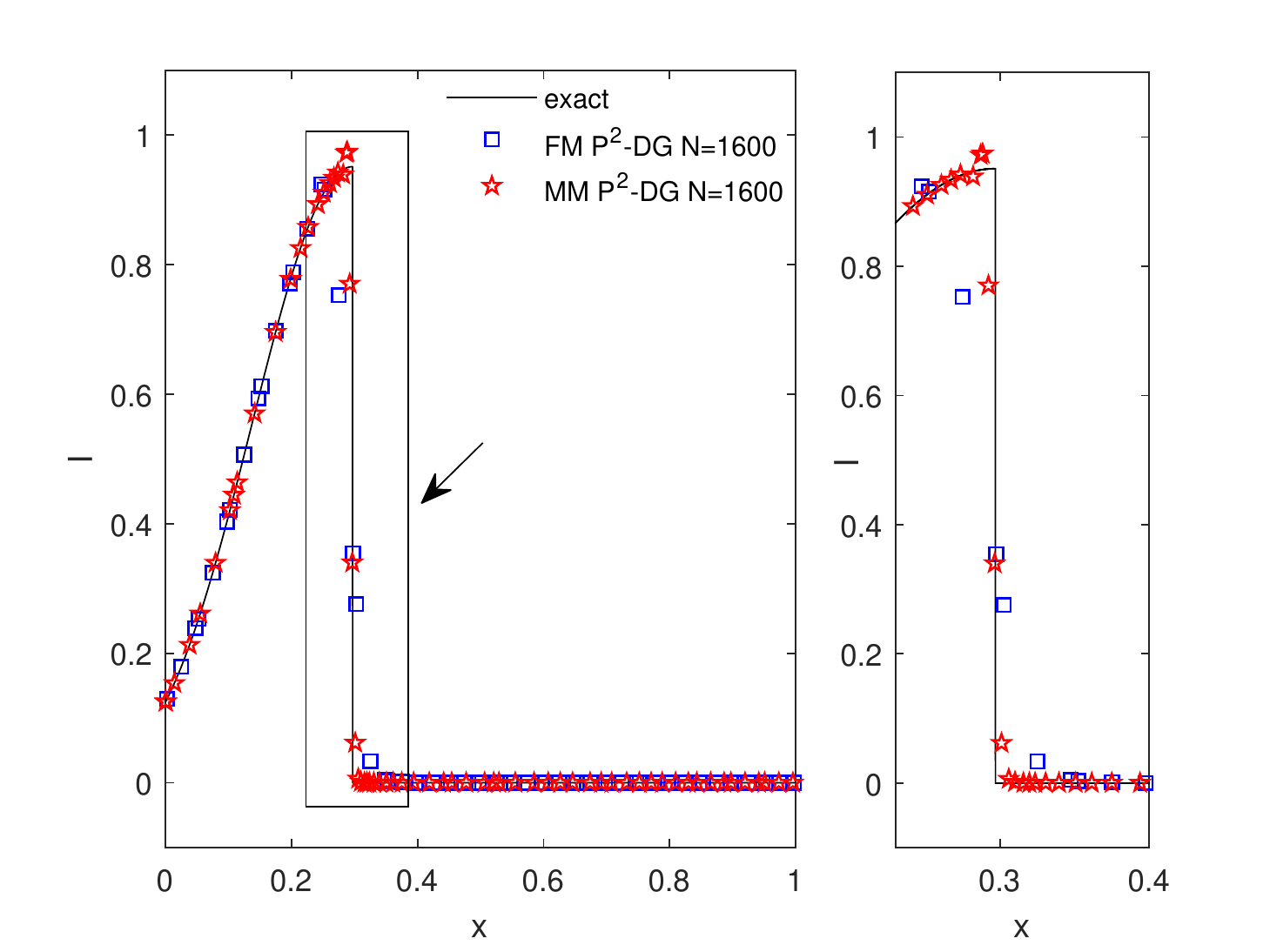}}
\subfigure[MM 1600 vs FM 57600]{
\includegraphics[width=0.30\textwidth]{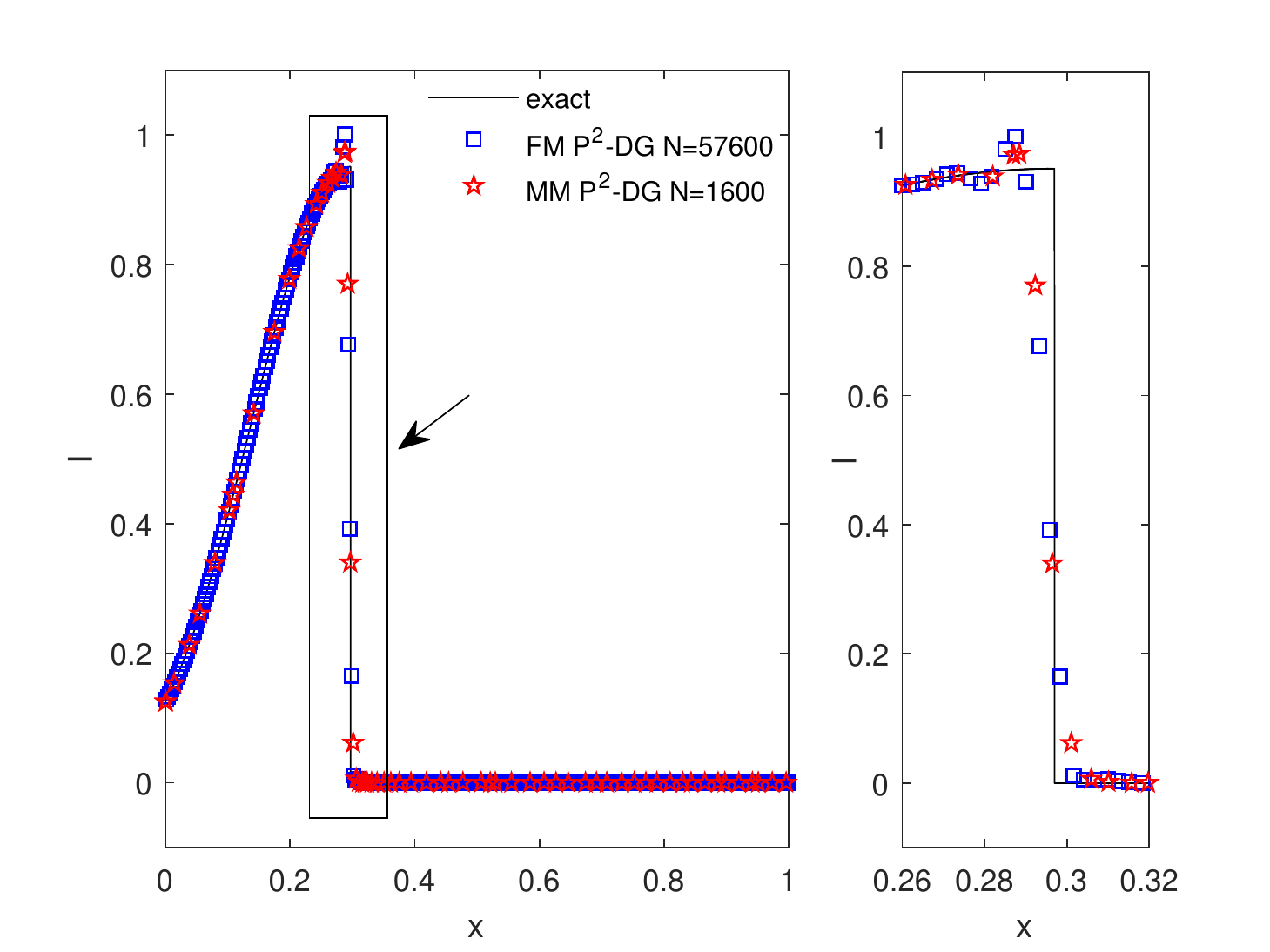}}
\subfigure[Moving Mesh at $t$=0.1]{
\includegraphics[width=0.30\textwidth]{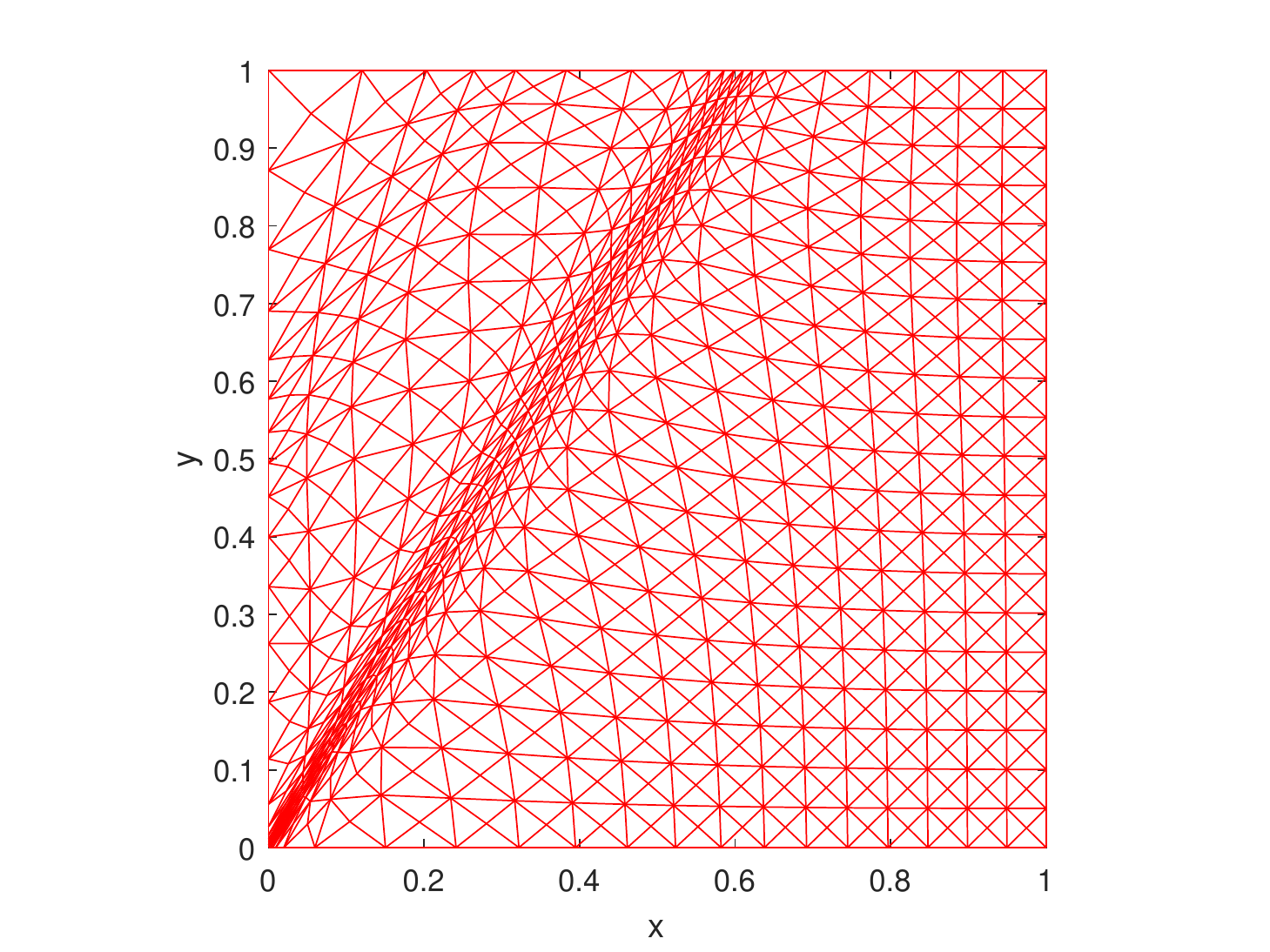}}
\caption{Example \ref{Ex2-RTE-2d}. The radiative intensity contours (and mesh) at $t=0.1$
is obtained by $P^2$-DG method with PP limiter. (d) and (e): The radiative intensity cut along the line $y=0.495$.}
\label{Fig:Ex2-2d-p2-mm}
\end{figure}

Unless otherwise stated, we use the Gauss-Legendre $P_8$
and Legendre-Chebyshev $P_8$-$T_8$ rules to discretize angular variables
for one- and two-dimensional problems, respectively,
and take the final time $T = 0.1$ and the time step size $\Delta t=2\times10^{-4}$.
For mesh movement, we take $\tau=0.01$.
The photon speed is $c=3.0\times 10^8$.
For the cases with exact solutions, the error in the computed solution is measured
in the (global) $L^1 $ and $L^\infty $ norm, i.e.,
$ \int_0^T \|e_h(\cdot,t)\|_{L^1}dt, ~ \int_0^T \|e_h(\cdot,t)\|_{L^\infty}dt.$
For two spatial dimensional examples, the initial triangular mesh is obtained by
dividing each element of a rectangular mesh into four triangles; cf. Fig.~\ref{Fig:Ex2-2d-p2-mm} (f).

% example 3 in 1D
\begin{example}\label{Ex3-RTE-1d}
(A discontinuous example of 1D RTE for the absorbing-scattering model.)
In this example we take the scattering coefficient $\sigma_s=1$, the extinction coefficient
and source term as
\begin{equation*}
\sigma_t =
\begin{cases}
1,~   &\text{for }0\leq x< 0.2\\
900,~ &\text{for }0.2\leq x< 0.6\\
90,~  &\text{for }0.6\leq x\leq1
\end{cases}
\quad \hbox{and}\quad
 q(x, \mu, t) =
\begin{cases}
 100 e^{-t},~   &\text{for }0\leq x< 0.2\\
 1,~ &\text{for }0.2\leq x< 0.6\\
 1000 e^{3t},~  &\text{for }0.6\leq x\leq1 .
\end{cases}
\end{equation*}
The initial condition is $I(x,\mu,0)=15x$
and the boundary condition is
\begin{equation*}
\begin{cases}
I(0,\mu,t)=0, \quad \quad\quad~~ \text{for }0<\mu \leq 1,~~0< t\leq T
\\
I(1,\mu,t)=15+2t,\quad \text{for }-1 \leq\mu < 0,~~0< t\leq T.
\end{cases}
\end{equation*}
\end{example}

\begin{figure}[h]
\centering
\subfigure[$P^1$-DG: $\mu=0.5255$]{
\includegraphics[width=0.30\textwidth]{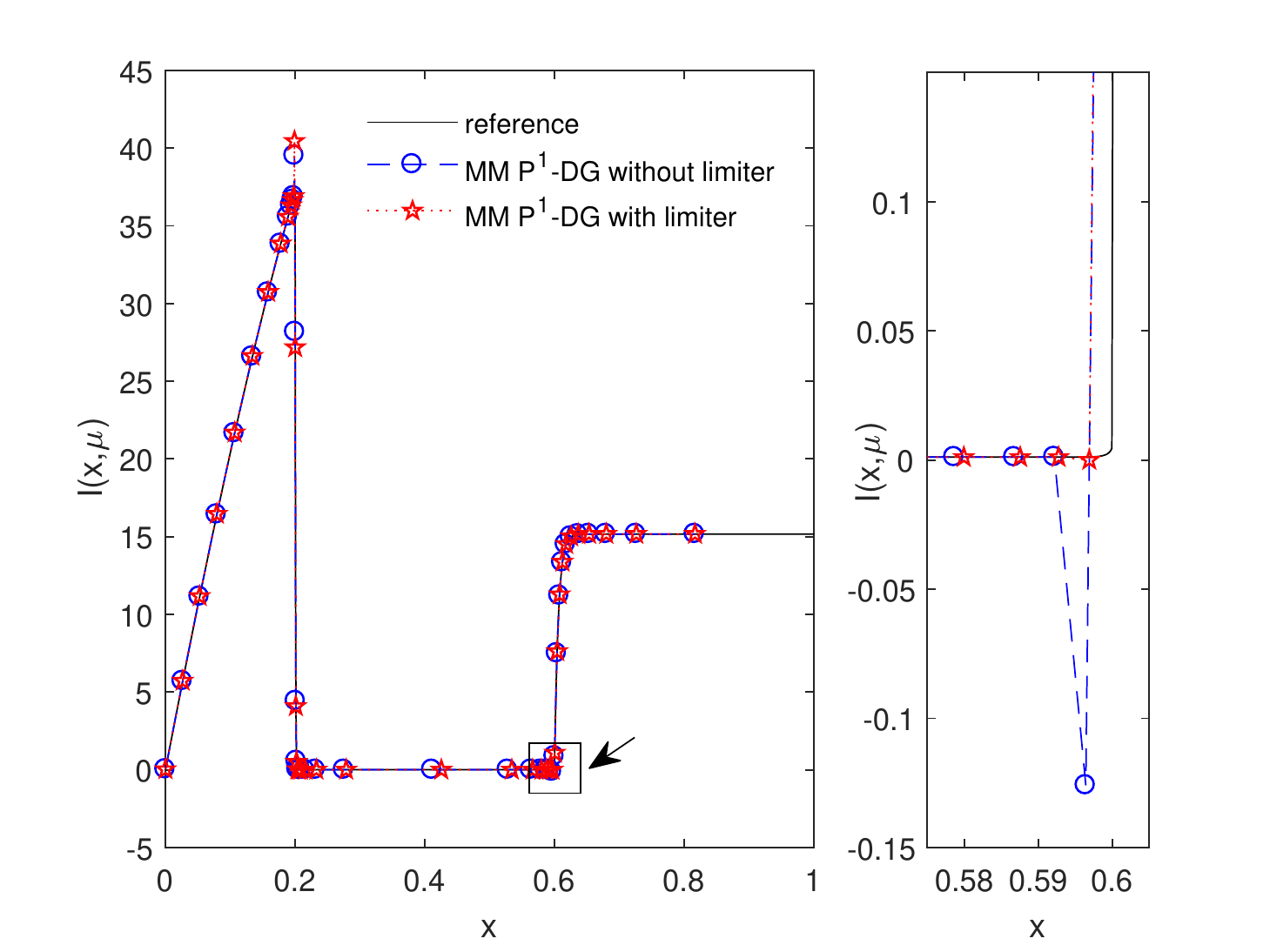}}
\subfigure[$P^1$-DG: $\mu=0.9603$]{
\includegraphics[width=0.30\textwidth]{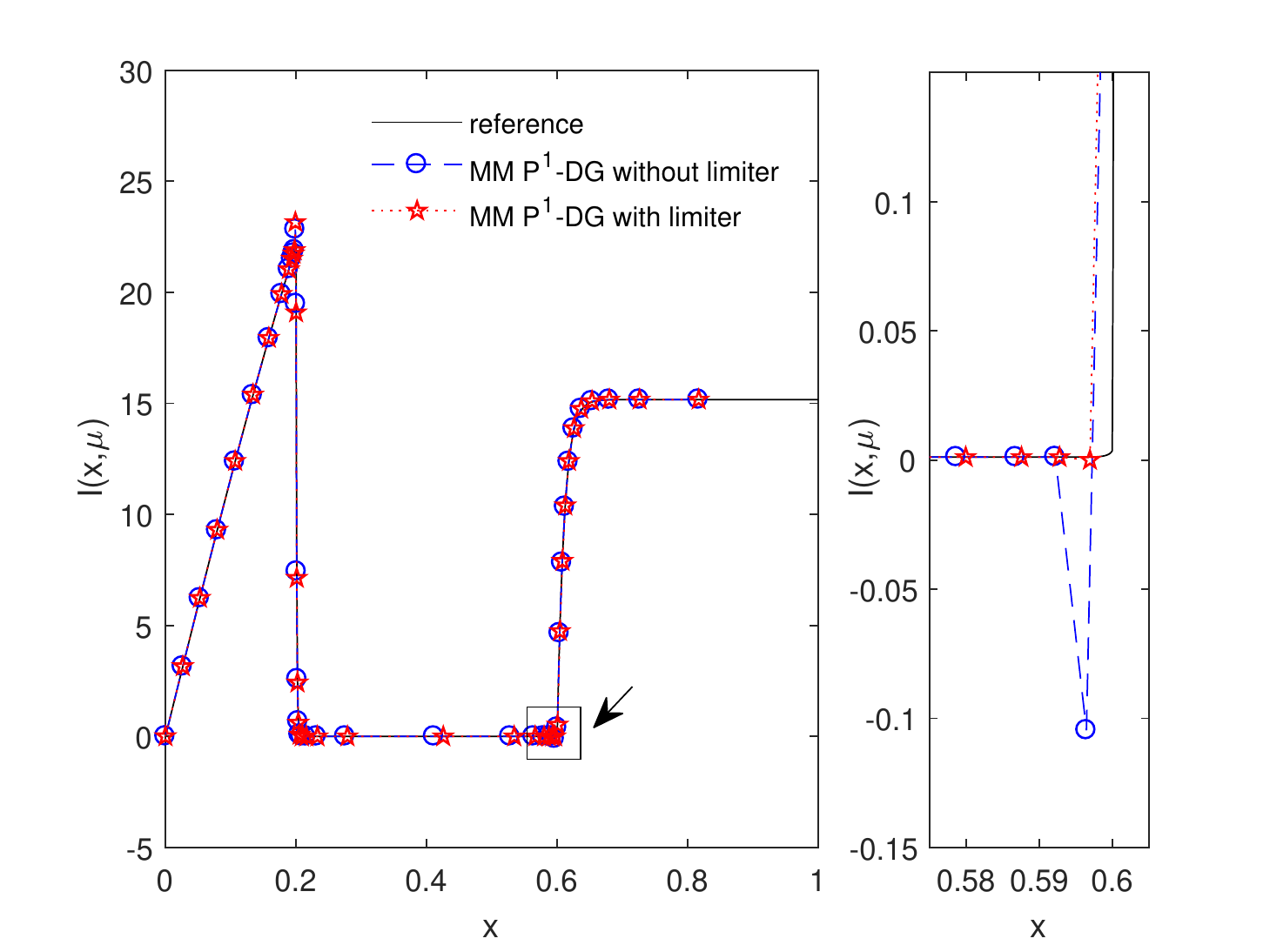}}
\\
\subfigure[$P^2$-DG: $\mu=0.5255$]{
\includegraphics[width=0.30\textwidth]{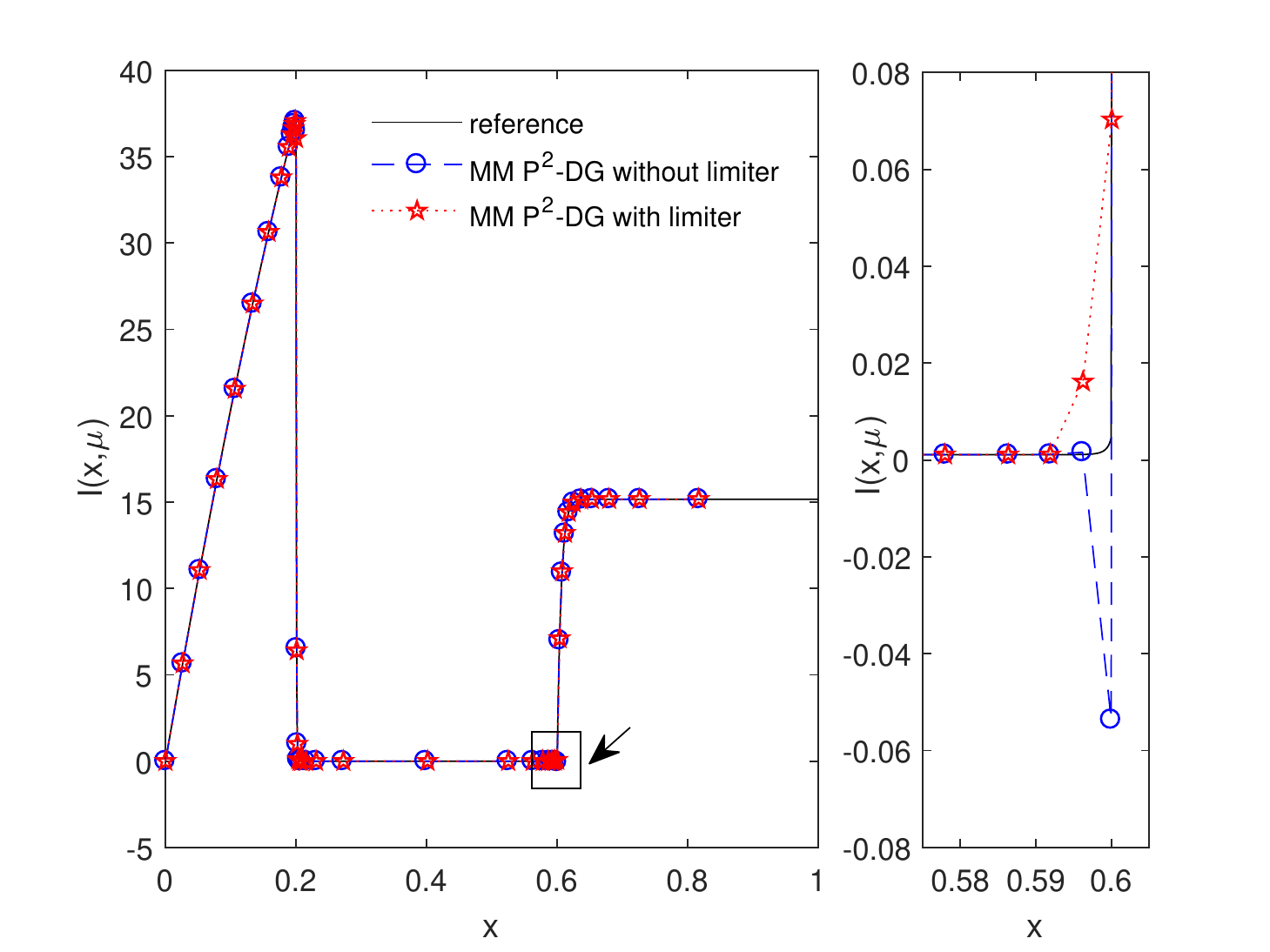}}
\subfigure[$P^2$-DG: $\mu=0.9603$]{
\includegraphics[width=0.30\textwidth]{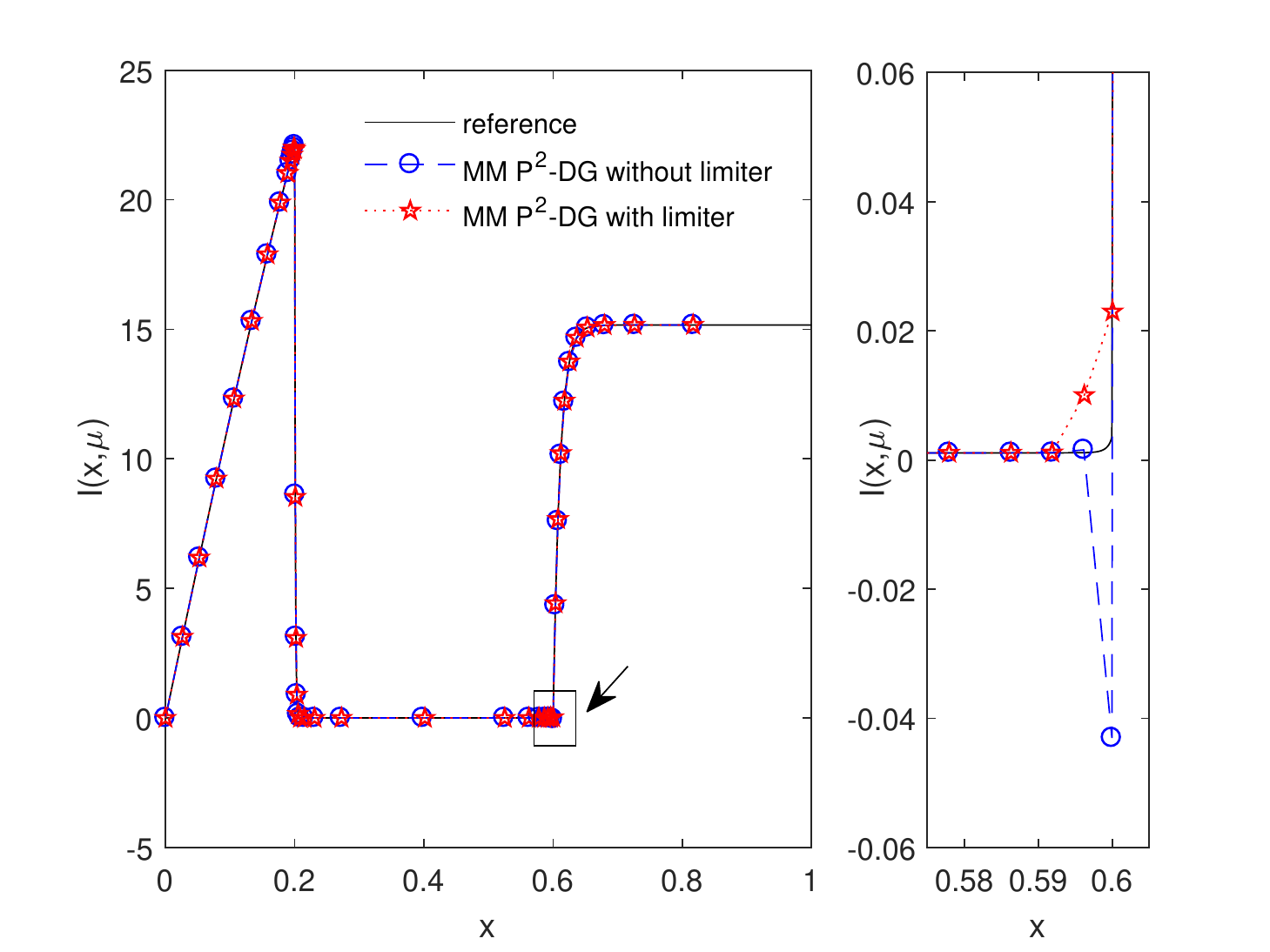}}
\caption{Example \ref{Ex3-RTE-1d}. The computed radiative intensity at the final time in the directions
$\mu=0.5255$ and $0.9603$ are obtained by the MM $P^1$-DG (Top) and $P^2$-DG (Bottom) methods with and without
PP limiter and with $N=40$. The dots represent the radiative intensity at the first Gauss-Lobatto points
on each cell.}
\label{Fig:Ex3-1d-p1and2-pp-u68}
\end{figure}

The solution of this problem has two sharp layers. Since its analytical form is unavailable,
for comparison purpose we take the numerical solution obtained by
the $P^2$-DG method with PP limiter and a fixed mesh of $N=10000$ as the reference solution.
The solutions at the final time in the directions $\mu =0.5255$ and $0.9603$
obtained by the moving mesh $P^1$-DG and $P^2$-DG methods ($N=40$) with and without PP limiter
are shown in Fig.~\ref{Fig:Ex3-1d-p1and2-pp-u68}.
We can see that the computed radiative intensity can have negative values for both $P^1$-DG and $P^2$-DG
and for fixed and moving meshes while those using the PP limiter can stay nonnegative.

\begin{figure}[h]
\centering
\subfigure[MM 80 vs FM 80]{
\includegraphics[width=0.30\textwidth]{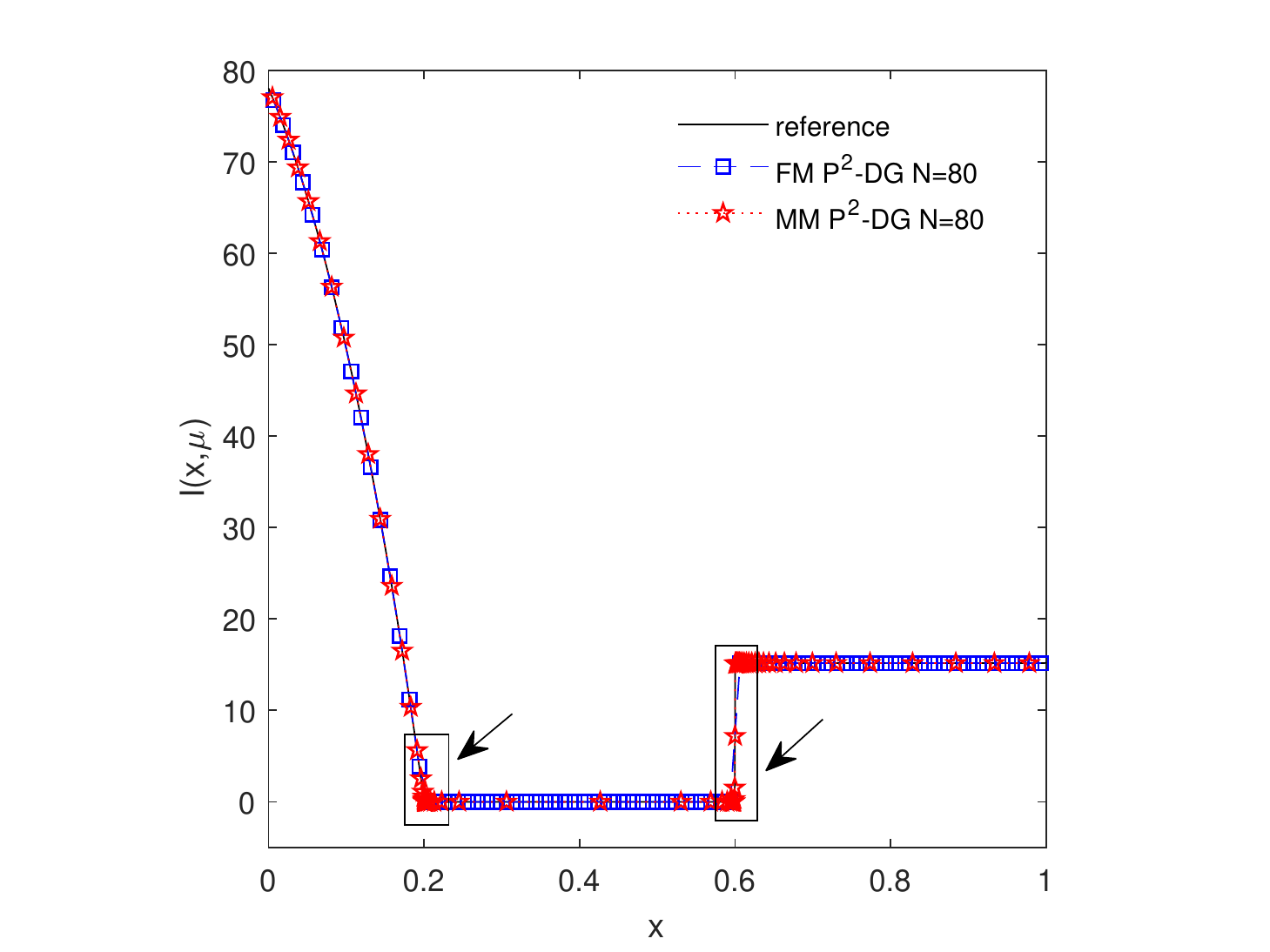}}
\subfigure[View (a) near $x$=0.2, 0.6]{
\includegraphics[width=0.30\textwidth]{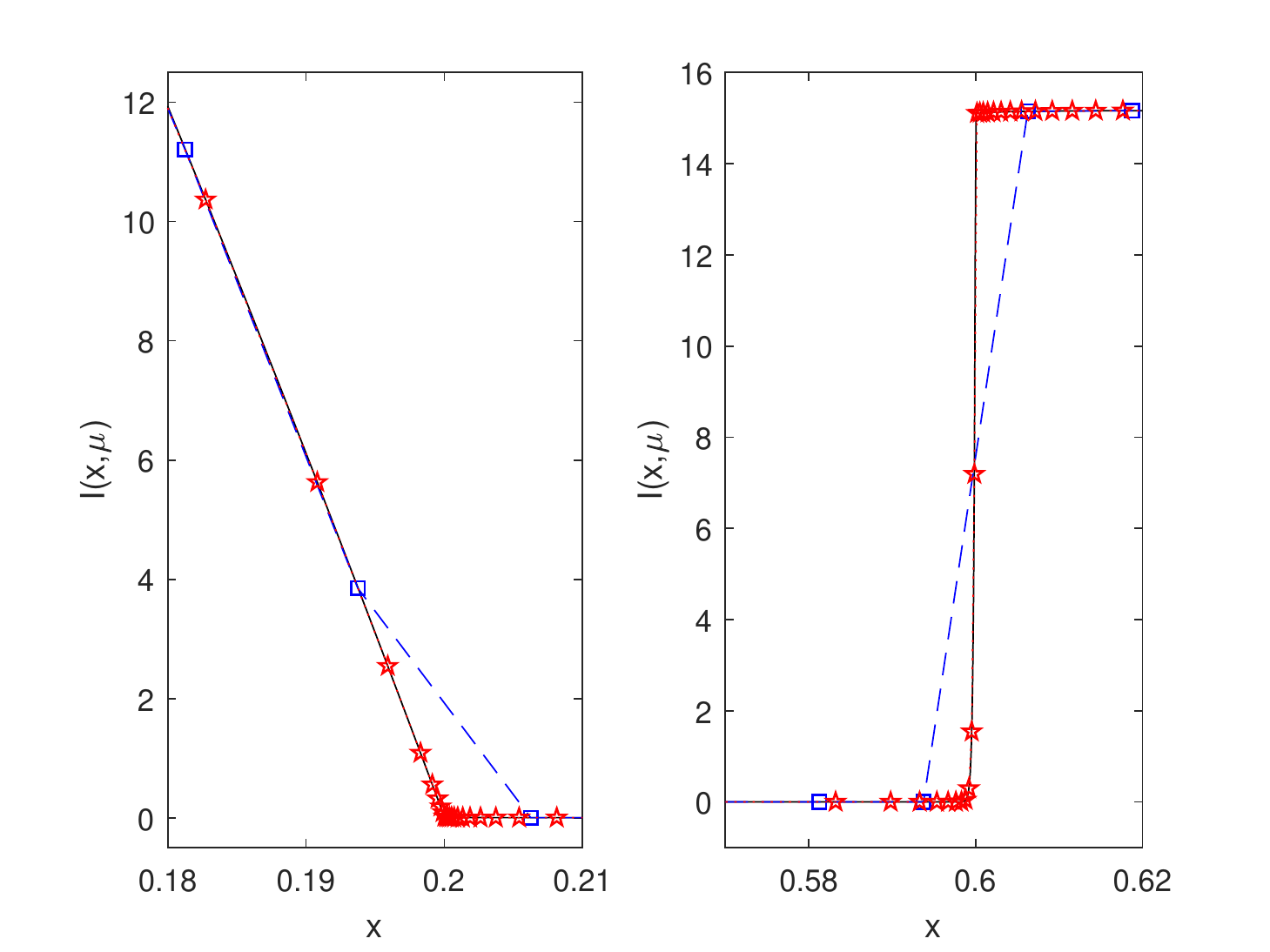}}
\\
\subfigure[MM 80 vs FM 1280]{
\includegraphics[width=0.30\textwidth]{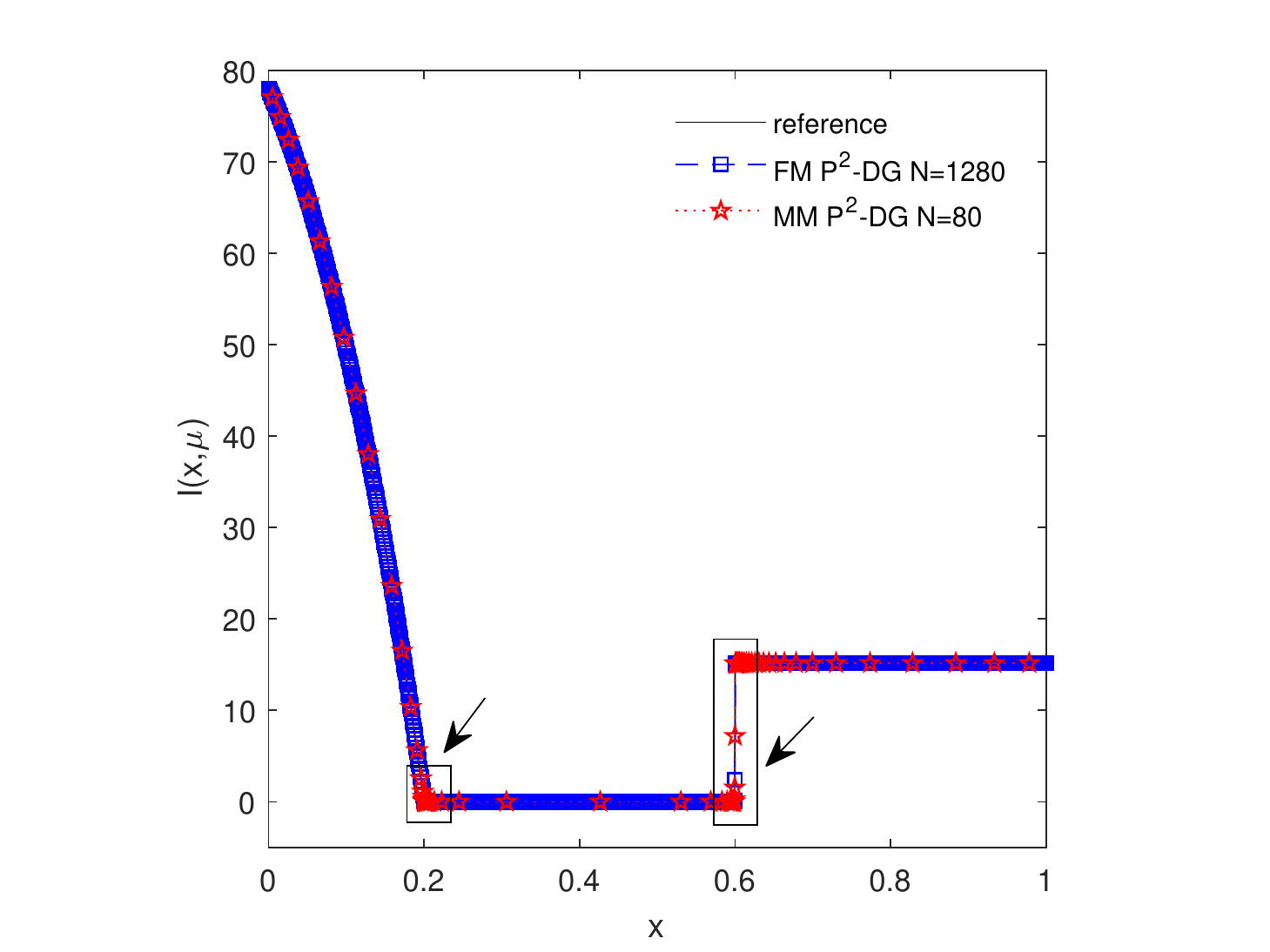}}
\subfigure[View (c) near $x$=0.2, 0.6]{
\includegraphics[width=0.30\textwidth]{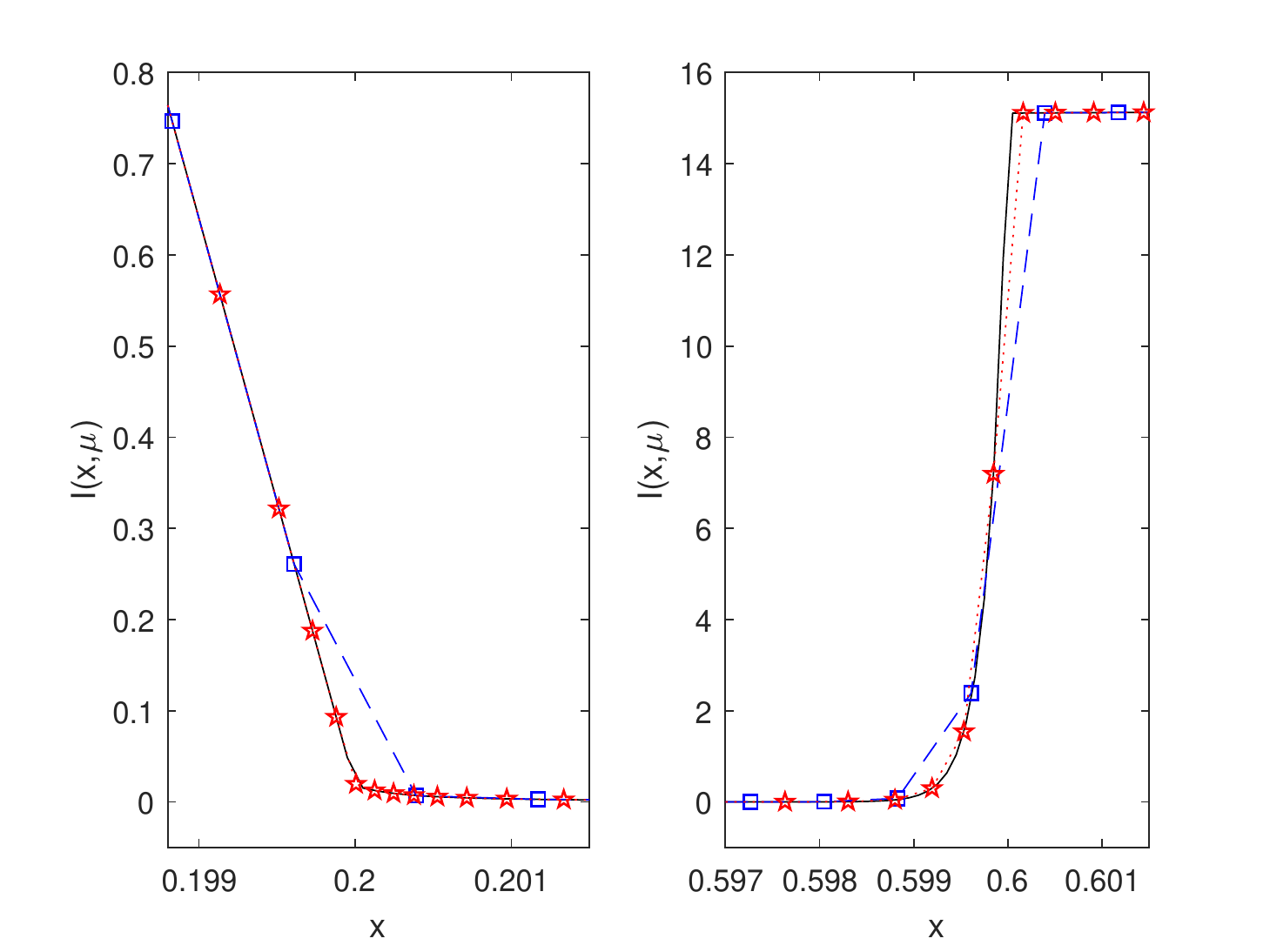}}
\caption{Example \ref{Ex3-RTE-1d}. The computed radiative intensity at the final time
in the direction $\mu= -0.1834$ obtained by the MM $P^2$-DG method ($N=80$) with PP limiter
is compared with those obtained by the FM $P^2$-DG method with PP limiter of $N=80$ and $1280$.
The dots represent the radiative intensity at the mid-points on each cell.}
\label{Fig:Ex3-1d-p2-mm-u4}
\end{figure}

\begin{figure}[h]
\centering
\subfigure[MM 80 vs FM 80]{
\includegraphics[width=0.30\textwidth]{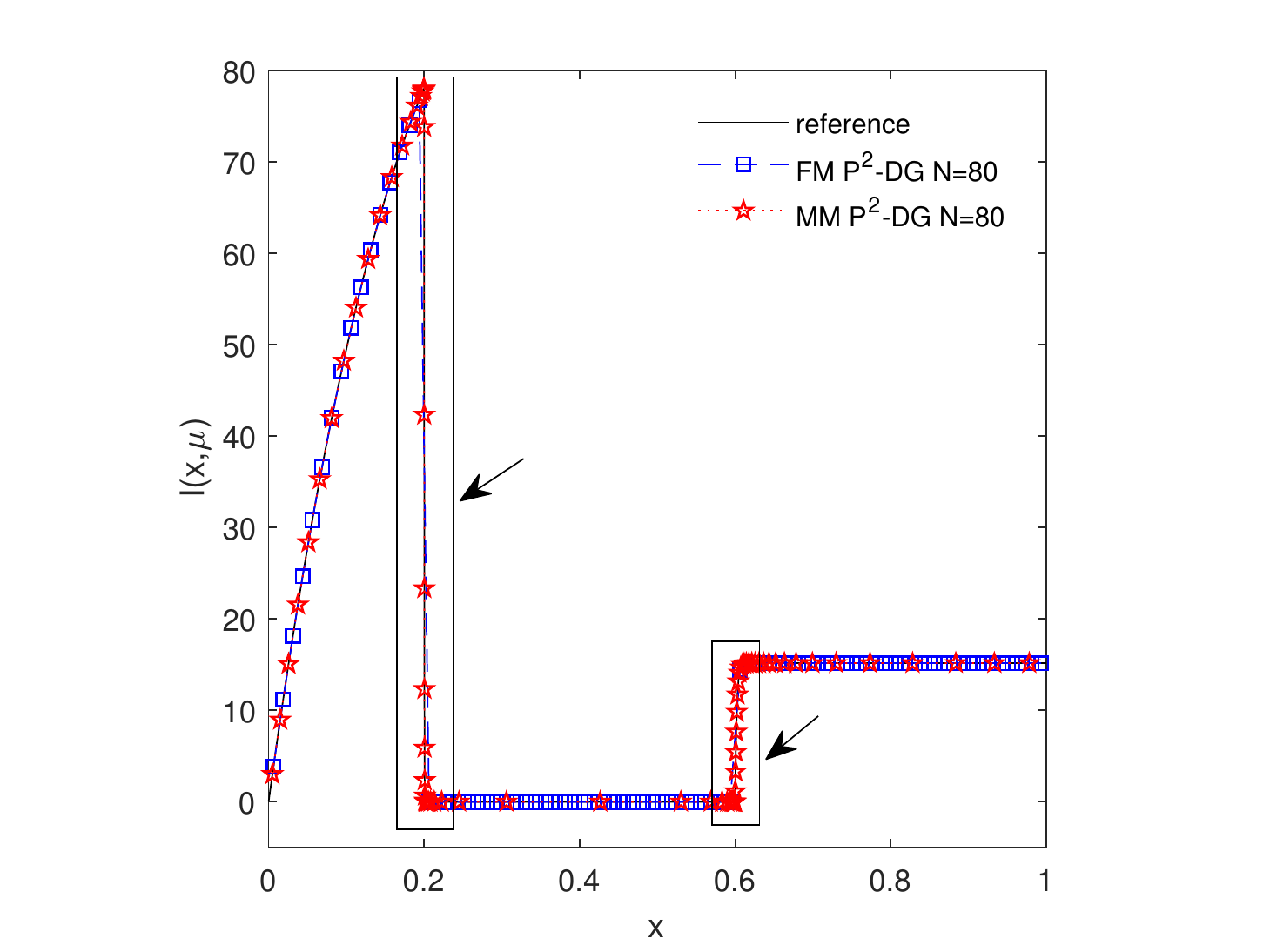}}
\subfigure[View (a) near $x$=0.2, 0.6]{
\includegraphics[width=0.30\textwidth]{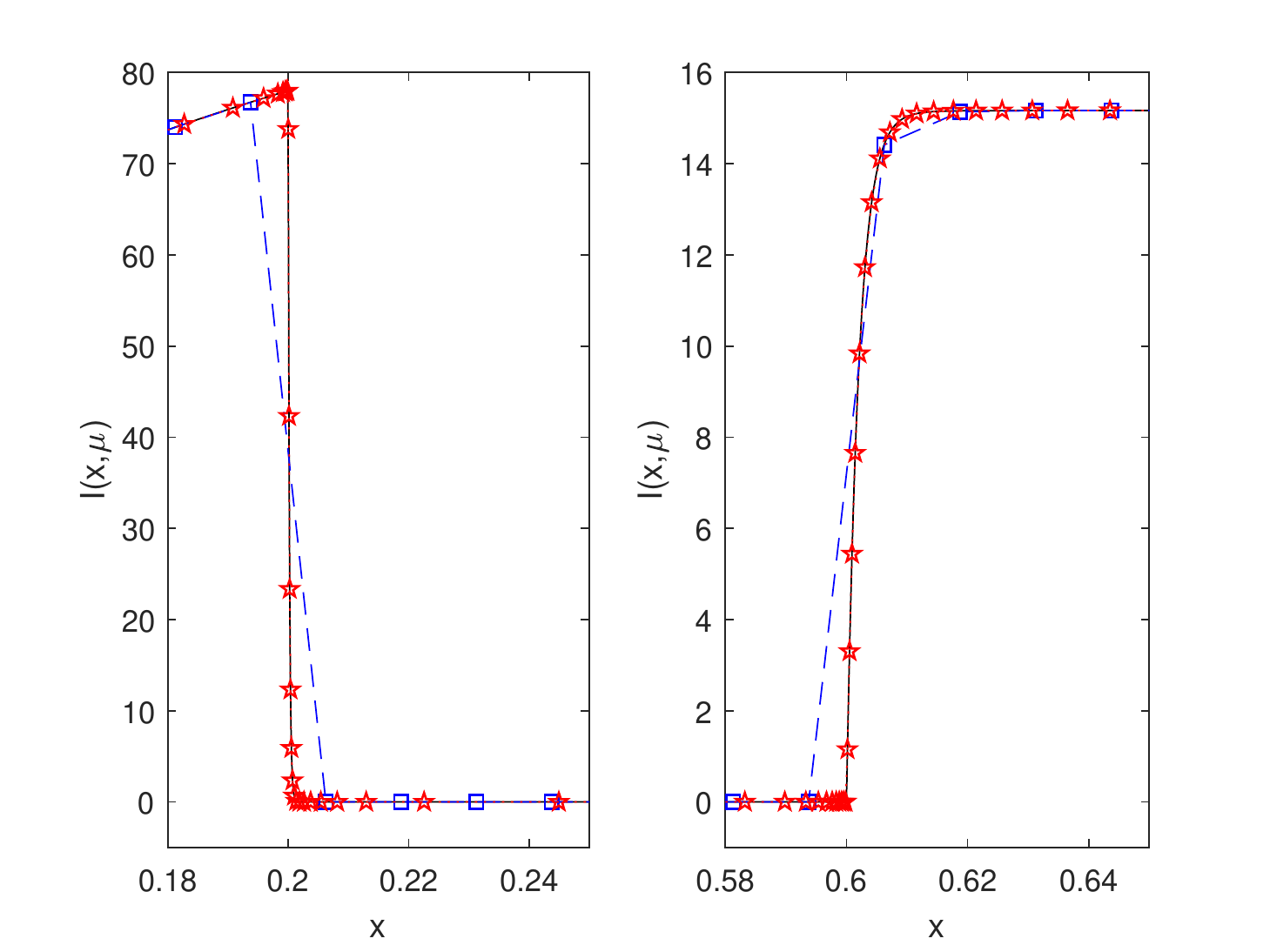}}
\\
\subfigure[MM 80 vs FM 1280]{
\includegraphics[width=0.30\textwidth]{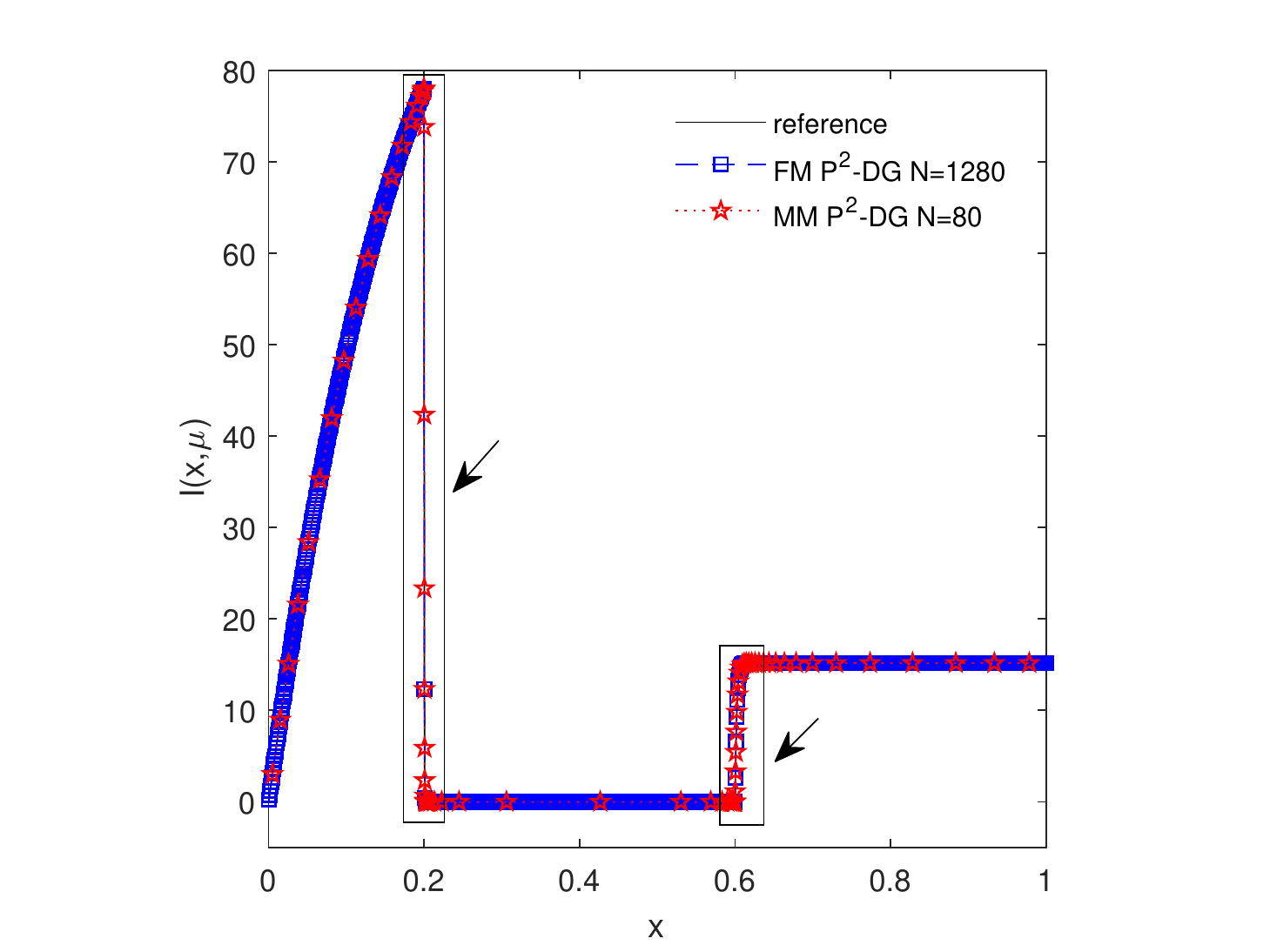}}
\subfigure[View (c) near $x$=0.2, 0.6]{
\includegraphics[width=0.30\textwidth]{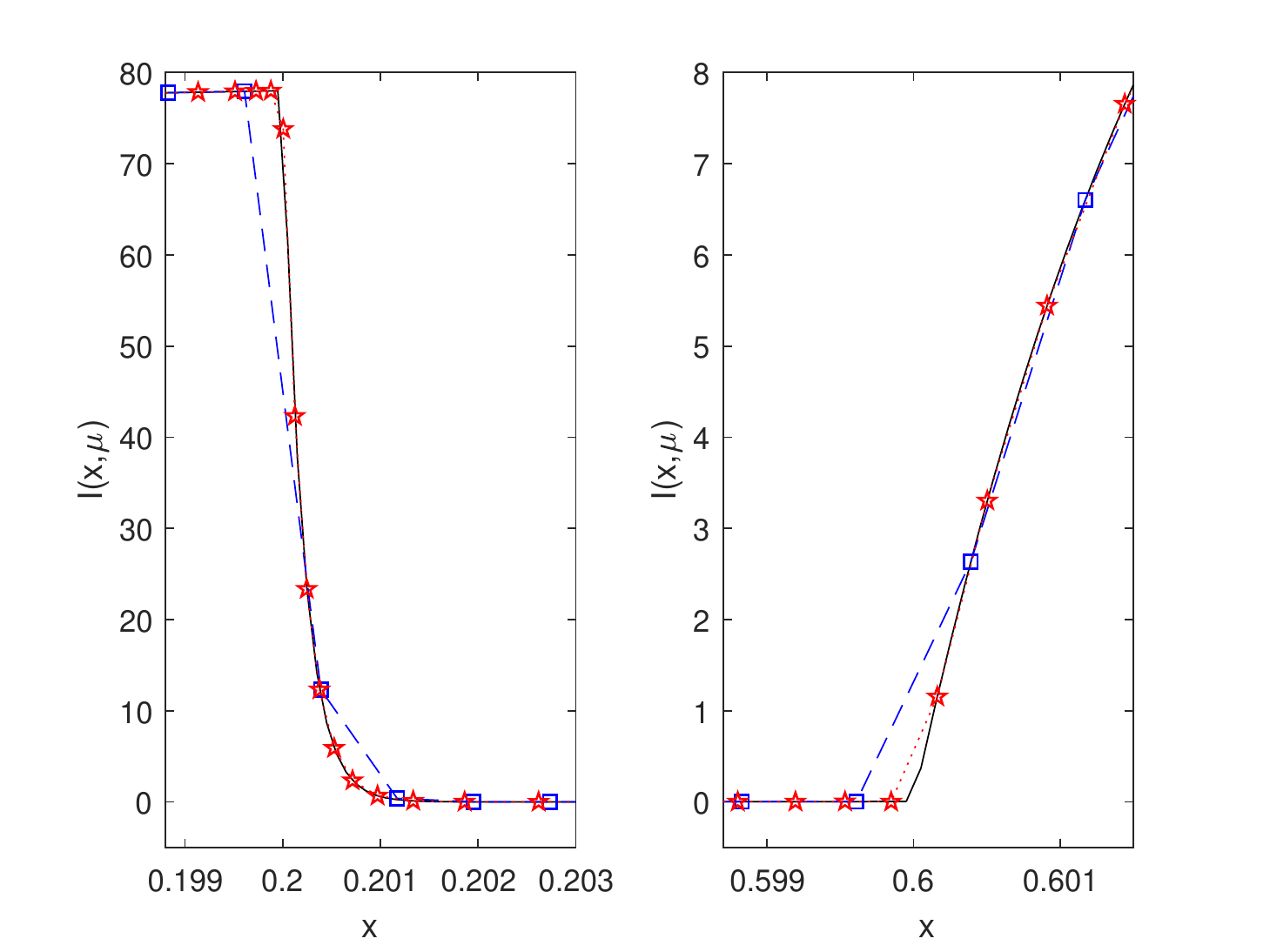}}
\caption{Example \ref{Ex3-RTE-1d}. The computed radiative intensity at the final time
in the direction $\mu=0.1834$ obtained by the MM $P^2$-DG method ($N=80$) with PP limiter
is compared with those obtained by the FM $P^2$-DG method with PP limiter of $N=80$ and $1280$.
The dots represent the radiative intensity at the mid-points on each cell.}
\label{Fig:Ex3-1d-p2-mm-u5}
\end{figure}

\begin{figure}[h]
\centering
\subfigure[Mesh trajectories, $P^2$-DG]{
\includegraphics[width=0.30\textwidth]{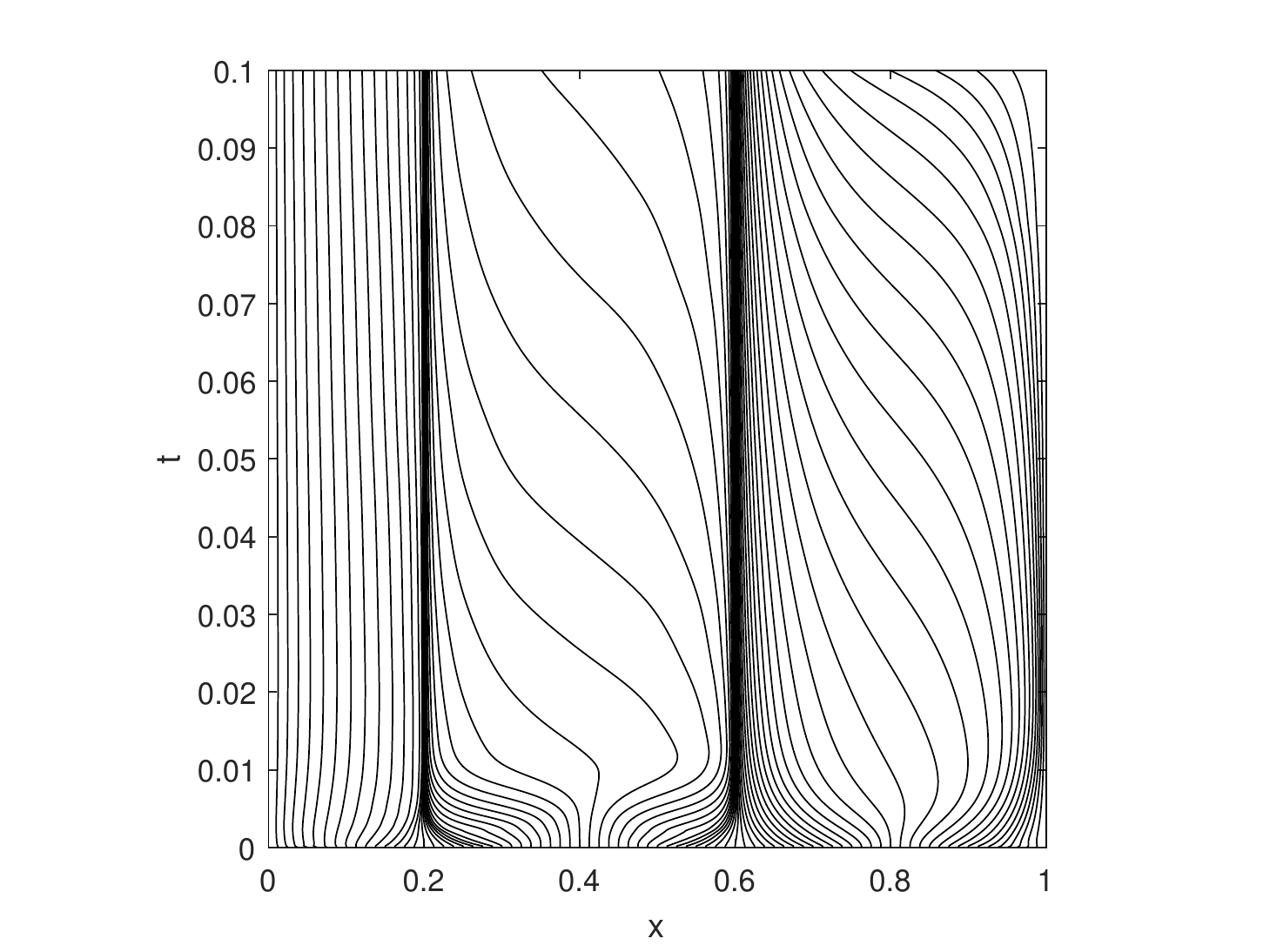}}
\subfigure[$N_{\varsigma}$, $P^1$-DG]{
\includegraphics[width=0.30\textwidth]{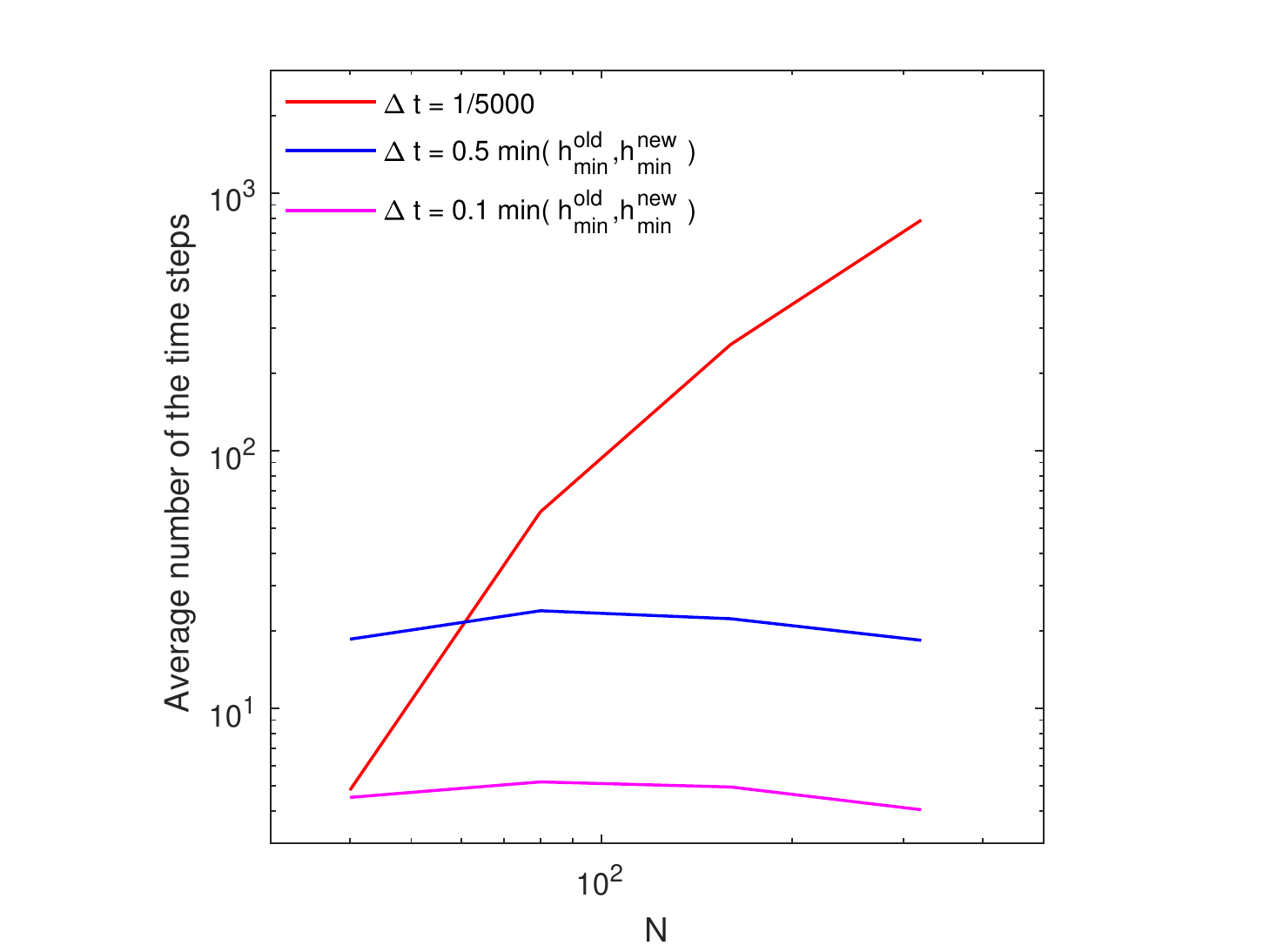}}
\subfigure[$N_{\varsigma}$, $P^2$-DG]{
\includegraphics[width=0.30\textwidth]{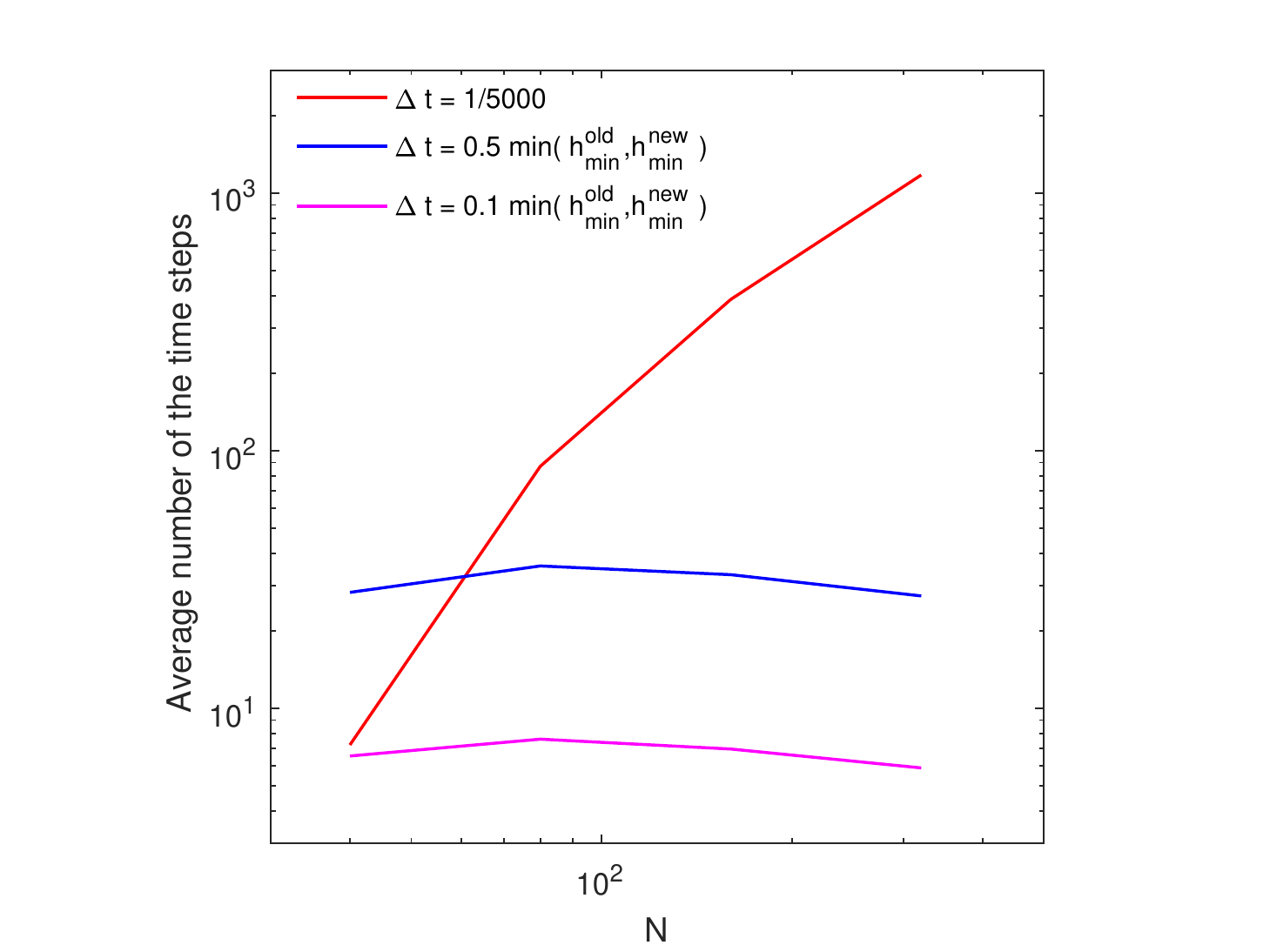}}
\caption{Example \ref{Ex3-RTE-1d}.
(a): The mesh trajectories are obtained by the MM $P^2$-DG method
with PP limiter and $N=80$. (b) and (c):
The average number of time steps used in the DG-interpolation for the MM-DG method with PP limiter.}
\label{Fig:Ex3-ld-Nsteps}
\end{figure}

The mesh trajectories for the MM $P^2$-DG method with PP limiter are shown
in Fig.~\ref{Fig:Ex3-ld-Nsteps} (a) which demonstrates the ability of the method to concentrate
mesh points in the regions of sharp layers.
The solution in the direction $\mu =-0.1834$ and 0.1834
obtained by the MM $P^2$-DG method ($N=80$) with PP limiter
is compared with the $P^2$-DG method with PP limiter and the fixed mesh of $N=80$
and $1280$ in Fig.~\ref{Fig:Ex3-1d-p2-mm-u4} and \ref{Fig:Ex3-1d-p2-mm-u5}, respectively.
The results show that the MM solution ($N=80$) is more accurate
than those with fixed meshes of $N=80$ and $1280$.
The figures also show that our MM method with PP limiter has the ability
to preserve the radiative intensity positivity.

To show the cost of the DG-interpolation in the MM-DG method with PP limiter,
we plot the average number of time steps $N_\varsigma$ in Fig.~\ref{Fig:Ex3-ld-Nsteps} (b,c).
One can see that $N_\varsigma$ is increasing as the mesh is being refined
when a fixed time step size $\Delta t=1/5000 = 2\times 10^{-4}$ is used.
On the other hand, $N_\varsigma$ stays almost constant when the time step size is chosen
as $\Delta t = 0.5 \min(h_{min}^{old},h_{min}^{new})$ and $0.1 \min(h_{min}^{old},h_{min}^{new})$
and is larger for the former than the latter.
These are consistent with the analysis in \S\ref{sec:Nsteps}.
%PP
%
%MM vs FM

%
%
% example 1 in 2D
\begin{example}\label{Ex1-RTE-2d}
(An accuracy test of 2D RTE for the absorbing-scattering model.)
In this example, we take $\sigma_t=22000$, $\sigma_s=1$. The source term and initial and boundary
conditions are chosen such that
the exact solution is given by
\[
I(x,y,\zeta,\eta,t)=e^{t}\big{(}(\zeta^2+\eta^2)\cos^4( \frac{\pi}{2}(x+y)) + 10^{-14}\big{)}.
\]
\end{example}

For this problem, the computed radiative intensity can have negative values
for both the $P^1$-DG and $P^2$-DG methods.
The error and convergence order for $P^2$-DG methods is shown
in Table \ref{tab:Ex1-2d-p2-error}.
(The results for $P^1$-DG method are omitted here to save space. They are similar to those for $P^2$-DG.)
We can see that the third-order convergence for $P^2$-DG
is achieved for fixed and moving meshes and with or without PP limiter.
The average number of time steps $N_\varsigma$ used in the DG-interpolation for the MM-DG method
is small (almost one) for relatively coarse meshes and then increases as the mesh is being
refined. This is because the mesh deformation over a time step (with a fixed time step size)
is small compared to the minimum element diameter for small $N$ and then becomes larger
for large $N$. This observation is consistent with that for the previous one-dimensional example and
the analysis in \S\ref{sec:Nsteps}.
\begin{table}[h]
\caption{Example \ref{Ex1-RTE-2d}. Error and convergence order for
three $P^{2}$-DG methods.}
%\vspace{3pt}
\centering
\label{tab:Ex1-2d-p2-error}
\begin{tabular}{ccccccc}
 \toprule
$N$&$L^1$-error & order &$L^{\infty}$-error  & order & $limiter(\%)$& $N_\varsigma$ \\
\midrule
~   &  \multicolumn{6}{c}{\em{FM $P^{2}$-DG method with PP limiter} }\\
\midrule
1600	&	7.400E-07&	 	&	1.031E-05&	 	&	5.00 	&	-	\\
6400	&	9.242E-08&	3.001 	&	1.295E-06&	2.993 	&	2.50 	&	-	\\
25600	&	1.152E-08&	3.004 	&	1.639E-07&	2.982 	&	1.25 	&	-	\\
57600	&	3.407E-09&	3.005 	&	4.909E-08&	2.973 	&	0.83 	&	-	\\
\midrule
~   &  \multicolumn{6}{c}{\em{MM $P^2$-DG method with PP limiter} }\\
\midrule
1600	&	8.163E-07&	 	&	2.561E-05&	 	&	5.00	&	1.02 	\\
6400	&	1.008E-07&	3.017 	&	3.711E-06&	2.787 	&	2.50	&	1.07 	\\
25600	&	1.177E-08&	3.099   &	3.538E-07&	3.391   &	1.25    &	1.17    \\
57600  &	3.441E-09&	3.033   &	9.399E-08&	3.269   &	0.83    &	1.67    \\
\midrule
~   &  \multicolumn{6}{c}{\em{MM $P^2$-DG method without PP limiter} }\\
\midrule
1600	&	8.117E-07&	 	&	2.561E-05&	 	&	-	&	1.02 	\\
6400	&	1.007E-07&	3.011 	&	3.711E-06&	2.787 	&	-	&	1.07 	\\
25600	&	1.177E-08&	3.097   &	3.538E-07&	3.391   &	-   &	1.17    \\
57600   &	3.440E-09&	3.032   &	9.398E-08&	3.269   &	-   &	1.67    \\
 \bottomrule	
\end{tabular}
\end{table}

% example 2 in 2D
\begin{example}\label{Ex2-RTE-2d}
(A discontinuous example of 2D RTE for the transparent model.)
In this test, we take $\sigma_t=0$,
$\sigma_s=0$, $q = 0$, $\zeta = 0.3$, and $\eta = 0.5$.
The computational domain is $(0,1)\times (0,1)$.
The initial and boundary conditions are
\begin{align*}
& I(x,y,\zeta,\eta,0)=\begin{cases}
\varepsilon,\quad &\text{for }y < \frac{\eta}{\zeta}x
\\
\cos^6\big{(}\frac{\pi}{2}y\big{)}, \quad & \text{otherwise}
\end{cases}
\\
& I(0,y,\zeta,\eta,t)=\cos^6\Big{(}\frac{\pi}{2}y\Big{)}\cos^{10}(t),\qquad
I(x,0,\zeta,\eta,t)=\varepsilon,
\end{align*}
where $\varepsilon = 10^{-14}$.
The exact solution of this example is
\begin{equation*}
I(x,y,\zeta,\eta,t)=\begin{cases}
\varepsilon,
\quad &\text{for }y < \frac{\eta}{\zeta}x
\\
\cos^6\big{(}\frac{\pi}{2}(y-\frac{\eta}{\zeta}x )\big{)}\cos^{10}(t-\frac{x}{c\zeta} ),
\quad &\text{otherwise}
\end{cases}
\end{equation*}
which is discontinuous along $y=\frac{\eta}{\zeta}x$.
\end{example}

\begin{figure}[h]
\centering
\subfigure[]{
\includegraphics[width=0.30\textwidth]
{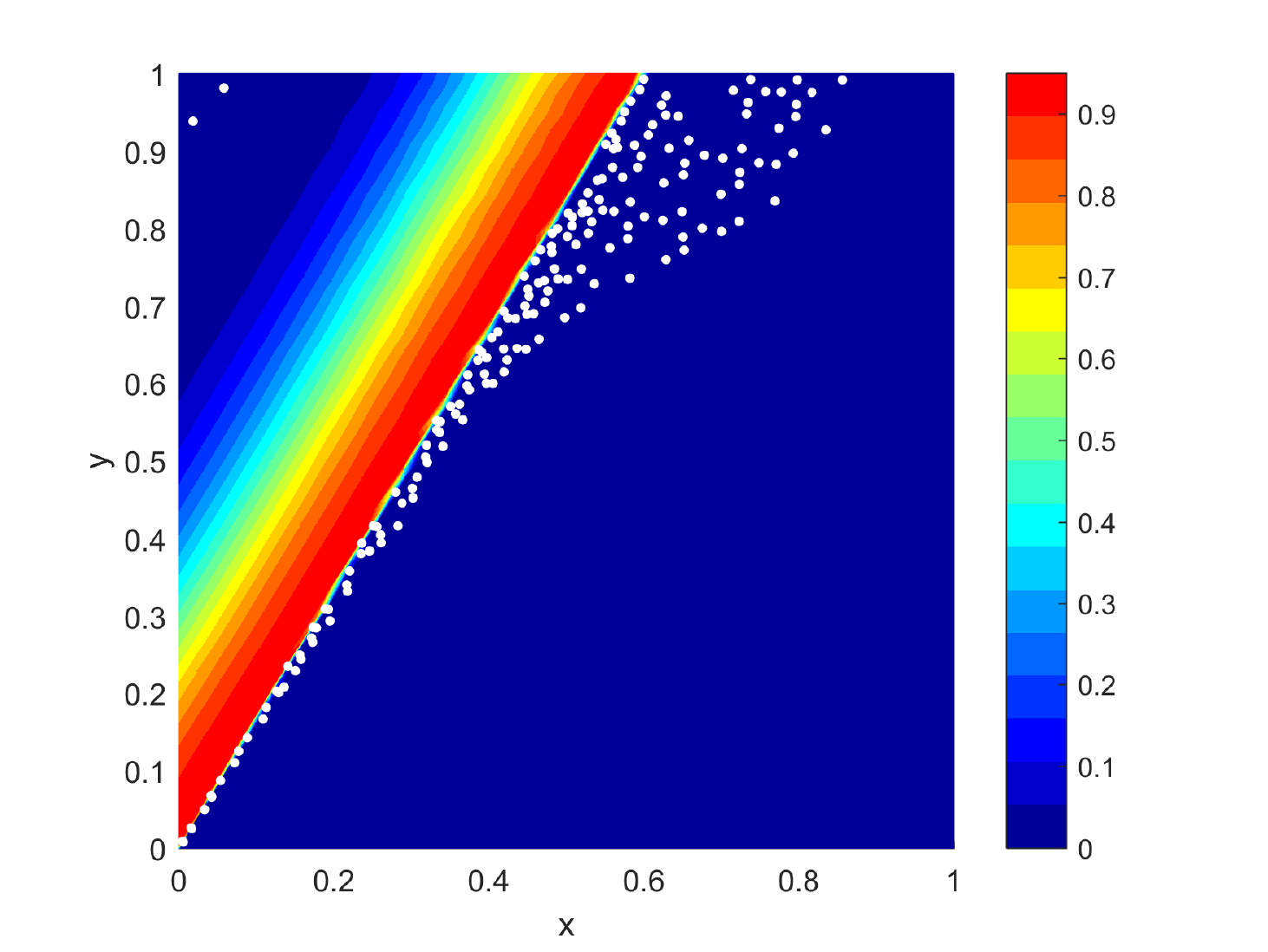}}
\subfigure[]{
\includegraphics[width=0.30\textwidth]
{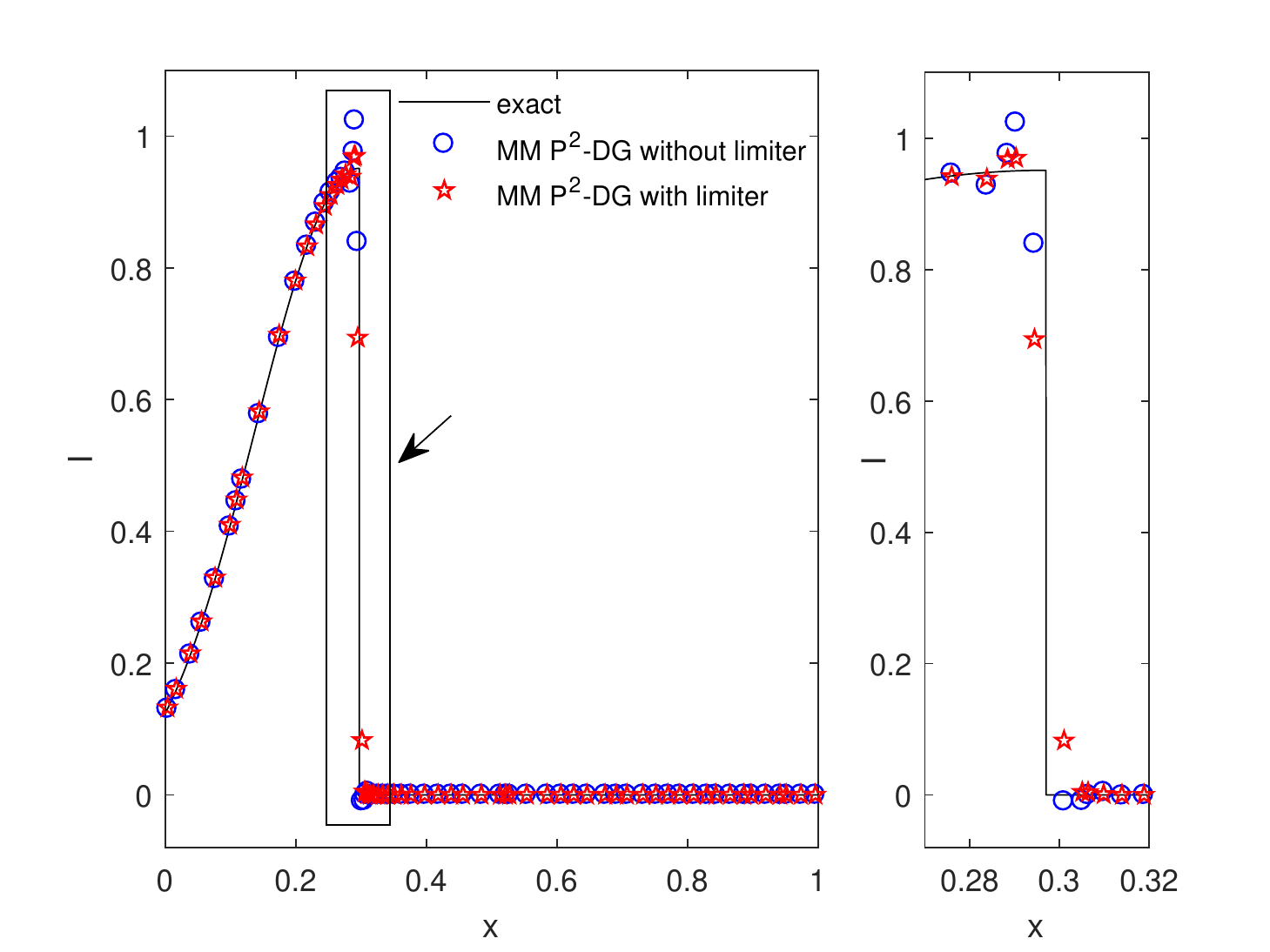}}
\caption{Example \ref{Ex2-RTE-2d}.
(a): The radiative intensity contours are obtained by the MM $P^2$-DG method with PP limiter ($N=1600$). The white dots represent the cells where the PP limiter has been applied.
(b): The radiative intensity cut along the line $y=0.495$ is obtained by the MM $P^2$-DG method
with and without PP limiter ($N=1600$).
}
\label{Fig:Ex2-2d-p2-pp}
\end{figure}

The radiative intensity obtained by the MM $P^2$-DG method with PP limiter
and $N=1600$ are plotted in Fig.~\ref{Fig:Ex2-2d-p2-pp} (a) and
the radiative intensity cut along the line $y=0.495$ obtained with and without PP limiter
is shown in Fig.~\ref{Fig:Ex2-2d-p2-pp} (b).
The cells where the PP limiter has been applied are marked with white dots.
The computed radiative intensity can have negative values for this example for the DG schemes
without PP limiter.

The contours of the radiative intensity obtained by the $P^2$-DG method with PP limiter on a moving mesh of $N=1600$ and fixed meshes of $N=1600$ and $57600$ are shown in Fig.~\ref{Fig:Ex2-2d-p2-mm} (a,b,c).
The corresponding cut along the line $y=0.495$ is plotted in Fig.~\ref{Fig:Ex2-2d-p2-mm} (d,e).
The results show that the MM solution ($N=1600$) is more accurate than that
with the fixed mesh of $N=1600$ and is comparable with that with the fixed mesh of $N=57600$.
The figures also show that our MM $P^2$-DG method with PP limiter
produces the positive radiative intensity.

\begin{figure}[h]
\centering
\subfigure[$L^1$ norm of error]{
\includegraphics[width=0.30\textwidth]{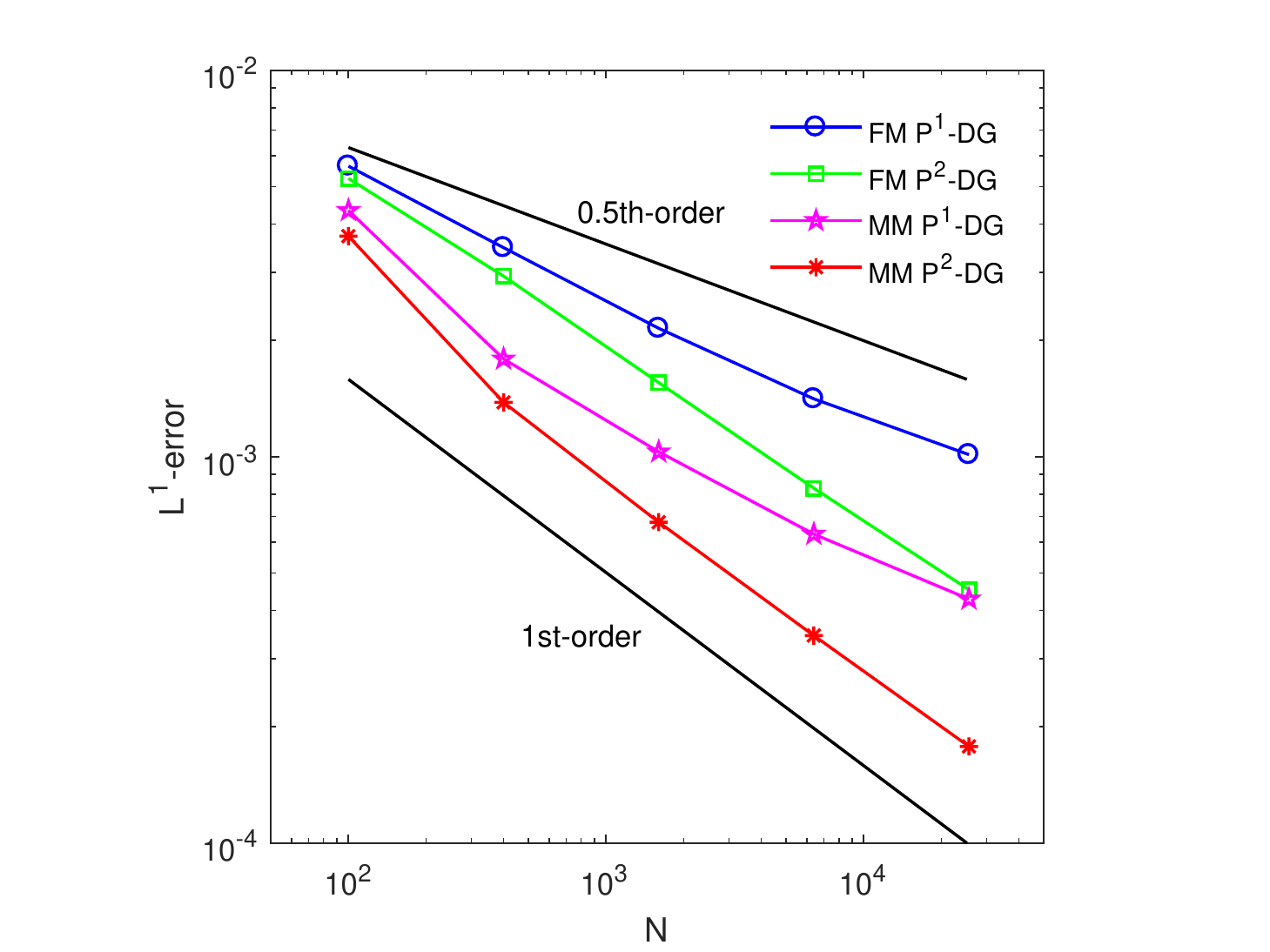}}
\subfigure[Number of time steps]{
\includegraphics[width=0.30\textwidth]{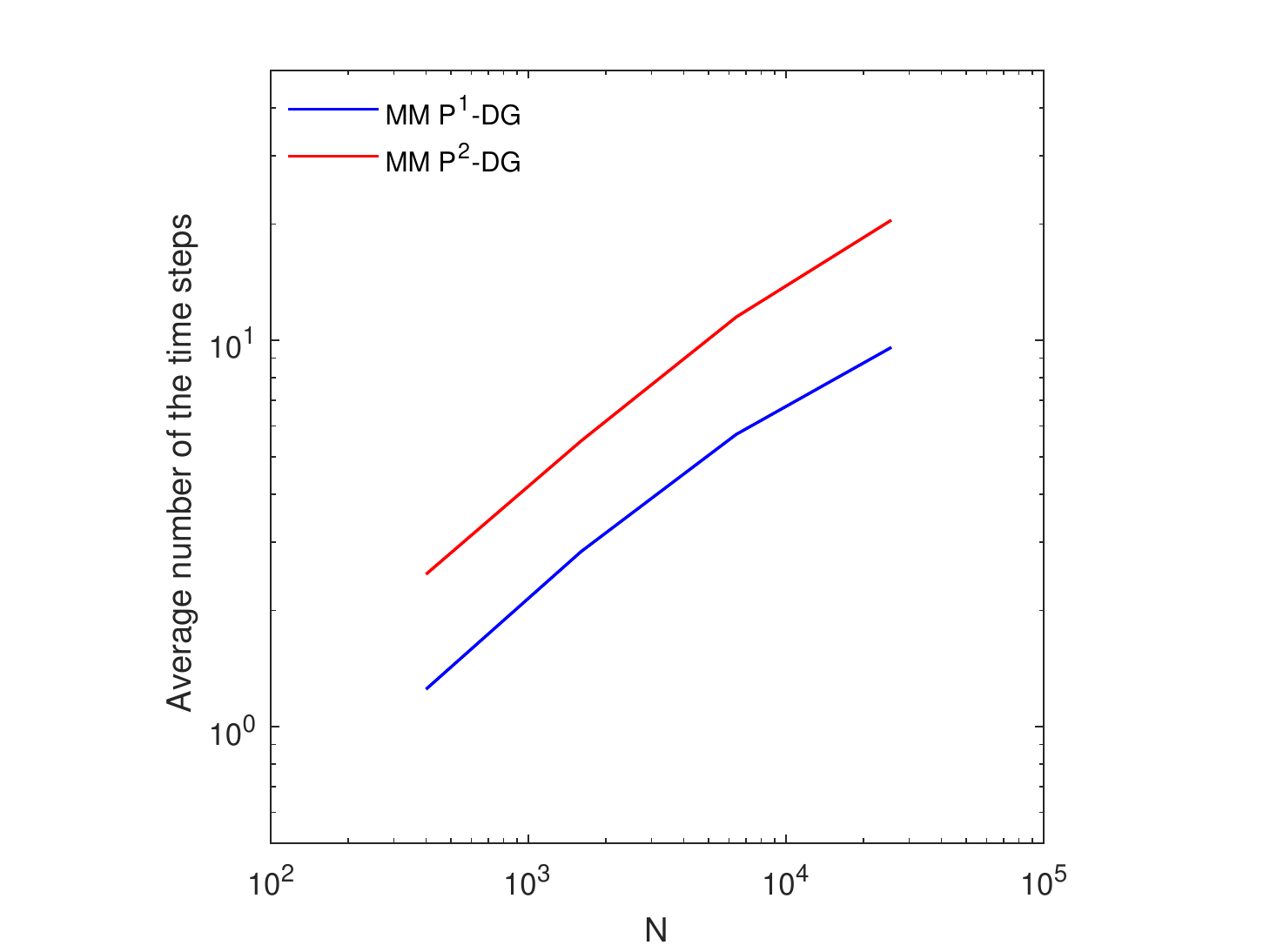}}
\subfigure[$L^1$ error vs CPU time]{
\includegraphics[width=0.30\textwidth]{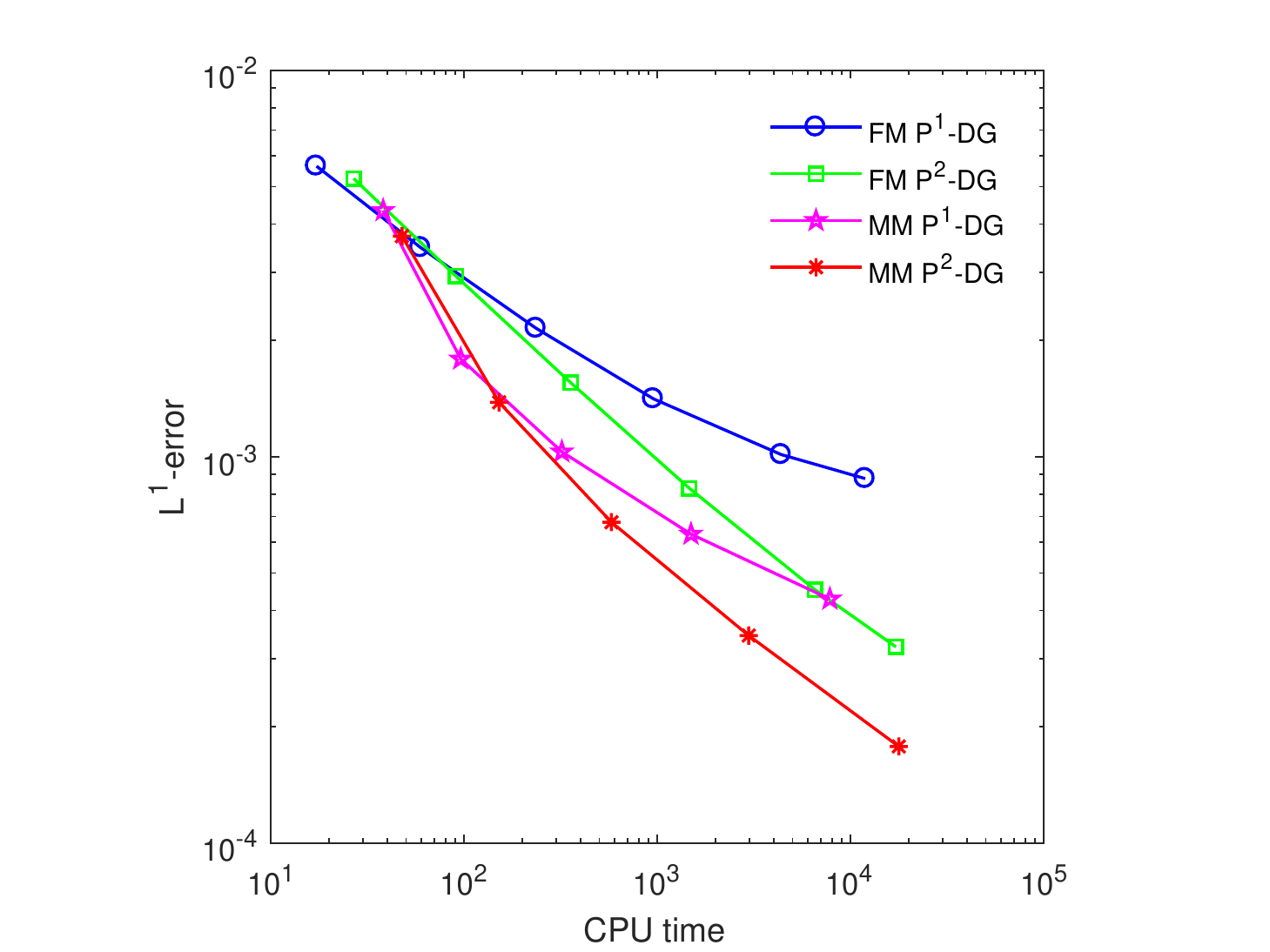}}
\caption{Example \ref{Ex2-RTE-2d}. Results are obtained by DG methods with PP limiter.
(a) Convergence history. (b): The average number of time steps used in the DG-interpolation
for the MM-DG computation. (c): The $L^1$ norm of the error is plotted against the CPU time (in seconds).
}\label{Fig:Ex2-2d-L1order-Nsteps-CPU}
\end{figure}

The error and convergence history in the $L^1$ norm are shown in Fig.~\ref{Fig:Ex2-2d-L1order-Nsteps-CPU} (a)
for the FM-DG and MM-DG methods with PP limiter.
One can see that both fixed and moving meshes lead to almost the same convergence order.
It is worth pointing out that we cannot expect the fixed/moving mesh DG can achieve the optimal
order for this problem since the solution is discontinuous.
The actual order of $P^1$-DG is about 0.5th and 1st for $P^2$-DG.
Moreover, the figures show that a moving mesh produces a more accurate solution than a fixed mesh
of the same number of elements for this example.

To show the efficiency of the MM-DG method with PP limiter,
we plot the average number of time steps used in the DG-interpolation in Fig.~\ref{Fig:Ex2-2d-L1order-Nsteps-CPU} (b),
which indicates that $N_{\varsigma}$ increases as the mesh is being refined.
We also plot the $L^1$ norm of the error against the CPU time in Fig.~\ref{Fig:Ex2-2d-L1order-Nsteps-CPU} (c).
It shows that the MM-DG method is more efficient than the fixed mesh method
and $P^2$-DG is more efficient than $P^1$-DG.

\begin{figure}[h]
\centering
\includegraphics[width=0.30\textwidth]{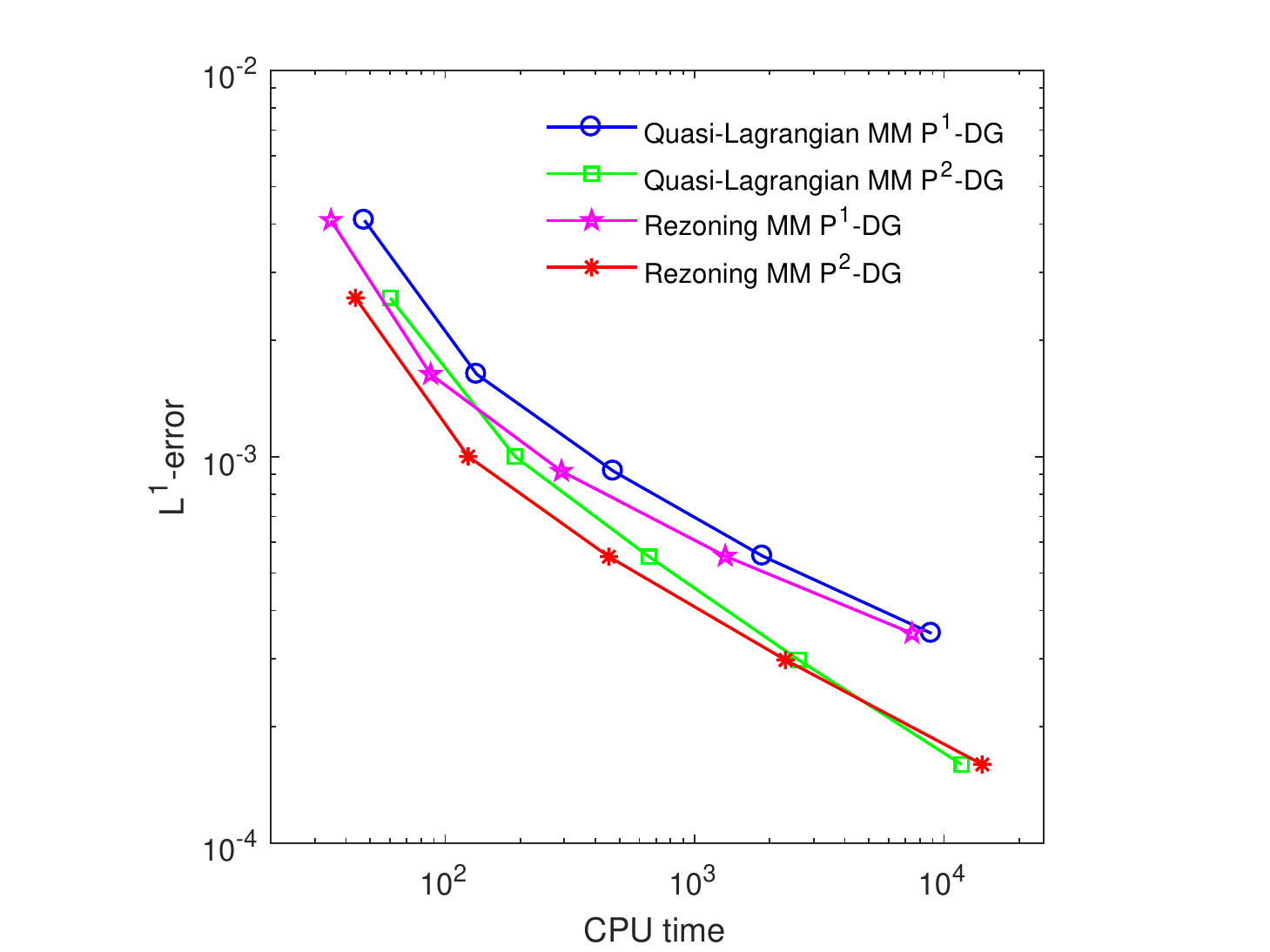}
\caption{Example \ref{Ex2-RTE-2d}.
The $L^1$ norm of error is plotted against the CPU time (in seconds)
for the rezoning MM-DG method (this work) and the quasi-Lagrangian MM-DG method
\cite{ZhangChengHuangQiu}.}
\label{Fig:Ex2-2d-Quasi-Reznoe}
\end{figure}

It is interesting to mention that a quasi-Lagrangian MM-DG method has been developed
in \cite{ZhangChengHuangQiu} for RTE. Compared to the method in the current work, it does not
require interpolation of the physical variables between old and new meshes although extra work
is needed to compute a convection term in the DG formulation of RTE that is caused by mesh movement.
It is unclear to us yet how to preserve the radiative intensity in the quasi-Lagrangian method,
which is an interesting future research topic.
To obtain a rough comparison, we plot in Fig.~\ref{Fig:Ex2-2d-Quasi-Reznoe} the $L^1$ norm of the error against CPU time for both the quasi-Lagrangian and rezoning MM-DG methods (without PP limiter).
We can see that both methods have comparable efficiency while the rezoning method
is slightly more efficient when the mesh is not very fine. As the mesh is being refined,
the DG-interpolation will need more steps and become more expensive, and then the quasi-Lagrangian method
becomes more efficient. It should also be pointed that this comparison is done with a fixed time step size.
The situation may be different when a variable time step size is used.

% section 7
\section{Conclusions}
\label{sec:conclusion}
In the previous sections we have presented a high-order DG-interpolation scheme for deforming unstructured
meshes based on the pseudo-time-dependent linear equation \eqref{pde}.
Such a scheme can be used for indirect ALE and
rezoning MM methods in numerical solution of partial differential equations.
We have shown that the scheme is conservative.
It is also positivity-preserving when a linear scaling limiter is used.
The scheme places no restrictions on the deformation of the old mesh to the new one.
The cost of the scheme has been investigated. The total cost of each use of the DG-interpolation
is $\mathcal{O}(N_v N_{\varsigma})$, where $N_v$ is the number of mesh vertices
and $N_{\varsigma}$ is the number of time steps used to integrate \eqref{pde} from $\varsigma = 0$
to $\varsigma = 1$. It is shown that $N_{\varsigma}$ depends on the magnitude of mesh deformation
relative to the size of mesh elements. It stays constant as the mesh is being refined
if the mesh deformation is in the order of the minimum element diameter, which is typical
in the MM solution of conservation laws with an explicit scheme.
On the other hand, $N_{\varsigma}$ will increase as the mesh is being refined
if the magnitude of the mesh deformation stays constant, a common situation
as in the MM solution of partial differential equations with a fixed time step or with an implicit scheme.
Moreover, the larger the mesh deformation is, the more time steps are needed.
Numerical examples in one and two dimensions have been presented to verify the convergence order,
mass conservation, positivity preservation, and cost analysis of the scheme.

As an application example, we have considered the use of the DG-interpolation scheme
in the rezoning MM-DG solution of RTE.
RTE has been discretized in our computation in angular directions using the discrete ordinate method,
in space using the DG method, and in time using the backward Euler scheme.
At each time step, the new mesh is generated using the MMPDE moving mesh method
and then the radiative intensity is interpolated from the old mesh to the new one using
the DG-interpolation scheme. Numerical results obtained for examples in one and two spatial dimensions
with various settings have demonstrated that the resulting rezoning MM-DG method
is 2nd-order with $P^1$-DG and 3rd-order with $P^2$-DG, more efficient than the method with
a fixed mesh, and able to preserve the positivity of the radiative intensity when the PP limiter is used.
It is also shown that the scheme is comparable in efficiency for not very fine meshes
with a quasi-Lagrangian MM-DG method developed in \cite{ZhangChengHuangQiu} for RTE
when a fixed time step size is used. It is still unclear if the latter can be made to preserve the positivity of the radiative positivity.

\section*{Acknowledgments}
M. Zhang and J. Qiu were supported partly by Science Challenge Project (China), No. TZ 2016002 and
National Natural Science Foundation--Joint Fund (China) grant U1630247.
This work was carried out while M. Zhang was visiting the Department of Mathematics, the University of Kansas
under the support by the China Scholarship Council (CSC: 201806310065).

% References

\end{document}